\documentclass[11pt]{article}
\usepackage{amsmath}
\usepackage{amssymb}
\usepackage[german,english]{babel}
\usepackage{url}
\usepackage{graphicx}
\usepackage{gastex}
\usepackage{longtable}
\usepackage{lscape}
\usepackage{tabularx}
\usepackage{multicol}
\usepackage{latexsym}
\usepackage{verbatim}
\usepackage{multirow}

\usepackage{graphicx}
\usepackage{epsfig,graphics,graphicx}
\usepackage{geometry}
\usepackage{color}
\hbadness=\maxdimen

\newtheorem{lem}{{\sc Lemma}}[section]
\newtheorem{thm}{{\sc Theorem}}[section]
\newtheorem{prop}{{\sc Proposition}}[section]

\newtheorem{claim}{{\sc Claim}}

\begin{document}
\title{\sc semi equivelar maps on the surface of euler characteristic $-1$}
\author{{Ashish Kumar Upadhyay and Anand Kumar Tiwari}\\[2mm]
{\normalsize Department of Mathematics}\\{\normalsize Indian Institute of Technology Patna }\\{\normalsize Patliputra Colony, Patna 800\,013,  India.}\\
{\small \{upadhyay, anand\}@iitp.ac.in}}
\maketitle

\vspace{-5mm}

\hrule

\begin{abstract} Semi-Equivelar maps are generalizations of Archimedean solids to the surfaces other than 2-sphere. In earlier work a complete classification of semi-equivelar map of type $(3^5, 4)$ on the surface of Euler characteristic $-1$ was given. In the meantime Karabas an Nedela classified vertex transitive semi-equivelar maps on the double torus. In this article we study the types of semi-equivelar maps on double torus that are also available on the surface of Euler characteristic $-1$. We classify them and show that none of them are vertex transitive. 
\end{abstract}

{\small

{\bf AMS classification\,: 52B70, 52C20}

{\bf Keywords\,:} Semi-Equivelar Maps
}

\bigskip

\hrule

\section{Introduction}\label{sec:introduction}
A triangulation of a surface is called $d$-covered if each edge of the triangulation is incident with a vertex of degree $d$. We got interested in studying the content presented in this article while attempting to answer a question of Negami and Nakamoto \cite{negami(2001)} about existence of $d$-covered triangulations for closed surfaces. We had answered their question in affirmative \cite{frank(2010)} for the surfaces of Euler characteristic $-127\leq \chi \leq -2$ and further became interested in looking at what happens for surfaces with $\chi = -1$. It was here that due to curvature considerations of this surface we had to construct a map on this surface which we named as Semi-equivelar map \cite{upadhyay1}. Such maps have also been studied in various forms (see \cite{altbrehm(1986)}, \cite{coxeter(1980),conder(1995),jones(1978),grunbaum(1987)}). In the meantime we came to learn that Nedela and Karabas {\cite{karabas(2007)}}, ${\cite{karabas(2012)}}$ have worked along similar lines and classified all the vertex transitive Archimedean maps on orientable surfaces of Euler characteristics $-2$, $-4$ and $-6$ (see also \cite{karabas(web)}). In particular, they have shown that there are seventeen isomorphism classes of Archimedean maps on orientable surface of Euler characteristic $-2$, out of which exactly fourteen are semi-equivelar maps with eleven distinct face sequences of types\,: $(3^5, 4)$, $(3^4, 4^2)$, $(3^4, 8)$, $(3^2, 4, 3, 6)$, $(3, 4^4)$, $(3, 4, 8, 4)$, $(3, 6, 4, 6)$, $(4^3, 6)$, $(4, 6, 16)$, $(4, 8, 12)$, $(6^2, 8)$. An orientable closed surface of Euler characteristic $-2$ is double cover of non orientable closed surface of Euler characteristic $-1$. This motivated us to explore the existence of above eleven types of semi-equivelar maps on non orientable surface of Euler characteristic $-1$. In the article \cite{upadhyay1} we have classified the semi-equivelar map of the type
$(3^5, 4)$ on this surface. Here, we investigate remaining types of semi-equivelar maps on this surface. In next few paragraphs we describe the definitions and terminologies used in this article. These definitions and terminologies are given in \cite{datta(2005)} and we are giving them here for the sake of ready reference. A standard reference on the subject of polyhedral maps is the article \cite{brehm(1997)} of Brehm and Schulte. For graph theory related terminologies one may also refer to\cite{bondy_murthy} and for topological preliminaries and terminologies one may refer to \cite{spanier(1981)}.

Throughout this article the term graph will mean a finite simple graph. A cycle of length $m$ or a {\em m-Cycle}, usually denoted by $C_m$, is by definition a connected 2-regular graph with $m$ vertices. A 2-dimensional {\em Polyhedral Complex} $K$ is a finite collection of $m_{i}$-cycles, where $\{m_{i}\colon 1 \leq i \leq n$ and $m_{i} \geq 3\} \subseteq \mathbb{N}$, together with vertices and edges of the cycles such that the non-empty intersection of any two cycles is either a vertex or is an edge. The cycles are called faces of $K$. The notations $V(K)$ and $EG(K)$ are used to denote the set of vertices and edges of $K$ respectively. A geometric object, called {\em geometric career} of $K$, denoted by $|K|$ can be associated to a polyhedral complex $K$ in the following manner\,: corresponding to each $m$-cycle $C_m$ in $K$, consider a $m$-gon $D_m$ whose boundary cycle is $C_m$. Then $|K|$ is union of all such $m$-gons. The complex $K$ is said to be connected (resp. compact or orientable) if $|K|$ is connected (resp. compact or orientable) topological space. A polyhedral complex $K$ is called a {\em Polyhedral 2-manifold} if for each vertex $v$ the faces containing $v$ are of the form $C_{m_1}, \ldots, C_{m_p}$ where $C_{m_1}\cap C_{m_2}, \ldots C_{m_{p-1}}\cap C_{m_p}$, and $C_{m_p}\cap C_{m_1}$ are edges for some $p \geq 3$. A connected polyhedral 2-manifold is called a {\em Polyhedral Map}. We will also use the term {\em map} for a polyhedral map. Among any two complexes $K_1$ and $K_2$ we define an isomorphism to be a bijective map $f\colon V(K_1) \longrightarrow V(K_2)$ for which $f(\sigma)$ is a face in $K_2$ if and only if $\sigma$ is a face in $K_1$. If $K_1 = K_2$ then $f$ is said to be an automorphism of $K_1$. The set of all automorphisms of a polyhedral complex $K$ form a group under the operation of composition of maps. This group is called the automorphism group of $K$. If this group acts transitively on the set $V(K)$ then the complex is called a {\em vertex transitive} complex. Some vertex transitive maps of Euler characteristic 0 have been studied in $\cite{babai(1991)}$ and many others in $\cite{brehm(2008)}, \cite{dattau(2005)}, \cite{dattau(2006)}$, \cite{datta_nandini} and $\cite{lutz(thesis)}$.


The {\em face sequence} (see \cite{upadhyay1}) of a vertex $v$ in a map is a finite cyclically ordered sequence $(a^p, b^q, ...., m^r)$ of powers of positive integers $a, b, ..., m \geq 3$ and $p, q, ..., r \geq 1$, such that through the vertex $v$, $p$ numbers of $C_a$, $q$ numbers of $C_b$, $\ldots$, $r$ numbers of $C_m$ are incident in the given cyclic order. A map $K$ is said to be {\em Semi-Equivelar} if face sequence of each vertex of $K$ is same. A SEM with face sequence $(a^p, b^q, ...., m^r)$, is also called SEM of type $(a^p, b^q, ...., m^r)$. In \cite{upadhyay0}, maps with face sequence $(3^3, 4^2)$ and $(3^2, 4, 3, 4)$ have been considered.

Let $EG(K)$ be the edge graph of a map $K$ and $V(K) = \{v_1, v_2, \ldots, v_n\}$. Let $L_{K}(v_i) = \{v_j \in V(K) \colon v_iv_j \in EG(K)\}$. For $0 \leq t \leq n$ define a graph $G_{t}(K)$ with $V(G_t(K)) = V(K)$ and $v_iv_j \in EG(G_t(K))$ if $|L_{K}(v_i) \bigcap L_{K}(v_j)| = t$. In other words the number of elements in the set $L_{K}(v_i) \bigcap L_{K}(v_j)$ is $t$. This graph was introduced in \cite{dattau(2006)} by B. Datta. Moreover if $K$ and $K'$ are two isomorphic maps then $G_i(K) \cong G_i(K')$ for each $i$. We have used these graphs in this article to determine whether two SEMs are isomorphic?

In the article \cite{upadhyay1} it has been shown that\,:

\begin{prop}\label{prop1} There exactly three non isomorphic semi equivelar maps of type $(3^5, 4)$ on the surface of Euler characteristic $-1$.
\end{prop}\hfill$\Box$.

In the present article we show\,:

\begin{lem}\label{lem1} If $K$ is a semi-equivelar map of type (3, 4, 8, 4) on the surface of Euler characteristic $-1$ then $K$ is isomorphic to $K_1(3, 4, 8, 4)$ or $K_2(3, 4, 8, 4)$ given in example Section \ref{example 4.1}.
\end{lem}

\begin{lem}\label{lem2}
If $M$ is a semi-equivelar map of type $(4, 6, 16)$ on the surface of Euler characteristic $-1$, then $M$ is isomorphic to $M_1(4, 6, 16)$ or $M_2(4, 6, 16)$ given in example Section \ref{example 4.1}.
\end{lem}

\begin{lem}\label{lem3}
 If $N$ is a semi-equivelar map of type $(6^2, 8)$ on the surface of Euler characteristic $-1$, then $N$ is isomorphic to $N_1(6^2, 8)$ or $N_2(6^2, 8)$ given in example Section \ref{example 4.1}.
\end{lem}

Thus from the above Proposition \ref{prop1} and Lemma \ref{lem1}, \ref{lem2}, \ref{lem3} it follows that\,:
\begin{thm}\label{thm1}
There are at least nine non-isomorphic semi-equivelar maps on the surface of Euler characteristic $-1$.
\end{thm}\hfill$\Box$

In the article we also show that\,:
\begin{thm}\label{thm2} There exist no semi-equivelar maps of types $(3^4, 8)$, $(3^2, 4, 3, 6)$, $(3, 6, 4, 6)$, $(4^3, 6)$ and $(4, 8, 12)$ on the surface of Euler characteristics -1.
\end{thm}

The article is organized in the following manner. In next section, we present examples of semi-equivelar maps on the surface of Euler characteristic $-1$. In the section, we describe the results and their proofs. We conclude the article by presenting a tabulated list of semi-equivelar maps on the surface of Euler characteristic $-1$.

\newpage

\section{Examples: Semi-equivelar maps on the surface of Euler characteristic $-1$}
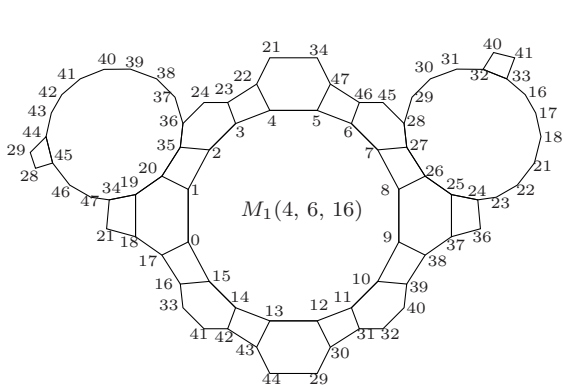
\label{example 4.1}

\begin{picture}(10, 10)(-40, 20)
\setlength{\unitlength}{3.5mm}
\drawpolygon(4, 1)(3.2, 2.5)(2.2, 3.5)(.9, 4)(-.9, 4)(-2.2, 3.5)(-3.2, 2.5)(-4, 1)(-4, -1)(-3.2, -2.5)(-2.2, -3.5)(-.9, -4)(.9, -4)(2.2, -3.5)(3.2, -2.5)(4, -1)

\drawpolygon(4, 1)(5, 1.5)(6, .8)(6, -.8)(5, -1.5)(4, -1)
\drawpolygon(-4, 1)(-5, 1.5)(-6, .8)(-6, -.8)(-5, -1.5)(-4, -1)

\drawpolygon(.9, 4)(1.4, 5)(.9, 6)(-.9, 6)(-1.4, 5)(-.9, 4)
\drawpolygon(.9, -4)(1.4, -5)(.9, -6)(-.9, -6)(-1.4, -5)(-.9, -4)

\drawpolygon(3.2, 2.5)(2.2, 3.5)(2.5, 4.3)(3.4, 4.3)(4.2, 3.5)(4.3, 2.6)

\drawpolygon(-3.2, 2.5)(-2.2, 3.5)(-2.5, 4.3)(-3.4, 4.3)(-4.2, 3.5)(-4.3, 2.6)

\drawpolygon(-3.2, -2.5)(-2.2, -3.5)(-2.5, -4.3)(-3.4, -4.3)(-4.2, -3.5)(-4.3, -2.6)

\drawpolygon(3.2, -2.5)(2.2, -3.5)(2.5, -4.3)(3.4, -4.3)(4.2, -3.5)(4.3, -2.6)

\drawline[AHnb=0](1.4, 5)(2.5, 4.3)
\drawline[AHnb=0](-1.4, 5)(-2.5, 4.3)
\drawline[AHnb=0](-1.4, -5)(-2.5, -4.3)
\drawline[AHnb=0](1.4, -5)(2.5, -4.3)

\drawline[AHnb=0](4.3, 2.6)(5, 1.5)
\drawline[AHnb=0](-4.3, 2.6)(-5, 1.5)
\drawline[AHnb=0](4.3, -2.6)(5, -1.5)
\drawline[AHnb=0](-4.3, -2.6)(-5, -1.5)



\drawpolygon(6, -.8)(6, .8)(7, .6)(7.1, -.5)

\drawpolygon(-6, -.8)(-6, .8)(-7, .6)(-7.1, -.5)

\drawpolygon(7, .6)(6, .8)(5, 1.5)(4.3, 2.6)(4.2, 3.5)(4.5, 4.5)(5.2, 5.2)(6.2, 5.57)(7.2, 5.56)(8.1, 5.2)(8.8, 4.6)(9.2, 3.9)(9.4, 3)(9.15, 2)(8.5, 1.16)(7.7, .7)

\drawpolygon(-7, .6)(-6, .8)(-5, 1.5)(-4.3, 2.6)(-4.2, 3.5)(-4.5, 4.5)(-5.2, 5.2)(-6.2, 5.57)(-7.2, 5.56)(-8.1, 5.2)(-8.8, 4.6)(-9.2, 3.9)(-9.4, 3)(-9.15, 2)(-8.5, 1.16)(-7.7, .7)

\drawpolygon(7.2, 5.56)(8.1, 5.2)(8.4, 6)(7.6, 6.2)
\drawpolygon(-9.4, 3)(-9.15, 2)(-9.8, 1.8)(-10, 2.4)

\put(-10.9, 2.3){\scriptsize \tiny 29}
\put(-10.4, 1.3){\scriptsize \tiny 28}

\put(8.4, 6){\scriptsize \tiny 41}
\put(7.2, 6.35){\scriptsize \tiny 40}

\put(3.3, .8){\scriptsize \tiny 8}
\put(2.7, 2.2){\scriptsize \tiny 7}

\put(1.9, 3){\scriptsize \tiny 6}
\put(.7, 3.4){\scriptsize \tiny 5}
\put(-1, 3.4){\scriptsize \tiny 4}
\put(-2.2, 3){\scriptsize \tiny 3}
\put(-3.1, 2.2){\scriptsize \tiny 2}
\put(-3.9, .8){\scriptsize \tiny 1}
\put(-3.9, -1.1){\scriptsize \tiny 0}
\put(-3.1, -2.5){\scriptsize \tiny 15}
\put(-2.4, -3.3){\scriptsize \tiny 14}
\put(-1.1, -3.8){\scriptsize \tiny 13}
\put(.6, -3.6){\scriptsize \tiny 12}
\put(1.5, -3.3){\scriptsize \tiny 11}
\put(2.2, -2.6){\scriptsize \tiny 10}
\put(3.4, -.99){\scriptsize \tiny 9}

\put(4.4, 2.5){\scriptsize \tiny 27}
\put(4.3, 3.4){\scriptsize \tiny 28}
\put(5, 1.5){\scriptsize \tiny 26}
\put(5.8, 1){\scriptsize \tiny 25}

\put(-5.2, 2.5){\scriptsize \tiny 35}
\put(-5.1, 3.4){\scriptsize \tiny 36}
\put(-5.9, 1.6){\scriptsize \tiny 20}
\put(-6.6, 1.1){\scriptsize \tiny 19}

\put(6.6, .75){\scriptsize \tiny 24}
\put(7.5, .3){\scriptsize \tiny 23}
\put(8.45, .8){\scriptsize \tiny 22}
\put(9.1, 1.7){\scriptsize \tiny 21}
\put(9.5, 2.8){\scriptsize \tiny 18}
\put(9.3, 3.7){\scriptsize \tiny 17}
\put(8.9, 4.5){\scriptsize \tiny 16}
\put(8.3, 5.1){\scriptsize \tiny 33}
\put(6.7, 5.15){\scriptsize \tiny 32}
\put(5.6, 5.75){\scriptsize \tiny 31}
\put(4.6, 5.45){\scriptsize \tiny 30}
\put(4.6, 4.3){\scriptsize \tiny 29}

\put(-7.28, .85){\scriptsize \tiny 34}
\put(-8.1, .35){\scriptsize \tiny 47}
\put(-9.2, .7){\scriptsize \tiny 46}
\put(-9.1, 2.1){\scriptsize \tiny 45}
\put(-10.3, 3){\scriptsize \tiny 44}
\put(-10.1, 3.8){\scriptsize \tiny 43}
\put(-9.7, 4.5){\scriptsize \tiny 42}
\put(-8.97, 5.2){\scriptsize \tiny 41}
\put(-7.45, 5.75){\scriptsize \tiny 40}
\put(-6.4, 5.65){\scriptsize \tiny 39}
\put(-5.2, 5.25){\scriptsize \tiny 38}
\put(-5.35, 4.3){\scriptsize \tiny 37}

\put(1.4, 5){\scriptsize \tiny 47}

\put(.6, 6.1){\scriptsize \tiny 34}

\put(-1.2, 6.2){\scriptsize \tiny 21}
\put(-2.3, 5.1){\scriptsize \tiny 22}
\put(-3, 4.6){\scriptsize \tiny 23}
\put(-3.9, 4.5){\scriptsize \tiny 24}


\put(1.4, -5.4){\scriptsize \tiny 30}

\put(.6, -6.4){\scriptsize \tiny 29}

\put(-1.2, -6.4){\scriptsize \tiny 44}
\put(-2.2, -5.3){\scriptsize \tiny 43}
\put(-3, -4.8){\scriptsize \tiny 42}
\put(-5.1, -3.6){\scriptsize \tiny 33}
\put(-3.95, -4.7){\scriptsize \tiny 41}

\put(2.4, -4.7){\scriptsize \tiny 31}

\put(3.3, -4.7){\scriptsize \tiny 32}
\put(4.3, -3.7){\scriptsize \tiny 40}

\put(4.4, -2.9){\scriptsize \tiny 39}

\put(-5.3, -2.8){\scriptsize \tiny 16}
\put(-5.9, -1.9){\scriptsize \tiny 17}
\put(5.1, -1.9){\scriptsize \tiny 38}

\put(5.8, -1.3){\scriptsize \tiny 37}
\put(-6.6, -1.2){\scriptsize \tiny 18}

\put(6.8, -1){\scriptsize \tiny 36}
\put(-7.5, -.98){\scriptsize \tiny 21}

\put(3.2, 4.4){\scriptsize \tiny 45}
\put(2.3, 4.45){\scriptsize \tiny 46}

\put(-2, 0){\scriptsize $M_1$(4, 6, 16)}

\end{picture}

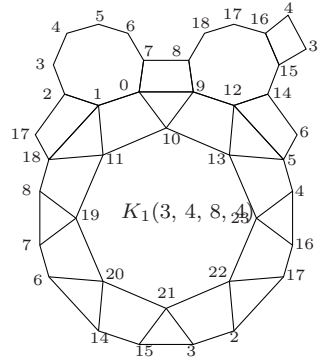
\begin{picture}(10, 10)(-105, 10)
\setlength{\unitlength}{3mm}
\drawpolygon(4, 0)(2.8, 2.8)(0, 4)(-2.8, 2.8)(-4, 0)(-2.8, -2.8)(0, -4)(2.8, -2.8)

\drawpolygon(5.6, 1.2)(5.2, 2.6)(3, 5)(1.2, 5.6)(-1.2, 5.6)(-3, 5)(-5.2, 2.6)(-5.6, 1.2)(-5.6, -1.2)(-5.2, -2.6)(-3, -5)(-1.2, -5.6)(1.2, -5.6)(3, -5)(5.2, -2.6)(5.6, -1.2)

\drawline[AHnb=0](0, 4)(1.2, 5.6)
\drawline[AHnb=0](0, 4)(-1.2, 5.6)
\drawline[AHnb=0](2.8, 2.8)(3, 5)
\drawline[AHnb=0](2.8, 2.8)(5.2, 2.6)

\drawline[AHnb=0](-2.8, 2.8)(-3, 5)
\drawline[AHnb=0](-2.8, 2.8)(-5.2, 2.6)

\drawline[AHnb=0](-2.8, -2.8)(-3, -5)
\drawline[AHnb=0](-2.8, -2.8)(-5.2, -2.6)

\drawline[AHnb=0](2.8, -2.8)(3, -5)
\drawline[AHnb=0](2.8, -2.8)(5.2, -2.6)

\drawline[AHnb=0](4, 0)(5.6, 1.2)
\drawline[AHnb=0](4, 0)(5.6, -1.2)

\drawline[AHnb=0](-4, 0)(-5.6, 1.2)
\drawline[AHnb=0](-4, 0)(-5.6, -1.2)

\drawline[AHnb=0](0, -4)(1.2, -5.6)
\drawline[AHnb=0](0, -4)(-1.2, -5.6)

\drawpolygon(1.2, 5.6)(1, 7)(-1, 7)(-1.2, 5.6)
\drawpolygon(1.2, 5.6)(1, 7)(1.7, 8.2)(3, 8.6)(4.4, 8.2)(5, 6.8)(4.5, 5.5)(3, 5)

\drawpolygon(-1.2, 5.6)(-1, 7)(-1.7, 8.2)(-3, 8.6)(-4.4, 8.2)(-5, 6.8)(-4.5, 5.5)(-3, 5)

\drawpolygon(5.2, 2.6)(3, 5)(4.5, 5.5)(5.8, 3.7)

\drawpolygon(-5.2, 2.6)(-3, 5)(-4.5, 5.5)(-5.8, 3.7)

\drawpolygon(5, 6.8)(4.4, 8.2)(5.4, 9)(6.2, 7.5)

\put(-.2, 3.3){\scriptsize \tiny 10}
\put(1.8, 2.4){\scriptsize \tiny 13}
\put(2.85, -.2){\scriptsize \tiny 23}
\put(1.85, -2.6){\scriptsize \tiny 22}
\put(-.4, -3.55){\scriptsize \tiny 21}
\put(-2.7, -2.7){\scriptsize \tiny 20}
\put(-3.8, -.1){\scriptsize \tiny 19}
\put(-2.7, 2.4){\scriptsize \tiny 11}
\put(5.35, 2.35){\scriptsize \tiny 5}

\put(5.7, 1){\scriptsize \tiny 4}
\put(5.7, -1.35){\scriptsize \tiny 16}
\put(5.3, -2.75){\scriptsize \tiny 17}
\put(2.8, -5.5){\scriptsize \tiny 2}
\put(.9, -6.1){\scriptsize \tiny 3}
\put(-1.4, -6.1){\scriptsize \tiny 15}
\put(-3.4, -5.55){\scriptsize \tiny 14}
\put(-5.9, -2.95){\scriptsize \tiny 6}
\put(-6.4, -1.4){\scriptsize \tiny 7}
\put(-6.4, 1.05){\scriptsize \tiny 8}
\put(-6.4, 2.5){\scriptsize \tiny 18}
\put(-6.9, 3.5){\scriptsize \tiny 17}
\put(-5.5, 5.3){\scriptsize \tiny 2}
\put(-5.8, 6.6){\scriptsize \tiny 3}
\put(-5.1, 8.25){\scriptsize \tiny 4}
\put(-3.15, 8.75){\scriptsize \tiny 5}
\put(-1.8, 8.3){\scriptsize \tiny 6}
\put(-1, 7.2){\scriptsize \tiny 7}
\put(.25, 7.2){\scriptsize \tiny 8}
\put(1.1, 8.45){\scriptsize \tiny 18}
\put(2.5, 8.85){\scriptsize \tiny 17}
\put(3.7, 8.6){\scriptsize \tiny 16}
\put(5.2, 9.1){\scriptsize \tiny 4}
\put(6.3, 7.4){\scriptsize \tiny 3}
\put(5, 6.3){\scriptsize \tiny 15}
\put(4.7, 5.35){\scriptsize \tiny 14}
\put(5.9, 3.6){\scriptsize \tiny 6}
\put(1.3, 5.65){\scriptsize \tiny 9}
\put(-2.1, 5.7){\scriptsize \tiny 0}
\put(-3.25, 5.3){\scriptsize \tiny 1}
\put(2.5, 5.45){\scriptsize \tiny 12}

\put(-2, 0){\scriptsize $K_1$(3, 4, 8, 4)}

\end{picture}


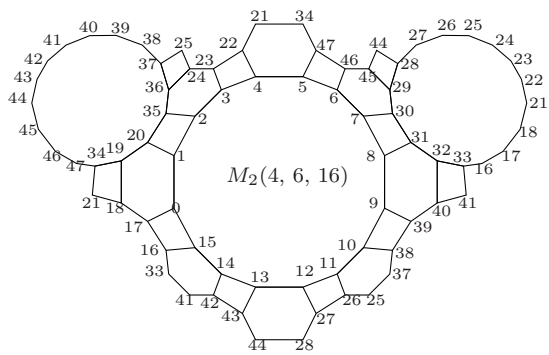
\begin{picture}(10, 10)(-30, 60)
\setlength{\unitlength}{3.5mm}
\drawpolygon(4, 1)(3.2, 2.5)(2.2, 3.5)(.9, 4)(-.9, 4)(-2.2, 3.5)(-3.2, 2.5)(-4, 1)(-4, -1)(-3.2, -2.5)(-2.2, -3.5)(-.9, -4)(.9, -4)(2.2, -3.5)(3.2, -2.5)(4, -1)

\drawpolygon(4, 1)(5, 1.5)(6, .8)(6, -.8)(5, -1.5)(4, -1)
\drawpolygon(-4, 1)(-5, 1.5)(-6, .8)(-6, -.8)(-5, -1.5)(-4, -1)

\drawpolygon(.9, 4)(1.4, 5)(.9, 6)(-.9, 6)(-1.4, 5)(-.9, 4)
\drawpolygon(.9, -4)(1.4, -5)(.9, -6)(-.9, -6)(-1.4, -5)(-.9, -4)

\drawpolygon(3.2, 2.5)(2.2, 3.5)(2.5, 4.3)(3.4, 4.3)(4.2, 3.5)(4.3, 2.6)

\drawpolygon(-3.2, 2.5)(-2.2, 3.5)(-2.5, 4.3)(-3.4, 4.3)(-4.2, 3.5)(-4.3, 2.6)

\drawpolygon(-3.2, -2.5)(-2.2, -3.5)(-2.5, -4.3)(-3.4, -4.3)(-4.2, -3.5)(-4.3, -2.6)

\drawpolygon(3.2, -2.5)(2.2, -3.5)(2.5, -4.3)(3.4, -4.3)(4.2, -3.5)(4.3, -2.6)

\drawline[AHnb=0](1.4, 5)(2.5, 4.3)
\drawline[AHnb=0](-1.4, 5)(-2.5, 4.3)
\drawline[AHnb=0](-1.4, -5)(-2.5, -4.3)
\drawline[AHnb=0](1.4, -5)(2.5, -4.3)

\drawline[AHnb=0](4.3, 2.6)(5, 1.5)
\drawline[AHnb=0](-4.3, 2.6)(-5, 1.5)
\drawline[AHnb=0](4.3, -2.6)(5, -1.5)
\drawline[AHnb=0](-4.3, -2.6)(-5, -1.5)

\drawpolygon(4.2, 3.5)(3.4, 4.3)(3.7, 5)(4.5, 4.5)

\drawpolygon(-4.2, 3.5)(-3.4, 4.3)(-3.7, 5)(-4.5, 4.5)

\drawpolygon(6, -.8)(6, .8)(7, .6)(7.1, -.5)

\drawpolygon(-6, -.8)(-6, .8)(-7, .6)(-7.1, -.5)

\drawpolygon(7, .6)(6, .8)(5, 1.5)(4.3, 2.6)(4.2, 3.5)(4.5, 4.5)(5.2, 5.2)(6.2, 5.57)(7.2, 5.56)(8.1, 5.2)(8.8, 4.6)(9.2, 3.9)(9.4, 3)(9.15, 2)(8.5, 1.16)(7.7, .7)

\drawpolygon(-7, .6)(-6, .8)(-5, 1.5)(-4.3, 2.6)(-4.2, 3.5)(-4.5, 4.5)(-5.2, 5.2)(-6.2, 5.57)(-7.2, 5.56)(-8.1, 5.2)(-8.8, 4.6)(-9.2, 3.9)(-9.4, 3)(-9.15, 2)(-8.5, 1.16)(-7.7, .7)

\put(3.3, .8){\scriptsize \tiny 8}
\put(2.7, 2.2){\scriptsize \tiny 7}

\put(1.9, 3){\scriptsize \tiny 6}
\put(.7, 3.4){\scriptsize \tiny 5}
\put(-1, 3.4){\scriptsize \tiny 4}
\put(-2.2, 3){\scriptsize \tiny 3}
\put(-3.1, 2.2){\scriptsize \tiny 2}
\put(-3.9, .8){\scriptsize \tiny 1}
\put(-4.1, -1.2){\scriptsize \tiny 0}
\put(-3.1, -2.5){\scriptsize \tiny 15}
\put(-2.4, -3.3){\scriptsize \tiny 14}
\put(-1.1, -3.8){\scriptsize \tiny 13}
\put(.6, -3.6){\scriptsize \tiny 12}
\put(1.5, -3.3){\scriptsize \tiny \tiny 11}
\put(2.2, -2.6){\scriptsize \tiny 10}
\put(3.4, -.99){\scriptsize \tiny 9}

\put(4.4, 2.5){\scriptsize \tiny 30}
\put(4.3, 3.4){\scriptsize \tiny 29}
\put(5, 1.5){\scriptsize \tiny 31}
\put(5.8, 1){\scriptsize \tiny 32}

\put(-5.2, 2.5){\scriptsize \tiny 35}
\put(-5.1, 3.4){\scriptsize \tiny 36}
\put(-5.8, 1.6){\scriptsize \tiny 20}
\put(-6.6, 1.1){\scriptsize \tiny 19}

\put(6.6, .75){\scriptsize \tiny 33}
\put(7.4, .3){\scriptsize \tiny 16}
\put(8.45, .8){\scriptsize \tiny 17}
\put(9.1, 1.7){\scriptsize \tiny 18}
\put(9.5, 2.8){\scriptsize \tiny 21}
\put(9.3, 3.7){\scriptsize \tiny 22}
\put(8.9, 4.5){\scriptsize \tiny 23}
\put(8.2, 5.2){\scriptsize \tiny 24}
\put(7, 5.75){\scriptsize \tiny 25}
\put(6, 5.75){\scriptsize \tiny 26}
\put(4.9, 5.4){\scriptsize \tiny 27}
\put(4.6, 4.3){\scriptsize \tiny 28}

\put(-7.3, .9){\scriptsize \tiny 34}
\put(-8.1, .4){\scriptsize \tiny 47}
\put(-9, .8){\scriptsize \tiny 46}
\put(-9.9, 1.8){\scriptsize \tiny 45}
\put(-10.3, 3){\scriptsize \tiny 44}
\put(-10.1, 3.8){\scriptsize \tiny 43}
\put(-9.7, 4.5){\scriptsize \tiny 42}
\put(-8.97, 5.2){\scriptsize \tiny 41}
\put(-7.6, 5.7){\scriptsize \tiny 40}
\put(-6.4, 5.65){\scriptsize \tiny 39}
\put(-5.2, 5.25){\scriptsize \tiny 38}
\put(-5.35, 4.3){\scriptsize \tiny 37}

\put(1.4, 5){\scriptsize \tiny 47}

\put(.6, 6.1){\scriptsize \tiny 34}

\put(-1.1, 6.1){\scriptsize \tiny 21}
\put(-2.3, 5.1){\scriptsize \tiny 22}
\put(-3.2, 4.45){\scriptsize \tiny 23}
\put(-3.5, 3.8){\scriptsize \tiny 24}
\put(-4, 5.1){\scriptsize \tiny 25}

\put(3.4, 5.1){\scriptsize \tiny 44}

\put(1.4, -5.4){\scriptsize \tiny 27}

\put(.6, -6.4){\scriptsize \tiny 28}

\put(-1.2, -6.4){\scriptsize \tiny 44}
\put(-2.2, -5.3){\scriptsize \tiny 43}
\put(-3, -4.8){\scriptsize \tiny 42}
\put(-5.1, -3.6){\scriptsize \tiny 33}
\put(-3.95, -4.7){\scriptsize \tiny 41}

\put(2.4, -4.7){\scriptsize \tiny 26}

\put(3.3, -4.7){\scriptsize \tiny 25}
\put(4.3, -3.7){\scriptsize \tiny 37}

\put(4.4, -2.9){\scriptsize \tiny 38}

\put(-5.3, -2.8){\scriptsize \tiny 16}
\put(-5.9, -1.9){\scriptsize \tiny 17}
\put(5.1, -1.9){\scriptsize \tiny 39}

\put(5.8, -1.3){\scriptsize \tiny 40}
\put(-6.6, -1.2){\scriptsize \tiny 18}

\put(6.8, -1){\scriptsize \tiny 41}
\put(-7.5, -.98){\scriptsize \tiny 21}

\put(3, 3.8){\scriptsize \tiny 45}
\put(2.3, 4.45){\scriptsize \tiny 46}

\put(-2, 0){\scriptsize $M_2$(4, 6, 16)}


\end{picture}

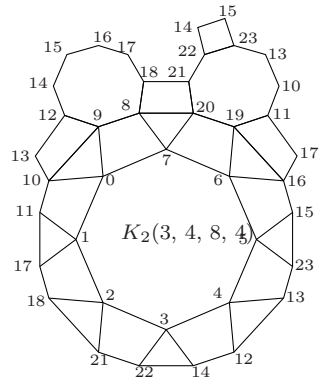
\begin{picture}(10, 10)(-105, 50)
\setlength{\unitlength}{3mm}
\drawpolygon(4, 0)(2.8, 2.8)(0, 4)(-2.8, 2.8)(-4, 0)(-2.8, -2.8)(0, -4)(2.8, -2.8)

\drawpolygon(5.6, 1.2)(5.2, 2.6)(3, 5)(1.2, 5.6)(-1.2, 5.6)(-3, 5)(-5.2, 2.6)(-5.6, 1.2)(-5.6, -1.2)(-5.2, -2.6)(-3, -5)(-1.2, -5.6)(1.2, -5.6)(3, -5)(5.2, -2.6)(5.6, -1.2)

\drawline[AHnb=0](0, 4)(1.2, 5.6)
\drawline[AHnb=0](0, 4)(-1.2, 5.6)
\drawline[AHnb=0](2.8, 2.8)(3, 5)
\drawline[AHnb=0](2.8, 2.8)(5.2, 2.6)

\drawline[AHnb=0](-2.8, 2.8)(-3, 5)
\drawline[AHnb=0](-2.8, 2.8)(-5.2, 2.6)

\drawline[AHnb=0](-2.8, -2.8)(-3, -5)
\drawline[AHnb=0](-2.8, -2.8)(-5.2, -2.6)

\drawline[AHnb=0](2.8, -2.8)(3, -5)
\drawline[AHnb=0](2.8, -2.8)(5.2, -2.6)

\drawline[AHnb=0](4, 0)(5.6, 1.2)
\drawline[AHnb=0](4, 0)(5.6, -1.2)

\drawline[AHnb=0](-4, 0)(-5.6, 1.2)
\drawline[AHnb=0](-4, 0)(-5.6, -1.2)

\drawline[AHnb=0](0, -4)(1.2, -5.6)
\drawline[AHnb=0](0, -4)(-1.2, -5.6)

\drawpolygon(1.2, 5.6)(1, 7)(-1, 7)(-1.2, 5.6)
\drawpolygon(1.2, 5.6)(1, 7)(1.7, 8.2)(3, 8.6)(4.4, 8.2)(5, 6.8)(4.5, 5.5)(3, 5)

\drawpolygon(-1.2, 5.6)(-1, 7)(-1.7, 8.2)(-3, 8.6)(-4.4, 8.2)(-5, 6.8)(-4.5, 5.5)(-3, 5)

\drawpolygon(1.7, 8.2)(3, 8.6)(2.6, 9.8)(1.4, 9.4)

\drawpolygon(5.2, 2.6)(3, 5)(4.5, 5.5)(5.8, 3.7)

\drawpolygon(-5.2, 2.6)(-3, 5)(-4.5, 5.5)(-5.8, 3.7)


\put(-.2, 3.3){\scriptsize \tiny 7}
\put(2.1, 2.4){\scriptsize \tiny 6}
\put(3.2, -.2){\scriptsize \tiny 5}
\put(2.05, -2.6){\scriptsize \tiny 4}
\put(-.3, -3.75){\scriptsize \tiny 3}
\put(-2.7, -2.7){\scriptsize \tiny 2}
\put(-3.8, -.1){\scriptsize \tiny 1}
\put(-2.7, 2.4){\scriptsize \tiny 0}
\put(5.35, 2.35){\scriptsize \tiny 16}

\put(5.7, 1){\scriptsize \tiny 15}
\put(5.7, -1.35){\scriptsize \tiny 23}
\put(5.3, -2.75){\scriptsize \tiny 13}
\put(3, -5.5){\scriptsize \tiny 12}
\put(1.1, -6.1){\scriptsize \tiny 14}
\put(-1.4, -6.1){\scriptsize \tiny 22}
\put(-3.4, -5.55){\scriptsize \tiny 21}
\put(-6.2, -2.95){\scriptsize \tiny 18}
\put(-6.8, -1.4){\scriptsize \tiny 17}
\put(-6.7, 1.05){\scriptsize \tiny 11}
\put(-6.4, 2.4){\scriptsize \tiny 10}
\put(-6.9, 3.5){\scriptsize \tiny 13}
\put(-5.7, 5.3){\scriptsize \tiny 12}
\put(-6.1, 6.7){\scriptsize \tiny 14}
\put(-5.45, 8.3){\scriptsize \tiny 15}
\put(-3.25, 8.75){\scriptsize \tiny 16}
\put(-2.2, 8.3){\scriptsize \tiny 17}
\put(-1.2, 7.2){\scriptsize \tiny 18}
\put(.05, 7.2){\scriptsize \tiny 21}
\put(.5, 8.15){\scriptsize \tiny 22}
\put(3.2, 8.7){\scriptsize \tiny 23}
\put(4.5, 8){\scriptsize \tiny 13}
\put(.3, 9.2){\scriptsize \tiny 14}
\put(2.3, 9.9){\scriptsize \tiny 15}
\put(5.1, 6.6){\scriptsize \tiny 10}
\put(4.7, 5.35){\scriptsize \tiny 11}
\put(5.9, 3.6){\scriptsize \tiny 17}
\put(1.3, 5.65){\scriptsize \tiny 20}
\put(-2, 5.75){\scriptsize \tiny 8}
\put(-3.25, 5.3){\scriptsize \tiny 9}
\put(2.6, 5.3){\scriptsize \tiny 19}

\put(-2, 0){\scriptsize $K_2$(3, 4, 8, 4)}


\end{picture}


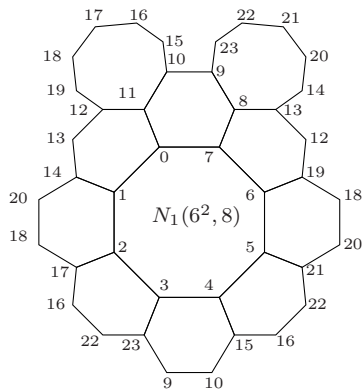
\begin{center}
\begin{picture}(75, 95)(5, 20)

\unitlength=5mm

\drawpolygon(2, .8)(.8, 2)(-.8, 2)(-2, .8)(-2, -.8)(-.8, -2)(.8, -2)(2, -.8)
\drawpolygon(2, .8)(3, 1.2)(4, .6)(4, -.6)(3, -1.2)(2, -.8)
\drawpolygon(-2, .8)(-3, 1.2)(-4, .6)(-4, -.6)(-3, -1.2)(-2, -.8)

\drawpolygon(.8, 2)(1.2, 3)(.6, 4)(-.6, 4)(-1.2, 3)(-.8, 2)
\drawpolygon(.8, -2)(1.2, -3)(.6, -4)(-.6, -4)(-1.2, -3)(-.8, -2)

\drawpolygon(1.2, 3)(.8, 2)(2, .8)(3, 1.2)(3.1, 2.2)(2.3, 3)
\drawpolygon(-1.2, 3)(-.8, 2)(-2, .8)(-3, 1.2)(-3.1, 2.2)(-2.3, 3)
\drawpolygon(1.2, -3)(.8, -2)(2, -.8)(3, -1.2)(3.1, -2.2)(2.3, -3)
\drawpolygon(-1.2, -3)(-.8, -2)(-2, -.8)(-3, -1.2)(-3.1, -2.2)(-2.3, -3)

\drawpolygon(.6, 4)(1.2, 3)(2.3, 3)(3, 3.5)(3.1, 4.4)(2.5, 5.2)(1.4, 5.3)(.7, 4.8)
\drawpolygon(-.6, 4)(-1.2, 3)(-2.3, 3)(-3, 3.5)(-3.1, 4.4)(-2.5, 5.2)(-1.4, 5.3)(-.7, 4.8)

\put(1.3, 3.1){\scriptsize \tiny 8}
\put(-1.9, 3.3){\scriptsize \tiny 11}
\put(-.8, 1.6){\scriptsize \tiny 0}
\put(-1.9, .6){\scriptsize \tiny 1}
\put(-1.9, -.8){\scriptsize \tiny 2}
\put(-.8, -1.8){\scriptsize \tiny 3}
\put(.4, -1.8){\scriptsize \tiny 4}
\put(1.5, -.8){\scriptsize \tiny 5}
\put(1.5, .6){\scriptsize \tiny 6}
\put(.4, 1.6){\scriptsize \tiny 7}
\put(4.1, .6){\scriptsize \tiny 18}
\put(3.1, 1.2){\scriptsize \tiny 19}

\put(3.2, 2.1){\scriptsize \tiny 12}
\put(2.5, 2.85){\scriptsize \tiny 13}
\put(3.1, 3.4){\scriptsize \tiny 14}

\put(-3.8, 2.1){\scriptsize \tiny 13}
\put(-3.2, 2.85){\scriptsize \tiny 12}
\put(-3.8, 3.4){\scriptsize \tiny 19}

\put(-3.7, -1.45){\scriptsize \tiny 17}
\put(-3.8, -2.3){\scriptsize \tiny 16}
\put(-2.9, -3.3){\scriptsize \tiny 22}
\put(-1.8, -3.4){\scriptsize \tiny 23}
\put(-4.8, .5){\scriptsize \tiny 20}
\put(-4.8, -.5){\scriptsize \tiny 18}

\put(-3.9, 1.2){\scriptsize \tiny 14}

\put(3.1, -1.4){\scriptsize \tiny 21}
\put(3.15, -2.3){\scriptsize \tiny 22}
\put(2.3, -3.3){\scriptsize \tiny 16}
\put(1.2, -3.4){\scriptsize \tiny 15}
\put(4.1, -.7){\scriptsize \tiny 20}

\put(-.7, -4.4){\scriptsize \tiny 9}
\put(.4, -4.4){\scriptsize \tiny 10}

\put(-.6, 4.2){\scriptsize \tiny 10}
\put(.7, 3.9){\scriptsize \tiny 9}

\put(-.65, 4.7){\scriptsize \tiny 15}
\put(.8, 4.5){\scriptsize \tiny 23}

\put(-1.6, 5.4){\scriptsize \tiny 16}
\put(1.25, 5.4){\scriptsize \tiny 22}

\put(-2.8, 5.4){\scriptsize \tiny 17}
\put(2.45, 5.3){\scriptsize \tiny 21}

\put(-3.9, 4.3){\scriptsize \tiny 18}
\put(3.2, 4.3){\scriptsize \tiny 20}

\put(-1, 0){\scriptsize $N_1(6^2, 8)$}
\end{picture}

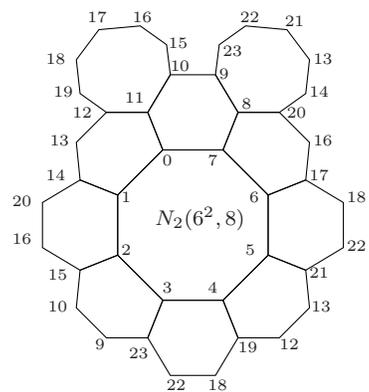
\begin{picture}(20, 10)(-50, 10)

\unitlength=5mm

\drawpolygon(2, .8)(.8, 2)(-.8, 2)(-2, .8)(-2, -.8)(-.8, -2)(.8, -2)(2, -.8)
\drawpolygon(2, .8)(3, 1.2)(4, .6)(4, -.6)(3, -1.2)(2, -.8)
\drawpolygon(-2, .8)(-3, 1.2)(-4, .6)(-4, -.6)(-3, -1.2)(-2, -.8)

\drawpolygon(.8, 2)(1.2, 3)(.6, 4)(-.6, 4)(-1.2, 3)(-.8, 2)
\drawpolygon(.8, -2)(1.2, -3)(.6, -4)(-.6, -4)(-1.2, -3)(-.8, -2)

\drawpolygon(1.2, 3)(.8, 2)(2, .8)(3, 1.2)(3.1, 2.2)(2.3, 3)
\drawpolygon(-1.2, 3)(-.8, 2)(-2, .8)(-3, 1.2)(-3.1, 2.2)(-2.3, 3)
\drawpolygon(1.2, -3)(.8, -2)(2, -.8)(3, -1.2)(3.1, -2.2)(2.3, -3)
\drawpolygon(-1.2, -3)(-.8, -2)(-2, -.8)(-3, -1.2)(-3.1, -2.2)(-2.3, -3)

\drawpolygon(.6, 4)(1.2, 3)(2.3, 3)(3, 3.5)(3.1, 4.4)(2.5, 5.2)(1.4, 5.3)(.7, 4.8)
\drawpolygon(-.6, 4)(-1.2, 3)(-2.3, 3)(-3, 3.5)(-3.1, 4.4)(-2.5, 5.2)(-1.4, 5.3)(-.7, 4.8)

\put(1.3, 3.1){\scriptsize \tiny 8}
\put(-1.8, 3.25){\scriptsize \tiny 11}
\put(-.8, 1.6){\scriptsize \tiny 0}
\put(-1.9, .6){\scriptsize \tiny 1}
\put(-1.9, -.8){\scriptsize \tiny 2}
\put(-.8, -1.8){\scriptsize \tiny 3}
\put(.4, -1.8){\scriptsize \tiny 4}
\put(1.4, -.8){\scriptsize \tiny 5}
\put(1.5, .6){\scriptsize \tiny 6}
\put(.4, 1.6){\scriptsize \tiny 7}
\put(4.1, .6){\scriptsize \tiny 18}
\put(3.1, 1.2){\scriptsize \tiny 17}

\put(3.2, 2.1){\scriptsize \tiny 16}
\put(2.5, 2.85){\scriptsize \tiny 20}
\put(3.1, 3.4){\scriptsize \tiny 14}

\put(-3.8, 2.1){\scriptsize \tiny 13}
\put(-3.2, 2.85){\scriptsize \tiny 12}
\put(-3.7, 3.4){\scriptsize \tiny 19}

\put(-3.85, -1.45){\scriptsize \tiny 15}
\put(-3.85, -2.3){\scriptsize \tiny 10}
\put(-2.6, -3.3){\scriptsize \tiny 9}
\put(-1.7, -3.55){\scriptsize \tiny 23}
\put(-4.8, .5){\scriptsize \tiny 20}
\put(-4.8, -.5){\scriptsize \tiny 16}

\put(-3.9, 1.2){\scriptsize \tiny 14}

\put(3.1, -1.4){\scriptsize \tiny 21}
\put(3.15, -2.3){\scriptsize \tiny 13}
\put(2.3, -3.3){\scriptsize \tiny 12}
\put(1.2, -3.4){\scriptsize \tiny 19}
\put(4.1, -.7){\scriptsize \tiny 22}

\put(-.7, -4.4){\scriptsize \tiny 22}
\put(.4, -4.4){\scriptsize \tiny 18}

\put(-.6, 4.1){\scriptsize \tiny 10}
\put(.7, 3.9){\scriptsize \tiny 9}

\put(-.65, 4.7){\scriptsize \tiny 15}
\put(.8, 4.5){\scriptsize \tiny 23}

\put(-1.6, 5.4){\scriptsize \tiny 16}
\put(1.25, 5.4){\scriptsize \tiny 22}

\put(-2.8, 5.4){\scriptsize \tiny 17}
\put(2.45, 5.3){\scriptsize \tiny 21}

\put(-3.9, 4.3){\scriptsize \tiny 18}
\put(3.2, 4.3){\scriptsize \tiny 13}
\put(-1, 0){\scriptsize $N_2 (6^2, 8)$}


\end{picture}

\end{center}

\newpage

\begin{claim}\label{claim 4.1}
$N_1(6^2, 8) \ncong N_2(6^2, 8)$ and $K_1(3, 4, 8, 4) \ncong K_2(3, 4, 8, 4)$, also
$N_1(6^2, 8)$, $N_2(6^2$, $8)$, $K_1(3, 4, 8, 4)$ and $K_2(3, 4, 8, 4)$ are not vertex transitive.
\end{claim}

\noindent{\bf Proof\,:} Consider the graphs $EG(G_{12}(N_1(6^2, 8)))$ = $\{[0, 7]$, $[3, 4]$, $[8, 13]$,
$[11, 12]$, $[15, 16]$, $[22, 23]\}$, $EG(G_{12}(N_2(6^2, 8)))$ = $\{[4, 5]$, $[18, 19]$, $[21, 22]\}$,
$EG(G_2(K_1(3, 4, 8, 4)))$ = $C_{12}($1, 10, 12, 6, 19, 18, 2, 21, 14, 5, 23, 17$) \cup C_{6}($4, 13, 9, 7, 20, 15$)$ and
$EG(G_2(K_2(3, 4, 8, 4)))$ = $C_{21}($0, 8, 21, 3, 12, 10, 1, 18, 20, 6, 15, 22, 2, 17, 19, 7, 9, 13, 5, 16, 11$)$.
From these graphs and discussions in Chapter 1 (page 12) it is evident that $N_1(6^2, 8)$ $\ncong$ $N_2(6^2, 8)$ and $K_1(3, 4, 8, 4)\ncong K_2(3, 4, 8, 4)$. Also from
these graphs one can deduce that above four maps are not vertex transitive. This proves the claim. $\hfill\Box$

\begin{claim}\label{claim 4.2}
$M_1(4, 6, 16)$ $\ncong$ $M_2(4, 6, 16)$ and $M_1(4, 6, 16)$, $M_2(4, 6, 16)$ are not vertex transitive.

\end{claim}

\noindent{\bf Proof\,:} Let $A(EG(M_1))$ and $A(EG(M_2))$ denote the adjacency matrices associated to edge graphs of
$M_1(4, 6, 16)$ and $M_2(4, 6, 16)$, respectively. Let $P_1(x)$ and $P_2(x)$ denote the characteristic polynomials of
$A(EG(M_1))$ and $A(EG(M_2))$ respectively. If the map $M_1(4, 6, 16)$ and $M_2(4, 6, 16)$ are isomorphic then $P_1(x)$ = $P_2(x)$, (see \cite{datta_nandini}). We have (using Maple)\,:

 $P_1(x)$= $x^{48} -73 x^{46}$ + $2454x^{44} - 50419x^{42}$ + $708648x^{40} - 63x^{39}$ + $3326x^{37}$ + $55370675x^{36} - 78998x^{35}
- 325536254x^{34}$ + $1117272x^{33}$ + $1488079446x^{32} - 10498532x^{31} - 5328759647x^{30}$ + $69274014x^{29}$ + $15001009001$ $x^{28}
- 330979906x^{27} - 33214008513x^{26}$ + $1164748518x^{25}$ + $57733175145x^{24} - 3045404365x^{23} - 78484320585 x^{22}$ +
$5935770108x^{21}$ + $82965261974x^{20} - 8621690840x^{19} - 67636071362x^{18}$ + $9302657658x^{17}$ + $42014823892 x^{16} -
7407374240x^{15} - 19530592234x^{14}$ + $4302417304x^{13}$ + $6604154516x^{12} - 1787400560x^{11} - 1549106652 x^{10}$ +
$513857976 x^9$ + $230136488x^8 - 96466160x^7 - 17066976x^6$ + $10545344x^5 - 49440x^4 - 495936x^3$ + $67264x^2 - 1920x$;

 $P_2(x)$ = $x^{48} - 72x^{46}$ + $2388x^{44}- 48424x^{42}$ + $672018x^{40} - 28x^{39} - 6770448x^{38}$ +
$1464x^{37}$ + $51267848x^{36} - 34548x^{35} - 298108536x^{34}$ + $486936x^{33}$ + $1348802145x^{32} - 4573164x^{31}$ $-
 4785171566x^{30}$ + $30247956x^{29}$ + $13360329054 x^{28} - 145305100 x^{27} - 29376425928x^{26}$ + $515828328x^{25}$ +
$50783351168x^{24} - 1365657624x^{23} - 68773076142x^{22}$ + $2706801464x^{21}$ + $72559583454x^{20} - 4017232620x^{19} - 59173427088x^{18}$ +
$4451481228x^{17}$ + $36880710516x^{16} - 3658879076x^{15}$ $-$ $17277557628x^{14}$ + $2204369472x^{13}$ + $5931587385x^{12} -
953952300x^{11} - 1432856946x^{10}$ + $286671228x^9$ + $226687857x^8 - 56423208x^7 - 20151768x^6$ + $6499968x^5$ + $573840x^4 -
330368x^3$ + $26880x^2$.

Therefore $M_1(4, 6, 16) \ncong M_2(4, 6, 16)$. Also, we have $EG(G_{15}(M_1(4, 6, 16)))$ = $EG$ $(G_{15}(M_2(4, 6, 16)))$ = $C_8($0, 2, 4, 6, 8, 10, 12, 14$)\cup C_8($1, 3, 5, 7, 9, 11, 13, 15$)\cup C_8($16, 18, 22, 24, 26, 28, 30, 32$) \cup C_8($17, 21, 23, 25, 27, 29, 31, 33$) \cup C_8($19, 35, 37, 39, 41, 43, 45, 47$) \cup C_8($20, 36, 38, 40, 42, 44, 46, 34$)$. Let $\alpha \in Aut(M_1(4, 6, 16))$ such that $\alpha(1) = 2$
then $\alpha$ induces an automorphism on $EG(G_{15}(M_1(4, 6, 16)))$. So $\alpha \{3, 15\} = \{0, 4\}$. This implies $\alpha (3) = 0$ or 4. But from the links of 1 and 2 it is easy to see that $\alpha(3) \neq 0$. So we have $\alpha (3) = 4$, this implies $\alpha(13) = 6$ and $\alpha(35) = 20$. From $\alpha(35) \mapsto 20$ we get $\alpha(42) = 43$. Now considering ${\textrm{lk}}(42)$ and ${\textrm{lk}}(43)$ and the map $\alpha(42) \mapsto 43$, we see that $\alpha (13) = 14$, a contradiction. Thus there is no automorphism which maps $1$ to $2$. Hence $M_1(4, 6, 16)$ is not vertex transitive. Similarly for $M_2(4, 6, 16)$ we get no automorphism such that $\alpha (1) = 2$. This proves the Claim \ref{claim 4.2}. $\hfill\Box$

\section{Enumeration of SEMs on the surface of Euler characteristic $-1$}

Considering Euler characteristic equation, it is easy to see that semi-equivelar maps of types
$(3^4, 4^2)$ and $(3, 4^4)$ do not exist on the surface of Euler characteristic $-1$. As, in these cases number of vertices required to complete a link of a vertex are more than the number of vertices of the SEMs. From the study of remaining eight types: $(3^4, 8)$, $(3^2, 4, 3, 6)$, $(3, 4, 8, 4)$, $(3, 6, 4, 6)$, $(4^3, 6)$,
$(4, 6, 16)$, $(4, 8, 12)$ and $(6^2, 8)$, we show the following \,:

\begin{lem}\label{lem 4.1}
There exists no SEM of type $(3^4, 8)$ on the surface of Euler characteristic $-1$.
\end{lem}

\noindent{\bf Proof\,:} Let $M$ be a SEM of type $(3^4, 8)$ on the surface of Euler characteristic $-1$. The notation
${\textrm{lk}}(i) = C_{10}([i_1, i_2, i_3, i_4, i_5, i_6, i_7], i_8, i_9, i_{10})$ for the link of a vertex $i$ will mean that $[i, i_1, i_{10}]$, $[i, i_9, i_{10}]$, $[i, i_8, i_9]$, $[i, i_7, i_8]$ form triangular faces and $[i, i_1, i_2, i_3, i_4, i_5, i_6, i_7]$ forms an octagonal face. If $|V|$ denotes the number of vertices in $V(M)$, $E(M)$ denotes the number of edges, $T(M)$ denotes the number of triangular faces and
$O(M)$ denotes the number of octagonal faces in map $M$, respectively, then it is easy to see that $E(M) = \frac{5|V|}{2}$, $T(M) = \frac{4|V|}{3}$ and $O(M) = \frac{|V|}{8}$. By Euler's equation we get, $-1$ = $|V|-\frac{5|V|}{2}+(\frac{4|V|}{3}+\frac{|V|}{8})$, $i.e.$ $-1 = |V|(\frac{-1}{24})$.
From the equation we see, if the the map exists then $|V| = 24$. Let $V = V(M) = \{$0, 1, \ldots, 23$\}$. Now, we prove the lemma by exhaustive search for all $M$.

Assume without loss of generality that ${\textrm{lk}}(0) = C_{10}([1, 2, 3, 4, 5, 6, 7], 8, 9, 10)$. This implies
${\textrm{lk}}(7) = C_{10}([6, 5, 4, 3, 2, 1, 0], 8, a, b)$ for some $a, b \in V$. One can see that $(a, b)\in\{$(10, 9), (11, 12)$\}$. In the first case when $(a, b) = (10, 9)$ then considering ${\textrm{ lk}}(10)$ we see that 1 lies in two octagonal faces, which is not allowed. In second case when $(a, b) = (11, 12)$ then we get ${\textrm{lk}}(7) = C_{10}([0, 1, 2, 3, 4, 5, 6], 12, 11, 8)$, ${\textrm{lk}}(6) = C_{10}([7, 0, 1, 2, 3, 4, 5], 14, 13, 12)$, ${\textrm{lk}}(5) = C_{10}([6, 7, 0, 1, 2, 3, 4], 16, 15, 14)$, ${\textrm{lk}}(4) = C_{10}([5, 6, 7, 0, 1, 2, 3], 18, 17, 16)$, ${\textrm{lk}}(3) = C_{10}([4, 5, 6, 7, 0, 1, 2], 20, 19, 18)$, ${\textrm{lk}}(2) = C_{10}([3, 4, 5, 6, 7, 0, 1], 22, 21, 20)$ and ${\textrm{lk}}(1) = C_{10}([2$, $3, 4, 5, 6, 7, 0], 10, 23, 22)$. This implies ${\textrm{lk}}(8) = C_{10}([9, c, d, e, f, g, h], 11, 7, 0)$ or ${\textrm{lk}}(8) = C_{10}([c$, $d, e, f, g, h, 11], 7, 0, 9)$ for some $c, d, e, f, g, h \in V$. But these two are isomorphic by the map (0, 7)(1, 6)(2, 5)(3, 4)(9, 11)(10, 12)(13, 23)(14, 22)(15, 21)(16, 20)(17, 19). Therefore, it is enough to consider ${\textrm{lk}}(8) = C_{10}([9, c, d, e, f, g, h], 11, 7, 0)$. Then we obtain the partial picture of the map
$M$ as shown in Figure I. Let $V(O_i)$, for $i$ = 1, 2, 3, denote the vertex set of an octagonal face $O_i$ then
it is easy to see that $V(O_1) = \{$0, 1, 2, 3, 4, 5, 6, 7$\}$, $V(O_2) = \{$8, 9, 13, 14, 17, 18, 21, 22$\}$ and
$V(O_3) = \{$10, 11, 12, 15, 16, 19, 20, 23$\}$. In this case we see that $(h, g)\in\{$(17, 18), (21, 22)$\}$.

If $(h, g) = (17, 18)$ then ${\textrm{lk}}(8) = C_{10}([9, c, d, e, f, 18, 17], 11, 7, 0)$, ${\textrm{lk}}(17) = C_{10}([8, 9$, $c, d, e, f, 18], 4, 16, 11)$ and ${\textrm{lk}}(18) = C_{10}([17, 8, 9, c, d, e, f], 19, 3, 4)$, where $f\in\{$13, 21$\}$. If $f = 13$ then $e = 14$ and $(c, d)\in\{$(21, 22), (22, 21)$\}$. In case $(c, d) = (21, 22)$, successively considering ${\textrm{lk}}(18)$, ${\textrm{lk}}(13)$ and ${\textrm{lk}}(14)$ we get $\deg(22) > 5$. A contradiction. On the other hand when $(c, d) = (22, 21)$ then successively considering ${\textrm{lk}}(18)$,
${\textrm{lk}}(13)$, ${\textrm{lk}}(14)$, ${\textrm{lk}}(21)$, ${\textrm{lk}}(22)$, ${\textrm{lk}}(9)$, ${\textrm{lk}}(8)$ and ${\textrm{lk}}(17)$, we get $C_4(0, 1, 23, 9) \subseteq {\textrm{lk}}(10)$. Again, a contradiction. Also for $f = 21$, considering ${\textrm{lk}}(18)$ we see ${\textrm{lk}}(21)$ can not be completed. So $(h, g) \neq (17, 18)$.

If $(h, g) = (21, 22)$ then $f\in\{$13, 17$\}$. In the first case when $f = 13$ then we have $e = 14$ and $(c, d) = (18, 17)$, now considering ${\textrm{lk}}(14)$ and ${\textrm{lk}}(17)$ successively we get a
triangular face [15, 16, 17] in $M$. This is not possible. So $f \neq 13$. On the other hand when $f = 17$ then $e = 18$ and $(c, d) = (14, 13)$, now successively considering ${\textrm{lk}}(18)$, ${\textrm{lk}}(13)$, ${\textrm{lk}}(14)$, ${\textrm{lk}}(9)$, ${\textrm{lk}}(8)$, ${\textrm{lk}}(21)$, ${\textrm{lk}}(22)$ and ${\textrm{lk}}(17)$, one can see that ${\textrm{lk}}(10)$ can not be
completed. So $(h, g) \neq (21, 22)$ and thus the lemma is proved. $\hfill\Box$

\vskip1.5in

\begin{center}
\begin{picture}(0, 0)(0, 7)
\setlength{\unitlength}{8.5mm}
\drawpolygon(2, 1)(1, 2)(-1, 2)(-2, 1)(-2, -1)(-1, -2)(1, -2)(2, -1)

\drawpolygon(3, 0)(3, 1.2)(2.5, 2)(1.5, 3)(0, 3)(-1.5, 3)(-2.5, 2)(-3, 1.2)(-3, 0)(-3, -1.2)(-2.5, -2)
(-1.5, -3)(0, -3)(1.5, -3)(2.5, -2)(3, -1.2)

\drawline[AHnb=0](2, 1)(3, 0)
\drawline[AHnb=0](2, -1)(3, 0)
\drawline[AHnb=0](2, 1)(3, 1.2)
\drawline[AHnb=0](2, 1)(2.5, 2)

\drawline[AHnb=0](1, 2)(1.5, 3)
\drawline[AHnb=0](1, 2)(0, 3)

\drawline[AHnb=0](-1, 2)(0, 3)
\drawline[AHnb=0](-1, 2)(-1.5, 3)

\drawline[AHnb=0](-1, 2)(-2.5, 2)
\drawline[AHnb=0](-2, 1)(-2.5, 2)

\drawline[AHnb=0](-2, 1)(-3, 1.2)
\drawline[AHnb=0](-2, 1)(-3, 0)

\drawline[AHnb=0](-2, -1)(-3, -1.2)
\drawline[AHnb=0](-2, -1)(-2.5, -2)

\drawline[AHnb=0](-1, -2)(-1.5, -3)
\drawline[AHnb=0](-1, -2)(0, -3)

\drawline[AHnb=0](1, -2)(1.5, -3)
\drawline[AHnb=0](1, -2)(2.5, -2)

\drawline[AHnb=0](2, -1)(2.5, -2)
\drawline[AHnb=0](2, -1)(3, -1.2)

\drawline[AHnb=0](1, 2)(2.5, 2)
\drawline[AHnb=0](-2, -1)(-3, 0)

\drawline[AHnb=0](-1, -2)(-2.5, -2)
\drawline[AHnb=0](0, -3)(1, -2)

\drawpolygon(3, 0)(3.5, -.5)(4.2, -.4)(4.7, .1)(4.7, .9)(4, 1.5)(3.4, 1.5)(3, 1.2)

\drawline[AHnb=0](3, -1.2)(3.5, -.5)

\drawpolygon(1.5, 3)(2.5, 2)(3, 2)(3.5, 2.4)(3.7, 3)(3.5, 3.5)(2.8, 3.8)(2, 3.6)

\drawline[AHnb=0](3, 2)(3.4, 1.5)
\drawline[AHnb=0](3, 1.2)(3, 2)

\drawpolygon(0, 3)(-1.5, 3)(-2, 3.5)(-2, 4.3)(-1.3, 4.8)(-.3, 4.8)(.4, 4.4)(.7, 3.7)

\drawline[AHnb=0](.7, 3.7)(1.5, 3)
\drawline[AHnb=0](.7, 3.7)(2, 3.6)

\drawpolygon(-2.5, 2)(-2.7, 2.7)(-3.3, 3)(-4, 3)(-4.6, 2.3)(-4.5, 1.7)(-4, 1.2)(-3, 1.2)

\drawline[AHnb=0](-2.7, 2.7)(-2, 3.5)
\drawline[AHnb=0](-2.7, 2.7)(-1.5, 3)

\drawpolygon(-3, 0)(-3.6, .3)(-4.3, .2)(-4.8, -.3)(-4.8, -1.2)(-4.2, -1.7)(-3.5, -1.7)(-3, -1.2)

\drawline[AHnb=0](-4, 1.2)(-3.6, .3)
\drawline[AHnb=0](-3.6, .3)(-3, 1.2)

\drawpolygon(-2.5, -2)(-3, -2.2)(-3.5, -2.6)(-3.6, -3.2)(-3.3, -3.8)(-2.3, -4)(-1.7, -3.7)(-1.5, -3)

\drawline[AHnb=0](-3, -2.2)(-3, -1.2)
\drawline[AHnb=0](-3, -2.2)(-3.5, -1.7)

\drawpolygon(0, -3)(-.5, -3.8)(-.4, -4.3)(.2, -4.7)(1, -4.7)(1.5, -4.4)(1.8, -3.7)(1.5, -3)

\drawline[AHnb=0](-.5, -3.8)(-1.5, -3)
\drawline[AHnb=0](-.5, -3.8)(-1.7, -3.7)

\drawpolygon(2.5, -2)(3, -1.2)(3.7, -1.3)(4.2, -1.6)(4.4, -2.3)(4.1, -2.8)(3.3, -3.1)(2.7, -2.9)

\drawline[AHnb=0](3.7, -1.3)(3.5, -.5)
\drawline[AHnb=0](3, -1.2)(3.5, -.5)

\drawline[AHnb=0](1.5, -3)(2.7, -2.9)
\drawline[AHnb=0](2.7, -2.9)(1.8, -3.7)

\put(-.98, 1.6){\scriptsize 0}

\put(.87, 1.6){\scriptsize 7}
\put(1.7, .8){\scriptsize 6}
\put(1.7, -.8){\scriptsize 5}
\put(.8, -1.8){\scriptsize 4}
\put(-.8, -1.8){\scriptsize 3}
\put(-1.9, -.8){\scriptsize 2}
\put(-1.9, .8){\scriptsize 1}

\put(-1.5, 3.2){\scriptsize 9}
\put(-.1, 3.2){\scriptsize 8}
\put(1.67, 2.9){\scriptsize 11}
\put(2.4, 2.2){\scriptsize 12}
\put(3.1, 1){\scriptsize 13}
\put(3.1, 0){\scriptsize 14}
\put(2.9, -1.6){\scriptsize 15}
\put(2.6, -2.2){\scriptsize 16}
\put(1.1, -3.4){\scriptsize 17}
\put(-.1, -3.4){\scriptsize 18}
\put(-2.1, -3.2){\scriptsize 19}
\put(-2.8, -2.5){\scriptsize 20}
\put(-3.6, -1.27){\scriptsize 21}
\put(-3.6, -.2){\scriptsize 22}
\put(-3.4, 1.3){\scriptsize 23}
\put(-3.1, 1.9){\scriptsize 10}

\put(0, 0){\scriptsize $O_1$}

\put(-1, 4){\scriptsize $O_2$}

\put(-3.9, 2){\scriptsize $O_3$}
\put(-4.3, -.8){\scriptsize $O_2$}
\put(-3, -3.3){\scriptsize $O_3$}
\put(.3, -4.1){\scriptsize $O_2$}
\put(3.5, -2.4){\scriptsize $O_3$}
\put(4, .4){\scriptsize $O_2$}
\put(2.9, 3){\scriptsize $O_3$}

\put(.7, 3.79){\scriptsize $h$}
\put(.5, 4.4){\scriptsize $g$}
\put(-.3, 4.9){\scriptsize $f$}
\put(-1.4, 4.9){\scriptsize $e$}
\put(-2.3, 4.2){\scriptsize $d$}
\put(-2.35, 3.5){\scriptsize $c$}

\put(-5, -5.8){\mbox{Figure I: Semi-equivelar map M of type $(3^4, 8)$}}

\end{picture}
\end{center}

\vspace{5cm}

\begin{lem}\label{lem 4.2}
There exists no SEM of type $(3^2, 4, 3, 6)$ on the surface of Euler characteristic $-1$.

\end{lem}

\noindent{\bf Proof \,:} Let $G$ be a SEM of type $(3^2, 4, 3, 6)$ on the surface of Euler characteristic $-1$.
The notation ${\textrm{lk}}(i) = C_{11}([i_1, i_2, i_3, i_4, i_5], i_6, i_7, i_8, i_9)$ for the link of $i$ will mean that $[i, i_1, i_9]$, $[i, i_5, i_6]$, $[i, i_6, i_7]$ form triangular faces, $[i, i_7, i_8, i_9]$ forms quadrangular face and $[i, i_1, i_2, i_3, i_4$ $, i_5]$ forms hexagonal face. Let $|V|$ denote the number of vertices in $V(G)$. If $E(G)$, $T(G)$, $Q(G)$ and $H(G)$ denote the number of edges, number of triangular faces, number of quadrangular faces and number of hexagonal faces in the map $G$, respectively, then it is easy to see that $E(G) = \frac{5|V|}{2}$, $T(G) = \frac{3|V|}{3}$, $Q(G) = \frac{|V|}{4}$ and $H(G) = \frac{|V|}{6}$. By Euler's equation we see, if the map exists then $|V| = 12$. Let $V = V(G) = \{$0, 1, \ldots, 11$\}$. Now, we prove the lemma by exhaustive search for all $G$. Assume that ${\textrm{lk}}(0) = C_{11}([1, 2, 3, 4, 5], 6, 7, 8, 9)$ then ${\textrm{lk}}(7) = C_{11}([a, b, c, d, e], 6, 0, 9, 8)$ or ${\textrm{lk}}(7) = C_{11}([6, a, b, c, d], e, 8, 9, 0)$ for some $a, b, c, d, e \in V$. But, for both the cases we need more than twelve vertices to complete ${\textrm{lk}}(7)$. This is not allowed. So the map does not exist. $\hfill\Box$

\begin{lem}\label{lem 4.3}
There exists no SEM of type $(3, 6, 4, 6)$ on the surface of Euler characteristic $-1$.

\end{lem}

\noindent {\bf Proof\,:} Let $E$ be a SEM of type $(3, 6, 4, 6)$ on the surface of Euler characteristic $-1$.
The notation ${\textrm{lk}}(i) = C_{11}([i_1, i_2, i_3, i_4, i_5], [i_6, i_7, i_8, i_9, i_{10}], i_{11})$ for the link of $i$ will mean that $[i, i_5, i_6]$ forms triangular face, $[i, i_1, i_{11}, i_{10}]$ forms quadrangular face and $[i, i_1, i_2, i_3, i_4, i_5]$, $[i, i_6, i_7, i_8, i_9, i_{10}]$ form hexagonal faces. Let $|V|$ denote the number of vertices in $V(E)$. If $E(E)$, $T(E)$, $Q(E)$ and $H(E)$ denote the number of edges, number of triangular faces, number of quadrangular faces and number of hexagonal faces, respectively, then we see that $E(E) = \frac{4|V|}{2}$, $T(E) = \frac{|V|}{3}$, $Q(E) = \frac{|V|}{4}$ and $H(E) = \frac{2|V|}{6}$. By Euler's equation we
see, if the map exists then $|V| = 12$. For this, let $V = V(E) = \{$0, 1, \ldots, 11$\}$. Now, we prove the lemma by exhaustive search for all $E$. For this assume that ${\textrm{lk}}(0) = C_{11}([1, 2, 3, 4, 5], [6, 7, 8, 9, 10], 11)$, then ${\textrm{lk}}(1) = C_{11}([0, 5, 4, 3, 2], [a$, $b, c, d, 11], 10)$ for some $a, b, c, d \in V$. Now, it is easy to see that ${\textrm{lk}}(1)$ can not be completed, as $a$, $b$, $c$, $d$ have no suitable values in $V(E)$. Therefore the required map does not exist.
Thus the lemma is proved.$\hfill\Box$

\begin{lem}\label{lem 4.4}
There exists no SEM of type $(4^3, 6)$ on the surface of Euler characteristic $-1$.

\end{lem}

\noindent{\bf Proof\,:} Let $F$ be a SEM of type $(4^3, 6)$ on the surface of Euler characteristic $-1$.
The notation ${\textrm{lk}}(i) = C_{11}([i_1, i_2, i_3, i_4, i_5], i_6, i_7, i_8, i_9, i_{10})$ for the link of $i$ will mean that $[i, i_1, i_{10}, i_{9}]$, $[i, i_5, i_6, i_7]$, $[i, i_7, i_8, i_9]$ form quadrangular faces and $[i, i_1, i_2, i_3, i_4, i_5]$ forms hexagonal face. Let $|V|$ denote the number of vertices in $V(F)$. If $E(F)$, $Q(F)$ and $H(F)$ denote the number of edges, number of quadrangular faces and number of hexagonal faces, respectively, then $E(F) = \frac{4|V|}{2}$, $Q(F) = \frac{3|V|}{4}$ and $H(F) = \frac{|V|}{6}$. By Euler's equation we see if the map exists then $|V| = 12$. For this, let $V = V(F) = \{$0, 1, \ldots, 11$\}$. Now we prove the lemma by exhaustive search for all $F$. Assume that ${\textrm{lk}}(0) = C_{11}([1, 2, 3, 4, 5], 6, 7, 8, 9, 10)$. This implies,
${\textrm{lk}}(7) = C_{11}([b, c, d, e, 6], 5, 0, 9, 8, a)$ or ${\textrm{lk}}(7) = C_{11}([b, c, d, e, 8], 9, 0, 5, 6, a)$ for some $a, b, c, d, e \in V$. Then for both the cases of ${\textrm{lk}}(7)$ we need more than twelve vertices to complete. But this is not allowed. So we do not get the required map. Thus the lemma is proved. $\hfill\Box$

\begin{lem}\label{lem 4.5}
There exists no SEM of type $(4, 8, 12)$ on the surface of Euler characteristic $-1$.

\end{lem}

\noindent{\bf Proof\,:} Let $M$ be a SEM of type $(4, 8, 12)$ on the surface of Euler characteristic $-1$.
The notation ${\textrm{lk}}(i) = C_{18}([{\textbf{\textit{i}}}_{\textbf{1}}, i_2, i_3, i_4, i_5, i_6, i_7, i_8, i_9, i_{10}, {\textbf{\textit{i}}}_{\textbf{11}}], i_{12}, [{\textbf{\textit{i}}}_{\textbf{13}}, i_{14}, i_{15}, i_{16}, i_{17}, i_{18}])$ for the link of $i$ will mean that $[i, i_{11}, i_{12}, i_{13}]$, $[i, i_1, i_{18}, i_{17}, i_{16}, i_{15}, i_{14}, i_{13}]$ and $[i, i_1, i_2, i_3, i_4, i_5, i_6, i_7$, $i_8, i_9, i_{10}, i_{11}]$ form 4-gonal, 8-gonal and 12-gonal faces. If $|V|$, $E(M)$, $Q(M)$, $O(M)$ and $R(M)$ denote the number of vertices, number of edges, number of 4-gonal faces, number of 8-gonal faces and number of 12-gonal faces in $M$, respectively, then we see that $E(M) = \frac{3|V|}{2}$, $Q(M) = \frac{|V|}{4}$, $O(M) = \frac{|V|}{8}$ and $R(M) = \frac{|V|}{12}$. By Euler's equation we see, if the map exists then $|V| = 24$. For this, let $V = V(M) = \{$0, 1, \ldots, 23$\}$.
Now we proceed as follows. Assume that ${\textrm{lk}}(0) = C_{18}([{\textbf 1}, 2, 3, 4, 5, 6, 7, 8, 9, 10, {\textbf{11}}], 12, [{\textbf{13}}, 14, 15, 16, 17, 18])$. This implies ${\textrm{lk}}(1) = C_{18}([{\textbf 0}, 11, 10, 9, 8, 7, 6, 5, 4, 3, {\textbf 2}], 19, [{\textbf{18}}, 17, 16, 15, 14, 13])$ and ${\textrm{lk}}(2) = C_{18}$ $([{\textbf 3}, 4, 5, 6, 7, 8, 9, 10, 11, 0, {\textbf 1}], 18, [{\textbf{19}}, a, b, c, d, e])$ for some
$a, b, c, d, e \in V$. Observe that $a \in \{$12, 20$\}$. If $a = 12$ then $b = 13$, for otherwise $\deg(12) > 3$.
But, then 13 appears in two octagonal faces, which is not allowed. So we have $a = 20$, this implies $b \in \{$12, 21$\}$. If $b = 12$ then $c = 13$ and we get 13 in two octagonal faces. So $b = 21$, this implies $c \in \{$12, 22$\}$. In case $c = 12$, $d = 13$. This implies 13 appears in two octagonal faces. If $c = 22$ then $d = 23$, now we see that $e$ has no suitable value in $V$ so that ${\textrm{lk}}(2)$ can be completed. So, the required map does not exist.

\smallskip
\noindent{\bf Proof of Theorem \ref{thm2}\,:} The proof of Theorem \ref{thm2} follows from Lemmas \ref{lem 4.1}, \ref{lem 4.2}, \ref{lem 4.3}, \ref{lem 4.4} and \ref{lem 4.5}. \hfill $\Box$

\smallskip

\noindent {\bf Proof of Lemma \ref{lem1}\,:} Let $K$ be a SEM of type $(3, 4, 8, 4)$ on the surface of Euler characteristic $-1$. The notation ${\textrm{ lk}}(i) = C_{11}([i_1, i_2, i_3, i_4, i_5, i_6, i_7], i_8, [i_9, i_{10}], i_{11})$ for the link of $i$ will mean that $[i, i_9, i_{10}]$ forms triangular face, $[i, i_7, i_8, i_9]$, $[i, i_1, i_{11}, i_{10}]$ form quadrangular faces and $[i, i_1, i_2, i_3, i_4, i_5, i_6, i_7]$ forms octagonal face. Let $|V|$ denote the number of vertices in $V(K)$. If $E(K)$, $T(K)$, $Q(K)$ and $O(K)$ denote the number of edges, number of triangular faces, number of quadrangular faces and number of octagonal faces in the map $K$, respectively, then we see that $E(K) = \frac{4|V|}{2}$, $T(K) = \frac{|V|}{3}$, $Q(K) = \frac{2|V|}{4}$ and $H(K) = \frac{|V|}{8}$. By Euler's equation we see, if the map exists then $|V| = 24$. Let $V = V(K) = \{$0, 1, \ldots, 23$\}$. Now, we prove the result by exhaustive search for all $K$.

Let ${\textrm{lk}}(0) = C_{11}([1, 2, 3, 4, 5, 6, 7], 8, [9, 10], 11)$, this implies
${\textrm{lk}}(9) = C_{11}([b, c, d, e, f, g$, $8], 7, [0, 10], a)$ and ${\textrm{lk}}(10) = C_{11}([11, l, k, j, i, h, a], 12, [9, 0], 1)$ for some $a, b, c, d, e, f, g, h, i, j$, $k, l \in V$. Observe that $b = 12$ and $a = 13$, then octagonal faces of the map $K$ are, $O_1 = [0, 1, 2, 3, 4, 5, 6, 7]$, $O_2 = [8, 9, 12, c, d, e, f, g]$ and $O_3 = [13, 10, 11, l, k, j, i, h]$. As, these faces share no vertex with each other, successively we get
$c = 14$, $d = 15$, $e = 16$, $f = 17$, $g = 18$, $l = 19$, $k = 20$, $j = 21$, $i = 22$ and $h = 23$.
This implies ${\textrm{lk}}(9) = C_{11}([12, 14, 15, 16, 17, 18, 8], 7, [0, 10], 13)$,
${\textrm{lk}}(10) = C_{11}([11, 19, 20, 21, 22, 23, 13], 12, [9, 0], 1)$ and ${\textrm{lk}}(8) = C_{11}([18, 17, 16, 15, 14, 12, 9], 0, [7, x], y)$ for some $x, y \in V$. In this case, $(x, y) \in \{$(19, 11), (19, 20), (20, 19), (20, 21), (21, 20), (21, 22), (22, 21), (22, 23), (23, 13), (23, 22)$\}$. If $(x, y) = (23, 13)$ then considering ${\textrm{lk}}(8)$ and ${\textrm{lk}}(13)$ successively we see 12\,18 as an edge and a non-edge both. Also,
$(19, 20)\cong(23, 13)$; $(20, 19)\cong(22, 21)$ and $(20, 21)\cong(22, 23)$ by the map (0, 9)(1, 12)(2, 14)(3, 15)(4, 16)(5, 17) (6, 18)(7, 8)(11, 13)(19, 23)(20, 22); $(20, 19)\cong(21, 22)$ by the map (0, 8)(1, 18)(2, 17)(3, 16)(4, 15)(5, 14)(6, 12)(7, 9)(10, 21) (11, 22)(13, 20)(19, 23); $(19, 11)\cong(21, 20)$ by the map (0, 8)(1, 18)(2, 17)(3, 16)(4, 15)(5, 14)(6, 12)(7, 9)(10, 21)(11, 20)(13, 22). So, it is enough to search the map for $(x, y) \in \{$(19, 11), (20, 19), (20, 21), (23, 22)$\}$.

When $(x, y) = (20, 19)$ then ${\textrm{lk}}(8) = C_{11}([9, 12, 14, 15, 16, 17, 18], 19, [20, 7], 0)$, ${\textrm{lk}}(7) = C_{11}([0, 1, 2, 3, 4, 5, 6], 21, [20, 8], 9)$,
${\textrm{lk}}(20) = C_{11}([21, 22, 23, 13, 10, 11, 19], 18$, $[8, 7], 6)$ and
${\textrm{lk}}(21) = C_{11}([22, 23, 13, 10, 11, 19, 20], 7, [6, m], n)$ for some $m, n \in V$.
In this case $(m, n)\in\{$(14, 15), (15, 14), (15, 16), (16, 15), (16, 17), (17, 16)$\}$. But $(14, 15)\cong(16, 15)$ by the map (0, 6)(1, 5)(2, 4)(8, 20)(9, 21)(10, 16)(11, 17)(12, 22)(13, 15)(14, 23)(18, 19), so we consider the following subcases.

When $(m, n) = (15, 14)$ then successively considering ${\textrm{lk}}(21)$,
${\textrm{lk}}(15)$, ${\textrm{lk}}(6)$ and ${\textrm{lk}}(16)$ one can see that ${\textrm{lk}}(17)$ can not be completed. When $(m, n) = (15, 16)$ then completing ${\textrm{lk}}(21)$, ${\textrm{lk}}(15)$ and ${\textrm{lk}}(6)$ we get ${\textrm{lk}}(5) = C_{11}([4, 3, 2, 1, 0, 7, 6], 15, [14, 23], p)$ for some $p \in V$. Observe that, $p\in\{$13, 22$\}$. But, for both values of $p$ considering ${\textrm{lk}}(14)$ and ${\textrm{lk}}(23)$ successively we see that ${\textrm{lk}}(12)$ can not be completed. So $(m, n)\neq (15, 16)$. When $(m, n) = (16, 15)$ then completing ${\textrm{lk}}(21)$, ${\textrm{lk}}(16)$ and ${\textrm{lk}}(6)$, we get ${\textrm{lk}}(17) = C_{11}([18, 8, 9, 12, 14, 15, 16], 6, [5, 23], p)$ for some $p \in V$. Now, proceeding as in previous case, we see that the map does not exist. When $(m, n) = (16, 17)$ then ${\textrm{lk}}(21) = C_{11}([22, 23, 13, 10, 11, 19, 20], 7, [6, 16]$, $17)$. Now completing ${\textrm{lk}}(16)$ and ${\textrm{lk}}(6)$, we get ${\textrm{lk}}(5) = C_{11}([4, 3, 2, 1, 0, 7, 6], 16, [15, 23], p)$ for some $p \in \{$13, 22$\}$. In the first case when $p = 13$ then considering ${\textrm{lk}}(5)$ and ${\textrm{lk}}(13)$ successively we see that ${\textrm{lk}}(15)$ can not be completed while for $p = 22$, considering ${\textrm{lk}}(22)$, ${\textrm{lk}}(15)$ and ${\textrm{lk}}(13)$ successively we get $C_9(8, 9, 10, 13, 14, 15, 16, 17, 18) \subseteq {\textrm{lk}}(12)$. A contradiction. So, $(m, n) \neq (16, 17)$. When $(m, n) = (17, 16)$ then ${\textrm{lk}}(21) = C_{11}([22, 23, 13, 10, 11, 19, 20], 7, [6, 17]$, $16)$. Now successively considering ${\textrm{lk}}(6)$, ${\textrm{lk}}(17)$, ${\textrm{lk}}(18)$, ${\textrm{lk}}(5)$ and ${\textrm{lk}}(11)$ we see 14 as an edge and a non-edge both. So $(m, n)\neq (17, 16)$. Thus for $(x, y) = (20, 19)$ the required map does not exist.

\noindent{\bf Case 1\,:} If $(x, y) = (19, 11)$ then successively we get ${\textrm{lk}}(8) = C_{11}([18, 17, 16, 15, 14, 12, 9], 0$, $[7, 19], 11)$, ${\textrm{lk}}(7) = C_{11}([0, 1, 2, 3, 4, 5, 6], 20, [19, 8], 9)$, ${\textrm{lk}}(19) = C_{11}([20, 21, 22, 23, 13, 10, 11]$, $18, [8, 7], 6)$, ${\textrm{lk}}(11) = C_{11}([19, 20, 21, 22, 23, 13, 10], 0, [1, 18], 8)$, ${\textrm{lk}}(1) = C_{11}([0, 7, 6, 5, 4, 3, 2]$, $17, [18, 11], 10)$, ${\textrm{lk}}(18) = C_{11}([17, 16, 15, 14, 12, 9, 8], 19, [11, 1], 2)$ and ${\textrm{lk}}(6) = C_{11}([5, 4, 3, 2$, $1, 0, 7], 19, [20, m], n)$ for some $m, n \in V$. In this case, $(m, n)\in\{$(14, 12), (14, 15), (15, 14), (15, 16), (16, 15), (16, 17)$\}$. If $(m, n) = (16, 17)$ then considering ${\textrm{lk}}(17)$ we see 25 as an edge and a non-edge both and, if $(m, n) = (14, 15)$ then considering ${\textrm{lk}}(6)$ and ${\textrm{lk}}(20)$ successively we see that ${\textrm{lk}}(21)$ can not be completed. For the remaining values of $(m, n)$, we have following subcases.

\noindent{\bf Subcase 1.1\,:} When $(m, n) = (15, 14)$ then ${\textrm{lk}}(6) = C_{11}([5, 4, 3, 2, 1, 0, 7], 19, [20, 15], 14)$, ${\textrm{lk}}(15) = C_{11}([16, 17, 18, 8, 9, 12, 14], 5, [6, 20], 21)$ and ${\textrm{lk}}(20) = C_{11}([21, 22, 23, 13, 10, 11, 19]$, $7, [6, 15], 16)$. This implies ${\textrm{lk}}(21) = C_{11}([22, 23, 13, 10, 11, 19, 20], 15, [16, o], p)$ for some $o, p \in V$. Observe that $(o, p)\in\{$(3, 2), (3, 4), (4, 3), (4, 5)$\}$. In case $(o, p)\in\{$(3, 2), (3, 4)$\}$ considering ${\textrm{lk}}(21)$, ${\textrm{lk}}(16)$ and ${\textrm{lk}}(3)$ successively we see that ${\textrm{lk}}(4)$ or ${\textrm{lk}}(17)$ can not be completed. If $(o, p) = (4, 3)$ then considering ${\textrm{lk}}(21)$, ${\textrm{lk}}(4)$, ${\textrm{lk}}(16)$ successively we see that ${\textrm{lk}}(22)$ can not be completed. If $(o, p) = (4, 5)$ then successively considering ${\textrm{lk}}(21)$, ${\textrm{lk}}(5)$ and ${\textrm{lk}}(22)$, we get $C_9(9, 10, 11, 19, 20, 21, 22, 23, 12) \subseteq {\textrm{lk}}(13)$. A contradiction. So, $(m, n)\neq(15, 14)$

\noindent{\bf Subcase 1.2\,:} If $(m, n) = (14, 12)$ then successively we get ${\textrm{lk}}(6) = C_{11}([5, 4, 3, 2, 1, 0, 7], 19$, $[20, 14], 12)$, ${\textrm{lk}(14) = C_{11}}([15, 16, 17, 18, 8, 9, 12], 5, [6, 20], 21)$, ${\textrm{lk}}(20) = C_{11}([21, 22, 23, 13$, $10, 11, 19], 7, [6, 14], 15)$, ${\textrm{lk}}(12) = C_{11}([14, 15, 16, 17, 18, 9, 10], 13, [5, 6], 14)$, ${\textrm{lk}}(5) = C_{11}([6$, $7, 0, 1, 2, 3, 4], 23, [13, 12], 14)$, ${\textrm{lk}}(13) = C_{11}([23, 22, 21, 20, 19, 11, 10], 9, [12, 5], 4)$ and ${\textrm{lk}}(21)$ = $C_{11} ([22, 23, 13, 10, 11, 19, 20], 14, [15, o], p)$ for some $o, p \in V$. It is easy to see that $(o, p)\in\{$(3, 2), (3, 4)$\}$. In case $(o, p) = (3, 4)$, considering ${\textrm{lk}}(21)$ and ${\textrm{lk}}(4)$ successively we get $C_9(3, 4, 23, 13, 10, 11, 19, 20, 21)\subseteq{\textrm{lk}}(22)$. A contradiction. So $(o, p) = (3, 2)$ then completing successively we get ${\textrm{lk}}(16) = C_{11}([17, 18, 8, 9, 12, 14, 15], 3, [4, 23], 22)$, ${\textrm{lk}}(22) = C_{11}([23, 13, 10, 11, 19, 20, 21], 3, [2, 17], 16)$, ${\textrm{lk}}(4) = C_{11}([5, 6, 7, 0, 1, 2, 3], 15, [16, 23], 13)$, ${\textrm{lk}}$ $(23) = C_{11}([13, 10, 11, 19, 20, 21, 22], 17, [16, 4], 5)$, ${\textrm{lk}}(17) = C_{11}([18, 8, 9, 12, 14, 15, 16], 23$, $[22, 2], 1)$, ${\textrm{lk}}(1) = C_{11}$ $([2, 3, 4, 5, 6, 7, 0], 10, [11, 18], 17)$ and ${\textrm{lk}}(2) = C_{11}([3, 4, 5, 6, 7, 0, 1], 18$, $[17, 22], 21)$. This is $K_1(3, 4, 8, 4)$ as given in Section \ref{example 4.1}.

\noindent{\bf Subcase 1.3\,:} When $(m, n) = (15, 16)$ then ${\textrm{lk}}(6) = C_{11}([7, 0, 1, 2, 3, 4, 5], 16, [15, 20], 19)$. This implies ${\textrm{lk}}(15) = C_{11}([16, 17, 18, 8, 9, 12, 14], 21, [20, 6], 5)$, ${\textrm{lk}}(20) = C_{11}([21, 22, 23, 13$, $10, 11, 19], 7, [6, 15], 14)$ and ${\textrm{lk}}(14) = C_{11}([12, 9, 8, 18, 17, 16, 15], 20, [21, o], p)$ for some $o, p \in V$. Then, $(o, p)\in\{$(3, 2), (3, 4), (4, 3), (4, 5)$\}$. When $(o, p)\in\{$(4, 3), (4, 5)$\}$ then successively considering ${\textrm{lk}}(14)$, ${\textrm{lk}}(4)$ and ${\textrm{lk}}(21)$, it is easy to see that ${\textrm{lk}}(22)$ can not be completed. When $(o, p) = (3, 2)$ then ${\textrm{lk}}(14) = C_{11}([15, 16, 17, 18, 8, 9, 12], 2, [3, 21], 20)$, ${\textrm{lk}}(3) = C_{11}([4, 5, 6, 7, 0, 1, 2], 12, [14, 21], 22)$, ${\textrm{lk}}(21) = C_{11}([22, 23, 13, 10, 11, 19, 20], 15, [14, 3], 4)$ and ${\textrm{lk}}(22) = C_{11}([23, 13, 10, 11, 19, 20, 21], 3, [4, q], r)$ for some $q, r \in V$. This implies $q = 17$ and $r = 16$, now considering ${\textrm{lk}}(22)$, ${\textrm{lk}}(16)$, ${\textrm{lk}}(5)$ and ${\textrm{lk}}(23)$ successively we see that ${\textrm{lk}}(17)$ can not be completed. So $(o, p) = (3, 4)$ then ${\textrm{lk}}(3) = C_{11}([4, 5, 6, 7, 0, 1, 2], 22, [21, 14], 12)$, completing successively we get ${\textrm{lk}}(21) = C_{11}([22, 23, 13, 10, 11, 19, 20], 15, [14, 3], 2)$, ${\textrm{lk}}(2) = C_{11}([1, 0, 7, 6, 5, 4, 3], 21, [22, 17], 18)$, ${\textrm{lk}}(1) = C_{11}([2, 3, 4, 5, 6, 7, 0], 10, [11, 18], 17)$, ${\textrm{lk}}(22) = C_{11}([23, 13, 10, 11, 19, 20, 21], 3, [2, 17], 16)$,
${\textrm{lk}}(17) = C_{11}([18, 8, 9, 12, 14, 15, 16], 23, [22, 2], 1)$, ${\textrm{lk}}(23) = C_{11}([13, 10, 11, 19, 20, 21, 22], 17, [16, 5], 4)$, ${\textrm{lk}}(5) = C_{11}([6, 7, 0, 1, 2, 3, 4], 13, [23$, $16], 15)$, ${\textrm{lk}}(16) = C_{11}([17, 18, 8, 9, 12, 14, 15], 6, [5, 23], 22)$, ${\textrm{lk}}(13) = C_{11}([10, 11, 19, 20, 21$, $22, 23], 5, [4, 12], 9)$, ${\textrm{lk}}(4) = C_{11}([5, 6, 7, 0, 1, 2, 3], 14, [12, 13], 23)$, ${\textrm{lk}}(12) = C_{11}([14, 15, 16$, $17, 18, 8, 9], 10, [13, 4], 3)$. This is isomorphic to $K_2(3, 4, 8, 4)$, as given in Section \ref{example 4.1}, by the map (0, 23, 7, 22)(1, 13, 6, 21)(2, 10, 5, 20)(3, 11, 4, 19)(8, 14, 17, 9, 15, 18, 12, 16).

\noindent{\bf Subcase 1.4\,:} When $(m, n) = (16, 15)$ then successively we get ${\textrm{lk}}(6) = C_{11}([7, 0, 1, 2, 3, 4, 5]$, $15, [16, 20], 19)$, ${\textrm{lk}}(20) = C_{11}([21, 22, 23, 13, 10, 11, 19], 7, [6, 16], 17)$, ${\textrm{lk}}(16) = C_{11}([17, 18$, $8, 9, 12, 14, 15], 5, [6, 20], 21)$, ${\textrm{lk}}(2) = C_{11}([3, 4, 5, 6, 7, 0, 1], 18, [17, 21], 22)$, ${\textrm{lk}}(17)= C_{11}([18$, $8, 9$, $12, 14, 15, 16], 20, [21, 2], 1)$, ${\textrm{lk}}(21) = C_{11}([22, 23, 13, 10, 11, 19, 20], 16, [17, 2], 3)$ and ${\textrm{lk}}$ (15) = $C_{11}([14, 12, 9, 8, 18, 17, 16], 6, [5, o], p)$ for some $o, p \in V$. Observe that, $(o, p)\in\{$(23, 13), (23, 22)$\}$. But for $(o, p) = (23, 13)$, considering ${\textrm{lk}}(15)$ and ${\textrm{lk}}(13)$ successively we get $C_9(8, 9, 10, 13, 14, 15, 16, 17, 18)\subseteq {\textrm{lk}}(12)$. This is a contradiction. On the other hand when $(o, p) = (23, 22)$ then ${\textrm{lk}}(15) = C_{11}([16, 17, 18, 8, 9, 12, 14], 22, [23, 5], 6)$, ${\textrm{lk}}(5) = C_{11}([6, 7, 0, 1, 2, 3, 4], 13, [23, 15], 16)$, ${\textrm{lk}}(23) = C_{11}([13, 10, 11, 19, 20, 21, 22], 14, [15, 5], 4)$, completing successively we get ${\textrm{lk}}(13) = C_{11}([10, 11, 19, 20, 21, 22, 23], 5, [4, 12], 9)$, ${\textrm{lk}}(4)$ = $C_{11}$ $([3, 2, 1, 0, 7, 6, 5], 23, [13, 12], 14)$, ${\textrm{lk}}(12) = C_{11}([14, 15, 16, 17, 18, 8, 9], 10, [13, 4], 3)$, ${\textrm{lk}}(14) = C_{11}([15, 16, 17, 18, 8, 9, 12], 4, [3, 22], 23)$, ${\textrm{lk}}(3) = C_{11}([4, 5, 6, 7, 0, 1, 2], 21, [22, 14], 12)$ and
${\textrm{lk}}(22) = C_{11}([23, 13, 10, 11, 19, 20, 21], 2, [3, 14], 15)$. This is isomorphic to $K_1(3, 4, 8, 4)$ by the map (0, 11)(1, 10)(2, 13)(3, 23)(4, 22)(5, 21)(6, 20)(7, 19)(9, 18)(12, 17)(14, 16).

\noindent{\bf Case 2\,:} When $(x, y) = (20, 21)$ then ${\textrm{lk}}(8) = C_{11}([9, 12, 14, 15, 16, 17, 18], 21, [20, 7], 0)$, ${\textrm{lk}}(7) = C_{11}([0, 1, 2, 3, 4, 5, 6], 19, [20, 8], 9)$, ${\textrm{lk}}(20) = C_{11}([21, 22, 23, 13, 10, 11, 19], 6, [7, 8], 18)$ and ${\textrm{lk}}(18) = C_{11}([17, 16, 15, 14, 12, 9, 8], 20, [21, m], n)$ for some $m, n \in V$. In this case $(m, n)\in\{$(2, 1), (2, 3), (3, 2), (3, 4), (4, 3), (4, 5), (5, 4), (5, 6)$\}$. When $(m, n) = (2, 3)$ then considering ${\textrm{lk}}(18)$, ${\textrm{lk}}(2)$, ${\textrm{lk}}(21)$ and ${\textrm{lk}}(1)$ successively we see that 11\,22 is simultaneously an edge and a non-edge of $K$. When $(m, n) = (5, 4)$ then considering ${\textrm{lk}}(18)$, ${\textrm{lk}}(21)$ and ${\textrm{lk}}(6)$ successively we see that 19\,22 is both an edge and a non-edge of $K$. So, $(m, n)\neq (2, 3)$, $(5, 4)$. For the remaining values of $(m, n)$ we have following subcases.

When $(m, n) = (4, 5)$ then we have ${\textrm{lk}}(18) = C_{11}([8, 9, 12, 14, 15, 16, 17], 5, [4, 21], 20)$, ${\textrm{lk}}(4) = C_{11}([5, 6, 7, 0, 1, 2, 3], 22, [21, 18], 17)$, ${\textrm{lk}}(21) = C_{11}([22, 23, 13, 10, 11, 19, 20],8, [18$, $4], 3)$ and ${\textrm{lk}}(22) = C_{11}([23, 13, 10, 11, 19, 20, 21], 4, [3, o], p)$ for some $o, p \in V$. Observe that, $(o, p)\in\{$(14, 15), (15, 14), (15, 16), (16, 15)$\}$. If $(o, p) = (14, 15)$ then successively considering ${\textrm{lk}}(22)$, ${\textrm{lk}}(14)$, ${\textrm{lk}}(3)$, ${\textrm{lk}}(12)$, ${\textrm{lk}}(13)$ and ${\textrm{lk}}(23)$ we see that $\deg(1) > 4$. A contradiction. If $(o, p)\in\{$(15, 14), (15, 16)$\}$ then considering ${\textrm{lk}}(22)$, ${\textrm{lk}}(15)$ and ${\textrm{lk}}(3)$ successively we see that ${\textrm{lk}}(16)$ or ${\textrm{lk}}(2)$ can not be completed. If $(o, p) = (16, 15)$ then considering ${\textrm{lk}}(22)$, ${\textrm{lk}}(16)$ and ${\textrm{lk}}(3)$ successively we see that ${\textrm{lk}}(2)$ can not be completed. So, $(m, n) \neq (4, 5)$. When $(m, n) = (3, 2)$ then ${\textrm{lk}}(18) = C_{11}([8, 9, 12, 14, 15, 16, 17], 2, [3, 21], 20)$. This implies ${\textrm{lk}}(3) = C_{11}([4, 5, 6, 7, 0, 1, 2], 17, [18, 21], 22)$, ${\textrm{lk}}(21) = C_{11}([22, 23, 13, 10, 11, 19, 20], 8, [18, 3], 4)$ and ${\textrm{lk}}(22) = C_{11}([23, 13, 10, 11, 19, 20, 21], 3, [4, o], p)$ for some $o, p \in V$. Observe that, $(o, p)\in\{$(14, 12), (14, 15), (15, 14), (15, 16), (16, 15)$\}$. Now proceeding further as in previous case we get a contradiction for each value of $(o, p)$. So $(m, n)\neq(3, 2)$.

\noindent{\bf Subcase 2.1\,:} When $(m, n) = (2, 1)$ then successively we get ${\textrm{lk}}(18) = C_{11}([8, 9, 12, 14, 15$, $16, 17], 1, [2, 21], 20)$, ${\textrm{lk}}(21) = C_{11}([22, 23, 13, 10, 11, 19, 20], 8, [18, 2], 3)$, ${\textrm{lk}}(2) = C_{11}([3, 4, 5$, $6, 7, 0, 1], 17, [18, 21], 22)$, ${\textrm{lk}}(1) = C_{11}([2, 3, 4, 5, 6, 7, 0], 10, [11, 17], 18)$, ${\textrm{lk}}(11) = C_{11}([19, 20$, $21, 22, 23, 13, 10], 0, [1, 17], 16)$, ${\textrm{lk}}(17) = C_{11}([18, 8, 9, 12, 14, 15, 16], 19, [11, 1], 2)$ and ${\textrm{lk}}(3) = C_{11}([4, 5, 6, 7, 0, 1, 2], 21, [22, o], p)$ for some $o, p \in V$. In this case we have $(o, p)\in\{$(14, 12), (14, 15), (15, 14)$\}$.

If $(o, p) = (14, 15)$ then ${\textrm{lk}}(3) = C_{11}([4, 5, 6, 7, 0, 1, 2], 21, [22, 14], 15)$,now completing ${\textrm{lk}}(14)$ and ${\textrm{lk}}(22)$ we get ${\textrm{lk}}(23) = C_{11}([13, 10, 11, 19, 20, 21, 22], 14, [12, 5], r)$ for some $r \in V$. It is easy to see that $r \in \{$4, 6$\}$. If $r = 4$ then successively considering ${\textrm{lk}}(23)$, ${\textrm{lk}}(4)$ and ${\textrm{lk}}(15)$ we get $\deg(13) > 4$ and if $r = 6$ then considering ${\textrm{lk}}(23)$ and ${\textrm{lk}}(6)$ successively we see that 13\,19 is both an edge and a non-edge of $K$. So $(o, p) \neq (14, 15)$. When $(o, p) = (15, 14)$ then ${\textrm{lk}}(3) = C_{11}([4, 5, 6, 7, 0, 1, 2], 21, [22, 15], 14)$, completing ${\textrm{lk}}(15)$ and ${\textrm{lk}}(22)$ we get ${\textrm{lk}}(23) = C_{11}([13, 10, 11, 19, 20, 21, 22], 15, [16, 5], r)$ for some $r \in V$. Observe that $r \in \{$4, 6$\}$. Now proceeding as in previous case, we get a contradiction for each value of $r$. So $(o, p)\neq(15, 14)$.

If $(o, p) = (14, 12)$ then we see that ${\textrm{lk}}(3) = C_{11}([4, 5, 6, 7, 0, 1, 2], 21, [22, 14],12)$,
${\textrm{lk}}(14) = C_{11}([15, 16, 17, 18, 8, 9, 12], 4, [3, 22], 23)$, now completing successively we get
${\textrm{lk}}(12) = C_{11}([14$, $15, 16, 17, 18, 8, 9], 10, [13, 4], 3)$, ${\textrm{lk}}(13) = C_{11}([10, 11, 19, 20, 21, 22, 23], 5, [4, 12], 9)$, ${\textrm{lk}}(4) = C_{11}([5, 6, 7, 0, 1, 2, 3], 14, [12, 13],23)$,
${\textrm{lk}}(23) = C_{11}([13, 10, 11, 19, 20, 21, 22],14, [15, 5], 4)$, ${\textrm{lk}}$ $(5) = C_{11}([6, 7, 0, 1, 2, 3, 4], 13, [23, 15], 16)$, ${\textrm{lk}}(16) = C_{11}([17, 18, 8, 9, 12, 14, 15], 5, [6, 19], 11)$,
${\textrm{lk}}(6) = C_{11}([7, 0, 1, 2, 3, 4, 5], 15, [16, 19], 20)$ and ${\textrm{lk}}(19) = C_{11}([20, 21, 22, 23, 13, 10, 11], 17$, $[16, 6], 7)$. This is $K_2(3, 4, 8, 4)$ as given in Section \ref{example 4.1}.

\noindent{\bf Subcase 2.2\,:} When $(m, n) = (3, 4)$ then successively we get ${\textrm{lk}}(18) = C_{11}([8, 9, 12, 14, 15$, $16, 17], 4, [3, 21], 20)$, ${\textrm{lk}}(3) = C_{11}([4, 5, 6, 7, 0, 1, 2], 22, [21, 18], 17)$,
${\textrm{lk}}(21) = C_{11}([22, 23, 13$, $10, 11, 19, 20], 8, [18, 3], 2)$ and ${\textrm{lk}}(17) = C_{11}([16, 15, 14, 12, 9, 8, 18], 3, [4, o], p)$ for some $o, p \in V$. In this case $(o, p)\in\{$(13, 23), (19, 11), (23, 13), (23, 22)$\}$. If $(o, p) = (13, 23)$ then considering ${\textrm{lk}}(17)$ and ${\textrm{lk}}(13)$ successively we see that 12\,17 is both an edge and a non-edge of $K$. If $(o, p) = (19, 11)$ then successively considering ${\textrm{lk}}(17)$, ${\textrm{lk}}(11)$ and ${\textrm{lk}}(1)$ we see easily that ${\textrm{lk}}(4)$ can not be completed. If $(o, p) = (23, 13)$ then considering ${\textrm{lk}}(17)$ and ${\textrm{lk}}(13)$ we see that 12\,16 is both an edge and a non-edge of $K$. If $(o, p) = (23, 22)$ then ${\textrm{lk}}(17) = C_{11}([18, 8, 9, 12, 14, 15, 16], 22, [23, 4], 3)$. This implies ${\textrm{lk}}(22) = C_{11}([23, 13, 10, 11, 19, 20, 21]$, $3, [2, 16], 17)$, completing successively we get ${\textrm{lk}}(4) = C_{11}([5, 6, 7, 0, 1, 2, 3], 18, [17, 23], 13)$, ${\textrm{lk}}(13) = C_{11}([10, 11, 19, 20, 21, 22, 23], 4, [5, 12], 9)$, ${\textrm{lk}}(12) = C_{11}([14, 15, 16, 17, 18, 8, 9], 10$, $[13, 5], 6)$, ${\textrm{lk}}(6) = C_{11}([7, 0, 1, 2, 3, 4, 5], 12, [14, 19], 20)$, ${\textrm{lk}}(19) = C_{11}([20, 21, 22, 23, 13, 10$, $11], 15, [14, 6], 7)$, ${\textrm{lk}}(14) = C_{11}([15, 16, 17, 18, 8, 9, 12], 5, [6, 19], 11)$, ${\textrm{lk}}(5) = C_{11}([6, 7, 0, 1, 2$, $3, 4], 23, [13, 12], 14)$, ${\textrm{lk}}(23) = C_{11}([13, 10, 11, 19, 20, 21, 22], 16, [17, 4], 5)$, ${\textrm{lk}}(2) = C_{11}([3, 4, 5$, $6, 7, 0, 1], 15, [16, 22], 21)$, ${\textrm{lk}}(15) = C_{11}([16, 17, 18, 8, 9, 12, 14], 19, [11, 1], 2)$, ${\textrm{lk}}(11) = C_{11}([19$, $20, 21, 22, 23, 13, 10], 0, [1, 15], 14)$, ${\textrm{lk}}(1) = C_{11}([2, 3, 4, 5, 6, 7, 0], 10, [11, 15], 16)$ and ${\textrm{lk}}(16) = C_{11}([17, 18, 8, 9, 12, 14, 15], 1, [2, 22], 23)$. This is isomorphic to $K_2(3, 4, 8, 4)$ by the map (0, 20, 16, 3, 23, 12)(1, 21, 15, 2, 22, 14)(4, 13, 9, 7, 19, 17)(5, 10, 8, 6, 11, 18).

\noindent{\bf Subcase 2.3\,:} When $(m, n) = (4, 3)$ then successively we get ${\textrm{lk}}(18) = C_{11}([8, 9, 12, 14, 15$,
$16, 17], 3, [4, 21], 20)$, ${\textrm{lk}}(4) = C_{11}([5, 6, 7, 0, 1, 2, 3], 17, [18, 21], 22)$, ${\textrm{lk}}(21) = C_{11}([22, 23, 13$,
$10, 11, 19, 20], 8, [18, 4], 5)$ and ${\textrm{lk}}(22) = C_{11}([23, 13, 10, 11, 19, 20, 21], 4, [5, o], p)$ for some $o, p \in V$.
Then we see that $(o, p)\in\{$(14, 15), (15, 14), (15, 16), (16, 15), (16, 17)$\}$. If $(o, p) = (14, 15)$ then successively considering
${\textrm{lk}}(22)$, ${\textrm{lk}}(14)$, ${\textrm{lk}}(6)$ we see that ${\textrm{lk}}(12)$ can not be completed. If $(o, p) = (15, 14)$
then successively considering ${\textrm{lk}}(22)$, ${\textrm{lk}}(15)$, ${\textrm{lk}}(6)$, ${\textrm{lk}}(19)$ and ${\textrm{lk}}(17)$ we
see that ${\textrm{lk}}(11)$ can not be completed. If $(o, p) = (15, 16)$ then successively considering ${\textrm{lk}}(15)$,
${\textrm{lk}}(6)$, ${\textrm{lk}}(19)$ and ${\textrm{lk}}(12)$ we see that 11\,13 is both an edge and a non-edge of $K$. If $(o, p) = (16, 15)$
then considering ${\textrm{lk}}(22)$, ${\textrm{lk}}(16)$ and ${\textrm{lk}}(6)$ successively we get $\deg(17) > 4$. If $(o, p) = (16, 17)$
then ${\textrm{lk}}(22) = C_{11}([23, 13, 10, 11, 19, 20, 21], 4, [5, 16], 17)$, completing successively, we get ${\textrm{lk}}(17) = C_{11}
([18, 8, 9, 12, 14, 15, 16], 22, [23, 3], 4)$, ${\textrm{lk}}(4) = C_{11}([5, 6, 7, 0, 1, 2, 3], 17, [18, 21], 22)$, ${\textrm{lk}}(3) = C_{11}
([4, 5, 6, 7, 0, 1, 2], 13, [23, 17], 18)$, ${\textrm{lk}}(13) = C_{11}([10, 11, 19, 20, 21, 22, 23], 3, [2, 12], 9)$, ${\textrm{lk}}(2) =
C_{11}([3, 4, 5, 6, 7, 0, 1], 14, [12, 13], 23)$, ${\textrm{lk}}(12)$ $= C_{11}([14, 15, 16, 17, 18, 8, 9], 10, [13, 2], 1)$, ${\textrm{lk}}(23)
 = C_{11}([13, 10, 11, 19, 20, 21, 22], 16, [17$, $3], 2)$, ${\textrm{lk}}(16) = C_{11}([17, 18, 8, 9, 12, 14, 15], 6, [5, 22], 23)$,
${\textrm{lk}}(6) = C_{11}([7, 0, 1, 2, 3, 4, 5], 16, [15$ ,$19], 20)$, ${\textrm{lk}}(19) = C_{11}([20, 21, 22, 23, 13, 10, 11], 14, [15, 6], 7)$,
${\textrm{lk}}(14) = C_{11}([15, 16, 17, 18$, $8, 9, 12], 2, [1, 11], 19)$, ${\textrm{lk}}(11) = C_{11}([19, 20, 21, 22, 23, 13, 10], 0, [1, 14], 15)$
and ${\textrm{lk}}(1) = C_{11}([2$, $3, 4, 5, 6, 7, 0], 10, [11, 14], 12)$. This is isomorphic to $K_1(3, 4, 8, 4)$ by the map (0, 7, 6, 5, 4, 3,
2, 1)(8, 20, 14, 10)(9, 19, 12, 11)(13, 18, 21, 15)(16, 23, 17, 22).

\noindent{\bf Subcase 2.4\,:} When $(m, n) = (5, 6)$ then successively we get ${\textrm{lk}}(18) = C_{11}([8, 9, 2, 14$, $15, 16, 17], 6, [5, 21], 20)$,
${\textrm{lk}}(6) = C_{11}([7, 0, 1, 2, 3, 4, 5], 18, [17, 19], 20)$, ${\textrm{lk}}(17) = C_{11}([18, 8$, $9, 12, 14, 15, 16], 11, [19, 6], 5)$,
${\textrm{lk}}(11) = C_{11}([19, 20, 21, 22, 23, 13, 10], 0, [1, 16], 17)$, ${\textrm{lk}}(1) = C_{11}([2, 3, 4, 5, 6, 7, 0], 10, [11, 16], 15)$,
${\textrm{lk}}(16) = C_{11}([17, 18, 8, 9, 12, 14, 15], 2, [1, 11], 19)$, \linebreak ${\textrm{lk}}(19) = C_{11}([20, 21, 22, 23, 13, 10, 11], 16, [17, 6], 7)$,
${\textrm{lk}}(5) = C_{11}([6, 7, 0, 1, 2, 3, 4], 22, [21, 18]$, $17)$, ${\textrm{lk}}(21) = C_{11}$ $([22, 23, 13, 10, 11, 19, 20], 8, [18, 5], 4)$ and
${\textrm{lk}}(22) = C_{11}([23, 13, 10, 11, 19$, $20, 21], 5, [4, o], p)$ for some $o, p \in V$. Observe that  $(o, p)\in\{$(14, 12), (14, 15)$\}$. In case $(o, p) = (14, 12)$, we get $C_9(9, 10, 11, 19, 20, 21, 22, 23, 12)\subseteq {\textrm{lk}}(13)$. A contradiction. So $(o, p) = (14, 15)$ then ${\textrm{lk}}(22) = C_{11}([23, 13, 10, 11, 19, 20, 21], 5, [4, 14], 15)$. This implies ${\textrm{lk}}(15) = C_{11}([16, 17, 18, 8, 9, 12, 14], 22, [23, 2], 1)$, completing successively, we get ${\textrm{lk}}(2) = C_{11}([3, 4, 5, 6, 7, 0, 1], 16, [15, 23], 13)$,
${\textrm{lk}}(13) = C_{11}([10, 11, 19, 20, 21, 22, 23], 2, [3, 12], 9)$,\linebreak ${\textrm{lk}}(12) = C_{11}([14, 15, 16, 17, 18, 8, 9], 10, [13, 3], 4)$,
${\textrm{lk}}(4) = C_{11}([5, 6, 7, 0, 1, 2, 3], 12, [14, 22]$, $ 21)$ and ${\textrm{lk}}(14) = C_{11}([15, 16, 17, 18, 8, 9, 12], 3, [4, 22], 23)$. This is isomorphic to $K_1(3, 4, 8, 4)$ by the map (0, 6)(1, 5)(2, 4)(8, 19, 9, 20)(10, 14, 22, 17)(11, 12, 21, 18)(13, 15, 23, 16).

\noindent{\bf Case 3\,:} When $(x, y) = (23, 22)$ then we get ${\textrm{lk}}(8) = C_{11}([9, 12, 14, 15, 16, 17, 18], 22, [23, 7]$, $0)$, ${\textrm{lk}}(23) = C_{11}([13, 10, 11, 19, 20, 21, 22], 18, [8, 7], 6)$,
${\textrm{lk}}(7) = C_{11}([0, 1, 2, 3, 4, 5, 6], 13, [23, 8]$, $9)$. This implies
${\textrm{lk}}(13) = C_{11}([10, 11, 19, 20, 21, 22, 23], 7, [6, 12], 9)$,
${\textrm{lk}}(6) = C_{11}([7, 0, 1, 2, 3, 4$, $5], 14, [12, 13], 23)$,
${\textrm{lk}}(12) = C_{11}([14, 15, 16, 17, 18, 8, 9], 10, [13, 6], 5)$ and
${\textrm{lk}}(5) = C_{11}([4, 3, 2, 1$, $0, 7, 6], 12, [14, l], m)$ for some $m, l \in V$.
It is easy to see that $(l, m)\in\{$(19, 20), (20, 19), (20, 21), (21, 20), (21, 22)$\}$.

When $(l, m) = (21, 20)$ then ${\textrm{lk}}(5) = C_{11}([6, 7, 0, 1, 2, 3, 4], 20, [21, 14], 12)$. Now considering ${\textrm{lk}}(21)$ and ${\textrm{lk}}(14)$ successively we see that ${\textrm{lk}}(22)$ can not be completed. When $(l, m) = (20, 21)$ then ${\textrm{lk}}(5) = C_{11}([6, 7, 0, 1, 2, 3, 4], 21, [20, 14], 12)$. This implies
${\textrm{lk}}(14) = C_{11}([15, 16, 17, 18, 8, 9, 12], 6, [5, 20], 19)$,
${\textrm{lk}}(20) = C_{11}([21, 22, 23, 13, 10, 11, 19], 15, [14, 5], 4)$ and
${\textrm{lk}}(21) = C_{11}([22, 23, 13, 10, 11, 19, 20], 5, [4, n], o)$ for some $n, o \in V$.
Observe that $(n, o)\in\{$(16, 17), (17, 16)$\}$. If $(n, o) = (16, 17)$ then considering ${\textrm{lk}}(21)$ and ${\textrm{lk}}(22)$ successively we get $C_9(8, 9, 12, 14, 15, 16, 21, 22, 18) \subseteq {\textrm{lk}}(17)$ and if $(n, o) = (17, 16)$ then successively considering ${\textrm{lk}}(21)$, ${\textrm{lk}}(17)$ and ${\textrm{lk}}(4)$ we see easily that ${\textrm{lk}}(22)$ can not be completed.

This implies ${\textrm{lk}}(14) = C_{11}([15, 16, 17, 18, 8, 9, 12], 6, [5, 19], 11)$, completing successively we get ${\textrm{lk}}(11) = C_{11}([19, 20, 21, 22, 23, 13, 10],0, [1, 15], 14)$,
${\textrm{lk}}(1) = C_{11}([2, 3, 4, 5, 6, 7, 0], 10$, $[11, 15], 16)$,
${\textrm{lk}}(15) = C_{11}([16, 17, 18, 8, 9, 12, 14], 19, [11, 1], 2)$,
${\textrm{lk}}(16) = C_{11}([17, 18, 8, 9, 12$, $14, 15], 1, [2, 21], 20)$,
${\textrm{lk}}(2) = C_{11}([3, 4, 5, 6, 7, 0, 1], 15, [16, 21], 22)$,
${\textrm{lk}}(21) = C_{11}([22, 23, 13$, $10, 11, 19, 20], 17, [16, 2], 3)$,
${\textrm{lk}}(22) = C_{11}([23, 13, 10, 11, 19, 20, 21], 2, [3, 18], 8)$,
${\textrm{lk}}(18) = C_{11}$ $([8, 9, 12, 14, 15, 16, 17], 4, [3, 22], 23)$,
${\textrm{lk}}(3) = C_{11}([4, 5, 6, 7, 0, 1, 2], 21, [22, 18], 17)$,
${\textrm{lk}}(4) = C_{11}([5, 6, 7, 0, 1, 2, 3], 18, [17, 20], 19)$,
${\textrm{lk}}(17) = C_{11}$ $([18, 8, 9, 12, 14, 15, 16], 21, [20, 4], 3)$ and
${\textrm{lk}}(20) = C_{11}([21, 22, 23, 13, 10, 11, 19], 5, [4, 17], 16)$. This is isomorphic to $K_1(3, 4, 8, 4)$ by the map (0, 21, 8, 2, 19, 12, 4, 10, 15, 6, 23, 17)(1, 20, 9, 3, 11, 14, 5, 13, 16, 7, 22, 18).

\noindent{\bf Subcase 3.2\,:} When $(l, m) = (20, 19)$ then
${\textrm{lk}}(5) = C_{11}([6, 7, 0, 1, 2, 3, 4], 19, [20, 14], 12)$. This implies
${\textrm{lk}}(20) = C_{11}([21, 22, 23, 13, 10, 11, 19], 4, [5, 14], 15)$ and
${\textrm{lk}}(4) = C_{11}([3, 2, 1, 0$, $7, 6, 5], 20, [19, n], o)$ for some $n, o \in V$.
Observe that $(n, o)\in\{$(16, 15), (16, 17)$\}$. If $(n, o) = (16, 17)$ then
considering ${\textrm{lk}}(4)$, ${\textrm{lk}}(16)$ and ${\textrm{lk}}(11)$ successively we see that ${\textrm{lk}}(1)$ can not be completed. On the other hand when $(n, o) = (16, 15)$ then ${\textrm{lk}}(4) = C_{11}([5, 6, 7, 0, 1, 2, 3], 15$, $[16, 19], 20)$. This implies ${\textrm{lk}}(15) = C_{11}([16, 17, 18, 8, 9, 12, 14], 20, [21, 3], 4)$, completing successively we get ${\textrm{lk}}(21) = C_{11}([22, 23, 13, 10, 11, 19, 20], 14, [15, 3], 2)$, ${\textrm{lk}}(3) = C_{11}([4, 5, 6, 7, 0$, $1, 2], 22, [21, 15], 16)$,
${\textrm{lk}}(16) = C_{11}([17, 18, 8, 9, 12, 14, 15], 3, [4, 19], 11)$,
${\textrm{lk}}(11) = C_{11}([19, 20$, $21, 22, 23, 13, 10], 0, [1, 17], 16)$,
${\textrm{lk}}(1) = C_{11}([2, 3, 4, 5, 6, 7, 0], 10, [11, 17], 18)$,
${\textrm{lk}}(18) = C_{11}([8$, $9, 12, 14, 15, 16, 17], 1, [2, 22], 23)$,
${\textrm{lk}}(2) = C_{11}([3, 4, 5, 6, 7, 0, 1], 17, [18, 22], 21)$ and
${\textrm{lk}}(22) = C_{11}$ $([23, 13, 10, 11, 19, 20, 21], 3, [2, 18], 8)$. This is isomorphic to $K_1(3, 4, 8, 4)$ by the map (0, 9)(1, 12)(2, 14)(3, 15)(4, 16)(5, 17)(6, 18)(7, 8)(11, 13)(19, 23)(20, 22).

\noindent{\bf Subcase 3.3\,:} When $(l, m) = (21, 22)$ then ${\textrm{lk}}(5) = C_{11}([6, 7, 0, 1, 2, 3, 4], 22, [21, 14], 12)$. This implies ${\textrm{lk}}(22) = C_{11}([23, 13, 10, 11, 19, 20, 21], 5, [4, 18], 8)$,
${\textrm{lk}}(18) = C_{11}([8, 9, 12, 14, 15$, $16, 17], 3, [4, 22], 23)$ and
${\textrm{lk}}(17) = C_{11}([16, 15, 14, 12, 9, 8, 18], 4, [3, n], o)$ for some $n, o \in V$.
Then we see that $(n, o)\in\{$(19, 11), (19, 20), (20, 19), (20, 21)$\}$. If $(n, o) = (19, 20)$ then
successively considering ${\textrm{lk}}(17)$, ${\textrm{lk}}(19)$ and ${\textrm{lk}}(3)$ we see that ${\textrm{lk}}(11)$ can not be completed. If $(n, o) = (20, 19)$ then successively considering ${\textrm{lk}}(17)$, ${\textrm{lk}}(20)$ and ${\textrm{lk}}(3)$ we see that ${\textrm{lk}}(19)$ can not be completed. If $(n, o) = (20, 21)$ then successively considering ${\textrm{lk}}(17)$, ${\textrm{lk}}(21)$, ${\textrm{lk}}(3)$ and ${\textrm{lk}}(20)$ we see that ${\textrm{lk}}(16)$ can not be completed.
If $(o, p) = (19, 11)$ then ${\textrm{lk}}(17) = C_{11}([18, 8, 9, 12, 14, 15, 16], 11, [19, 3], 4)$,
${\textrm{lk}}(3) = C_{11}([4, 5, 6, 7, 0, 1, 2], 20, [19, 17], 18)$ and
${\textrm{lk}}$ $(19) = C_{11}([20, 21, 22, 23, 13, 10, 11], 16, [17, 3], 2)$, completing successively we get
${\textrm{lk}}(11) = C_{11}([19, 20, 21, 22, 23, 13, 10], 0, [1, 16], 17)$,
${\textrm{lk}}(1) = C_{11}([2, 3, 4, 5, 6, 7, 0], 10, [11, 16], 15)$,
${\textrm{lk}}$ $(16) = C_{11}([17, 18, 8, 9, 12, 14, 15], 2, [1, 11], 19)$,
${\textrm{lk}}(2) = C_{11}([3, 4, 5, 6, 7, 0, 1], 16, [15, 20], 19)$,
${\textrm{lk}}(15) = C_{11}([16, 17, 18, 8, 9, 12, 14], 21, [20, 2], 1)$,
${\textrm{lk}}(20) = C_{11}([21, 22, 23, 13, 10, 11, 19], 3$, $[2, 15], 14)$. This is isomorphic to $K_1(3, 4, 8, 4)$ by the map (0, 18, 20, 4, 14, 13)(1, 17, 21, 5, 12, 10)(2, 16, 22, 6, 9, 11)(3, 15, 23, 7, 8, 19).
Thus the Lemma \ref{lem1} is proved.

\noindent {\bf Proof of Lemma \ref{lem2}\,:} Let $M$ be a SEM of type $(4, 6, 16)$ on the surface of Euler characteristic $-1$. The notation ${\textrm{lk}}(i) = C_{20}([{\textbf{\textit{i}}}_{\textbf{1}}, i_2, i_3, i_4, i_5, i_6, i_7, i_8, i_9, i_{10}, i_{11}, i_{12}, i_{13}, i_{14}, {\textbf{\textit{i}}}_{\textbf{15}}, i_{16}$, $ {\textbf{\textit{i}}}_{\textbf{17}}, i_{18}, i_{19}, i_{20})$ for the link of $i$ will mean that $[i, i_{15}, i_{16}, i_{17}],[i, i_{1}, i_{20}, i_{19}, i_{18}, i_{17}]$ and $[i, i_{1}, i_{2}, i_{3}, i_{4}, i_{5}, i_6, i_7, i_8, i_9$ $, i_{10}, i_{11}, $ $i_{12}, i_{13}, i_{14}, i_{15}]$ form 4-gonal face (face with 4-gonal boundary), 6-gonal face (face with 6-gonal boundary) and 16-gonal face (face with 16-gonal boundary), respectively. Let $|V|$ denote the number of vertices in $V(M)$. If $E(M)$, $Q(M)$, $H(M)$ and $P(M)$ denote the number of edges, number of 4-gonal faces, number of 6-gonal faces and number of 16-gonal faces, respectively, then $E(M) = \frac{3|V|}{2}$, $Q(M) = \frac{|V|}{4}$, $H(M) = \frac{|V|}{6}$ and $P(M) = \frac{|V|}{16}$. By Euler's equation we see if the map exists then $|V| = 48$. For this, let $V = V(M) = \{$0, 1, \ldots, 47$\}$. Now, we prove the result by exhaustive search for all $M$.

Assume, ${\textrm{lk}}(0) = C_{20}([1, 2, 3, 4, 5, 6, 7, 8, 9, 10, 11, 12, 13, 14, {\textbf{15}}], 16, {\textbf{17}}, 18, 19, 20)$ then successively we get ${\textrm{lk}}(17) = C_{20}([{\textbf{18}}, 21, 22, 23, 24, 25, 26, 27, 28, 29, 30, 31, 32, 33, {\textbf{16}}], 15, {\textbf 0}$, $1, 20, 19)$, ${\textrm{lk}}(18) = C_{20}([{\textbf{17}}, 16, 33, 32, 31, 30, 29, 28, 27, 26, 25, 24, 23, 22, {\textbf{21}}], 34, {\textbf{19}}, 20, 1, 0)$, ${\textrm{lk}}(19)$ $= C_{20}([{\textbf{20}}, 35, 36, 37, 38, 39, 40, 41, 42, 43, 44, 45, 46, 47, {\textbf{34}}], 21, {\textbf{18}}, 17, 0, 1)$, ${\textrm{lk}}(20) = C_{20}([{\textbf{19}}, 34, 47, 46, 45, 44, 43, 42, 41, 40, 39, 38, 37, 36, {\textbf{35}}], 2, {\textbf 1}, 0, 17, 18)$ and ${\textrm{lk}}(1) = C_{20}([{\textbf{0}}$, $15, 14, 13, 12, 11, 10, 9, 8, 7, 6, 5, 4, 3, {\textbf{2}}], 35, {\textbf{20}}, 19, 18, 17)$. This implies ${\textrm{lk}}(2) = C_{20}([{\textbf{3}}, 4, 5$, $6, 7, 8, 9, 10, 11, 12, 13, 14, 15, 0, {\textbf{1}}], 20, {\textbf{35}}, 36, d, c)$ for some $c, d \in V$. Then we see that $(c, d) \in \{$(23, 24), (24, 23), (25, 26), (26, 25), (27, 28), (28, 27), (29, 30), (30, 29), (31, 32), (32, 31)$\}$. Observe that, $(29, 30)\cong(25, 26)$ by the map (0, 19)(1, 20)(2, 35)(3, 36)(4, 37)(5, 38)(6, 39)(7, 40)(8, 41)(9, 42)(10, 43) (11, 44)(12, 45)(13, 46)(14, 47)(15, 34)(16, 21)(17, 18)(22, 33)(23, 32)(24, 31)(25, 30)(26, 29)(27, 28); $(31, 32)\cong(23, 24)$ by the map (0, 3)(1, 2)(4, 15)(5, 14)(6, 13)(7, 12)(8, 11)(9, 10)(16, 22, 26, 30)(17, 23, 27, 31)(18, 24, 28, 32)(19, 36)(20, 35)(21, 25, 29, 33) (34, 37)(38, 47)(39, 46)(40, 45)(41, 44)(42, 43); $(30, 29)\cong(26, 25)$ by the map (0, 36)(1, 35)(2, 20)(3, 19)(4, 34)(5, 47) (6, 46)(7, 45)(8, 44)(9, 43)(10, 42)(11, 41)(12, 40)(13, 39)(14, 38)(15, 37)(16, 24, 30, 18, 26, 32, 22, 28)(17, 25, 31, 21, 27, 33, 23, 29). So we have $(c, d)\in\{$(24, 23), (25, 26), (27, 28), (28, 27), (30, 29), (31, 32), (32, 31)$\}$.

If $(c, d) = (24, 23)$ then successively considering ${\textrm{lk}}(2)$, ${\textrm{lk}}(3)$, ${\textrm{lk}}(23)$, ${\textrm{lk}}(24)$, ${\textrm{lk}}(35)$ and ${\textrm{lk}}(36)$ we see that ${\textrm{lk}}(21)$ and ${\textrm{lk}}(22)$ can not be completed. If $(c, d) = (25, 26)$ then successively considering ${\textrm{lk}}(2)$, ${\textrm{lk}}(3)$, ${\textrm{lk}}(25)$, ${\textrm{lk}}(26)$, ${\textrm{lk}}(35)$, ${\textrm{lk}}(36)$ we see that ${\textrm{lk}}(23)$ and ${\textrm{lk}}(24)$ can not be completed. If $(c, d) = (28, 27)$ then successively considering ${\textrm{lk}}(2)$, ${\textrm{lk}}(3)$, ${\textrm{lk}}(27)$, ${\textrm{lk}}(28)$, ${\textrm{lk}}(35)$, ${\textrm{lk}}(36)$ we see that ${\textrm{lk}}(25)$ and ${\textrm{lk}}(26)$ can not be completed.
If $(c, d) = (32, 31)$ then successively considering ${\textrm{lk}}(2)$, ${\textrm{lk}}(3)$, ${\textrm{lk}}(31)$, ${\textrm{lk}}(32)$, ${\textrm{lk}}(35)$, ${\textrm{lk}}(36)$, we see that ${\textrm{lk}}(16)$ and ${\textrm{lk}}(33)$ can not be completed. So we search for $(c, d) \in \{$(27, 28), (30, 29), (31, 32)$\}$.
\smallskip

\noindent{\bf Case 1\,:} If $(c, d) = (27, 28)$ then constructing successively we get ${\textrm{lk}}(2) = C_{20}([{\textbf{3}}, 4, 5, 6, 7$, $8, 9, 10, 11, 12, 13, 14, 15, 0, {\textbf{1}}], 20, {\textbf{35}}, 36, 28, 27)$, ${\textrm{lk}}(3) = C_{20}([{\textbf{2}}, 1, 0, 15, 14, 13, 12, 11, 10, 9$, $8, 7, 6, 5, {\textbf{4}}], 26, {\textbf{27}}, 28, 36, 35)$,
${\textrm{lk}}(27) = C_{20}([{\textbf{28}}, 29, 30, 31, 32, 33, 16, 17, 18, 21, 22, 23, 24, 25$, ${\textbf{26}}], 4, {\textbf{3}}, 2, 35, 36)$, ${\textrm{lk}}(28) = C_{20}([{\textbf{27}}, 26, 25, 24, 23, 22, 21, 18, 17, 16, 33, 32, 31, 30, {\textbf{29}}], 37, {\textbf{36}}$, $35, 2, 3)$, ${\textrm{lk}}(35)$ $= C_{20}([{\textbf{36}}, 37, 38, 39, 40, 41, 42, 43, 44, 45, 46, 47, 34, 19, {\textbf{20}}], 1, {\textbf{2}}, 3, 27, 28)$, ${\textrm{lk}}$ $(36) = C_{20}([{\textbf{35}}, 20, 19, 34, 47, 46, 45, 44, 43, 42, 41, 40, 39, 38, {\textbf{37}}], 29, {\textbf{28}}, 27, 3, 2)$ and ${\textrm{lk}}(21) = C_{20}([{\textbf{22}}, 23, 24, 25, 26, 27, 28, 29, 30, 31, 32, 33, 16, 17, {\textbf{18}}], 19, {\textbf{34}}, 47, k, j)$ for some $k, j \in V$. Observe that $(k, j) \in \{$(6, 7), (7, 6), (8, 9), (9, 8), (10, 11), (11, 10), (12, 13), (13, 12)$\}$. If $(k, j) = (7, 6)$ then completing ${\textrm{lk}}(6)$, ${\textrm{lk}}(7)$, ${\textrm{lk}}(21)$, ${\textrm{lk}}(22)$, ${\textrm{lk}}(34)$, ${\textrm{lk}}(47)$, ${\textrm{lk}}(4)$ it is easy to see that ${\textrm{lk}}(24)$ and ${\textrm{lk}}(25)$ can not be completed.
If $(k, j) \in \{$(8, 9), (11, 10)$\}$ then completing ${\textrm{lk}}(21)$, ${\textrm{lk}}(22)$, ${\textrm{lk}}(34)$,
${\textrm{lk}}(47)$, ${\textrm{lk}}(k)$ and ${\textrm{lk}}(j)$ we see that ${\textrm{lk}}(15)$ can not be completed.
Also, $(8, 9)\cong (9, 8)$ by the map (0, 17)(1, 18)(2, 21)(3, 22)(4, 23)(5, 24)(6, 25)(7, 26)(8, 27)(9, 28)(10, 29)(11, 30)(12, 31)(13, 32) (14, 33)(15, 16)(34, 34)(36, 47)(37, 46)(38, 45)(39, 44)(40, 43)(41, 42). So we have $(k, j) \in \{$(6, 7), (10, 11), (12, 13), (13, 12)$\}$.

If $(k, j) = (12, 13)$ then successively considering ${\textrm{lk}}(12)$, ${\textrm{lk}}(13)$, ${\textrm{lk}}(22)$, ${\textrm{lk}}(21)$, ${\textrm{lk}}(34)$, ${\textrm{lk}}(47)$, ${\textrm{lk}}(23)$, ${\textrm{lk}}(24)$, ${\textrm{lk}}(33)$, ${\textrm{lk}}(16)$, ${\textrm{lk}}(15)$, ${\textrm{lk}}(14)$, ${\textrm{lk}}(4)$, ${\textrm{lk}}(5)$, ${\textrm{lk}}(31)$, ${\textrm{lk}}(32)$, ${\textrm{lk}}(25)$ and ${\textrm{lk}}(26)$ we see that ${\textrm{lk}}(11)$ can not be completed.

\smallskip

\noindent{\bf Subcase 1.1\,:} If $(k, j) = (6, 7)$ then successively we get
${\textrm{lk}}(6) = C_{20}([{\textbf{7}}, 8, 9, 10, 11, 12, 13$, $14, 15, 0, 1, 2, 3, 4, {\textbf{5}}], 46, {\textbf{47}}, 34, 21, 22)$, ${\textrm{lk}}(7) = C_{20}([{\textbf{6}}, 5, 4, 3, 2, 1, 0, 15, 14, 13, 12, 11, 10, 9$, ${\textbf{8}}],23, {\textbf{22}}, 21, 34, 47)$, ${\textrm{lk}}(21) = C_{20}([{\textbf{22}}, 23, 24, 25, 26, 27, 28, 29, 30, 31, 32, 33, 16, 17, {\textbf{18}}], 19$, ${\textbf{34}}, 47, 6, 7)$, ${\textrm{lk}}(22) = C_{20}([{\textbf{21}}, 18, 17, 16, 33, 32, 31, 30, 29, 28, 27, 26, 25, 24, {\textbf{23}}], 8, {\textbf{7}}, 6, 47, 34)$, ${\textrm{lk}}$ $(34) = C_{20}([{\textbf{47}}, 46, 45, 44, 43, 42, 41, 40, 39, 38, 37, 36, 35, 20, {\textbf{19}}], 18, {\textbf{21}}, 22, 7, 6)$, ${\textrm{lk}}(47) = C_{20}$ $([{\textbf{34}}, 19, 20, 35, 36, 37, 38, 39, 40, 41, 42, 43, 44, 45, {\textbf{46}}], 5, {\textbf{6}}, 7, 22, 21)$, ${\textrm{lk}}(4) = C_{20}([{\textbf{5}}, 6, 7$, $8, 9, 10, 11, 12, 13, 14, 15, 0, 1, 2, {\textbf{3}}], 27, {\textbf{26}}, 25, 45, 46)$, ${\textrm{lk}}(5) = C_{20}([{\textbf{4}}, 3, 2, 1, 0, 15, 14, 13, 12$, $11, 10, 9, 8, 7, {\textbf{6}}], 47, {\textbf{46}}, 45, 25, 26)$, ${\textrm{lk}}(45) = C_{20}([{\textbf{46}}, 47, 34, 19, 20, 35, 36, 37, 38, 39, 40, 41$, $42, 43, {\textbf{44}}],24, {\textbf{25}}, 26, 4, 5)$, ${\textrm{lk}}(46) = C_{20}([{\textbf{45}}, 44, 43, 42, 41, 40, 39, 38, 37, 36, 35, 20, 19, 34, {\textbf{47}}]$, $6, {\textbf{5}}, 4, 26, 25)$, ${\textrm{lk}}(25) = C_{20}([{\textbf{26}}, 27, 28, 29, 30, 31, 32, 33, 16, 17, 18, 21, 22, 23, {\textbf{24}}], 44, {\textbf{45}}, 46$, $5, 4)$, ${\textrm{lk}}(26)$ $= C_{20}([{\textbf{25}}, 24, 23, 22, 21, 18, 17, 16, 33, 32, 31, 30, 29, 28, {\textbf{27}}], 3, {\textbf{4}}, 5, 46, 45)$, ${\textrm{lk}}(8) = C_{20}([{\textbf{9}}, 10, 11, 12, 13, 14, 15, 0, 1, 2, 3, 4, 5, 6, {\textbf{7}}], 22, {\textbf{23}}, 24, 44, 43)$, ${\textrm{lk}}(9) = C_{20}([{\textbf{8}}, 7, 6, 5, 4, 3$, $2, 1, 0, 15, 14, 13, 12, 11, {\textbf{10}}], 42, {\textbf{43}}, 44, 24, 23)$, ${\textrm{lk}}(23) = C_{20}([{\textbf{24}}, 25, 26, 27, 28, 29, 30, 31, 32$, $33, 16, 17, 18, 21, {\textbf{22}}], 7, {\textbf{8}}, 9, 43, 44)$, ${\textrm{lk}}(24) = C_{20}([{\textbf{23}}, 22, 21, 18, 17, 16, 33, 32, 31, 30, 29, 28$, $27, 26, {\textbf{25}}], 45, {\textbf{44}}, 43, 9, 8)$, ${\textrm{lk}}(44) = C_{20}([{\textbf{43}}, 42, 41, 40, 39, 38, 37, 36, 35, 20, 19, 34, 47, 46, {\textbf{45}}]$, $25, {\textbf{24}}, 23, 8, 9)$ and ${\textrm{lk}}(15) = C_{20}([{\textbf{14}}, 13, 12, 11, 10, 9, 8, 7, 6, 5, 4, 3, 2, 1, {\textbf{0}}], 17, {\textbf{16}}, 33, m, n)$ for some $m, n \in V$. Observe that $(m, n) \in \{$(39, 40), (40, 39), (41, 42)$\}$. In case $(m, n) = (39, 40)$, completing ${\textrm{lk}}(14)$, ${\textrm{lk}}(15)$, ${\textrm{lk}}(16)$, ${\textrm{lk}}(33)$, ${\textrm{lk}}(39)$ and ${\textrm{lk}}(40)$ it is easy to see that ${\textrm{lk}}(30)$ and ${\textrm{lk}}(31)$ can not be completed. Also in case $(m, n) = (41, 42)$, completing ${\textrm{lk}}(14)$, ${\textrm{lk}}(15)$, ${\textrm{lk}}(16)$, ${\textrm{lk}}(33)$,
${\textrm{lk}}(41)$, ${\textrm{lk}}(42)$, ${\textrm{lk}}(12)$, ${\textrm{lk}}(13)$, ${\textrm{lk}}(24)$,
${\textrm{lk}}(43)$, ${\textrm{lk}}(44)$ we see that ${\textrm{lk}}(22)$ and ${\textrm{lk}}(23)$ can not be completed. So $(m, n) = (40, 39)$ then completing successively we get ${\textrm{lk}}(15) = C_{20}([{\textbf{14}}, 13, 12, 11, 10, 9, 8, 7, 6, 5, 4, 3, 2, 1, {\textbf{0}}], 17, {\textbf{16}}, 33, 40, 39)$, ${\textrm{lk}}(16) = C_{20}([{\textbf{33}}, 32, 31, 30, 29, 28, 27, 26, 25, 24, 23, 22, 21, 18, {\textbf{17}}], 0, {\textbf{15}}, 14, 39, 40)$, ${\textrm{lk}}(14) = C_{20}([{\textbf{15}}$, $0, 1, 2, 3, 4, 5, 6, 7, 8, 9, 10, 11, 12, {\textbf{13}}], 38, {\textbf{39}}, 40, 33, 16)$, ${\textrm{lk}}(33) = C_{20}([{\textbf{16}}, 17, 18, 21, 22, 23, 24$, $25, 26, 27, 28, 29, 30, 31, {\textbf{32}}], 41, {\textbf{40}}, 39, 14, 15)$, ${\textrm{lk}}(39) = C_{20}([{\textbf{40}}, 41, 42, 43, 44, 45, 46, 47, 34$, $19, 20, 35, 36, 37, {\textbf{38}}],13, {\textbf{14}}, 15, 16, 33)$, ${\textrm{lk}}(40) = C_{20}([{\textbf{39}}, 38, 37, 36, 35, 20, 19, 34, 47, 46, 45$, $44, 43, 42, {\textbf{41}}], 32, {\textbf{33}}, 16, 15, 14)$, ${\textrm{lk}}(10) = C_{20}([{\textbf{11}}, 12, 13, 14, 15, 0, 1, 2, 3, 4, 5, 6, 7, 8, {\textbf{9}}], 43$, ${\textbf{42}}, 41, 32, 31)$, ${\textrm{lk}}(11) = C_{20}([{\textbf{10}}, 9, 8, 7, 6, 5, 4, 3, 2, 1, 0, 15, 14, 13, {\textbf{12}}], 30, {\textbf{31}}, 32, 41, 42)$, ${\textrm{lk}}$ $(31) = C_{20}([{\textbf{32}}, 33, 16, 17, 18, 21, 22, 23, 24, 25, 26, 27, 28, 29, {\textbf{30}}], 12, {\textbf{11}}, 10, 42, 41)$, ${\textrm{lk}}(32) = C_{20}([{\textbf{31}}, 30, 29, 28, 27, 26, 25, 24, 23, 22, 21, 18, 17, 16, {\textbf{33}}], 40, {\textbf{41}}, 42, 10, 11)$, ${\textrm{lk}}(41) = C_{20}([{\textbf{42}}$, $43, 44, 45, 46, 47, 34, 19, 20, 35, 36, 37, 38, 39, {\textbf{40}}], 33, {\textbf{32}}, 31, 11, 10)$, ${\textrm{lk}}(42) = C_{20}([{\textbf{41}}, 40, 39$, $38, 37, 36, 35, 20, 19, 34, 47, 46, 45, 44, {\textbf{43}}], 9, {\textbf{10}}, 11, 31, 32)$, ${\textrm{lk}}(12) = C_{20}([{\textbf{13}}, 14, 15, 0, 1, 2, 3$, $4, 5, 6, 7, 8, 9, 10, {\textbf{11}}], 31$, ${\textbf{30}}, 29, 37, 38)$, ${\textrm{lk}}(13) = C_{20}([{\textbf{12}}, 11, 10, 9, 8, 7, 6, 5, 4, 3, 2, 1, 0, 15$, ${\textbf{14}}], 39, {\textbf{38}}, 37, 29, 30)$, ${\textrm{lk}}(29) = C_{20}([{\textbf{30}}, 31, 32, 33, 16, 17, 18, 21, 22, 23, 24, 25, 26, 27, {\textbf{28}}], 36$, ${\textbf{37}}, 38, 13, 12)$, ${\textrm{lk}}(30) = C_{20}$ $([{\textbf{29}}, 28, 27, 26, 25, 24, 23, 22, 21, 18, 17, 16, 33, 32, {\textbf{31}}], 11, {\textbf{12}}, 13$, $38, 37)$, ${\textrm{lk}}(37) = C_{20}([{\textbf{38}}, 39, 40, 41, 42, 43, 44, 45, 46, 47, 34, 19, 20, 35, {\textbf{36}}], 28, {\textbf{29}}, 30, 12, 13)$, ${\textrm{lk}}(38) = C_{20}([{\textbf{37}}, 36, 35, 20, 19, 34, 47, 46, 45, 44, 43, 42, 41, 40, {\textbf{39}}], 14, {\textbf{13}}, 12, 30, 29)$, ${\textrm{lk}}(43) = C_{20}([{\textbf{44}}, 45, 46, 47, 34, 19, 20, 35, 36, 37, 38, 39, 40, 41, {\textbf{42}}], 10, {\textbf{9}}, 8, 23, 24)$. This is isomorphic to $M_1(4, 6, 16)$, as given in Section \ref{example 4.1}, by the map (0, 9)(1, 8)(2, 7)(3, 6)(4, 5)(10, 15)(11, 14)(12, 13)(16, 39, 31, 42)(17, 38, 30, 43)(18, 37, 29, 44)(19, 25, 34, 24)(20, 26, 47, 23)(21, 36, 28, 45)(22, 35, 27, 46)(32, 41, 33, 40).

\smallskip

\noindent{\bf Subcase 1.2\,:} If $(k, j) = (10, 11)$ then constructing successively we get ${\textrm{lk}}(10) = C_{20}([{\textbf{11}}$, $12, 13, 14, 15, 0, 1, 2, 3, 4, 5, 6, 7, 8, {\textbf{9}}], 46, {\textbf{47}}, 34, 21, 22)$,
${\textrm{lk}}(11) = C_{20}([{\textbf{10}}, 9, 8, 7, 6, 5, 4, 3, 2, 1$, $0, 15, 14, 13, {\textbf{12}}], 23, {\textbf{22}}, 21, 34, 47)$, ${\textrm{lk}}(21) = C_{20}([{\textbf{22}}, 23, 24, 25, 26, 27, 28, 29, 30, 31, 32, 33$, $16, 17, {\textbf{18}}], 19, {\textbf{34}}, 47, 10, 11)$, ${\textrm{lk}}(22) = C_{20}([{\textbf{21}}, 18, 17, 16, 33, 32, 31, 30, 29, 28, 27, 26, 25, 24$, ${\textbf{23}}], 12, {\textbf{11}}, 10, 47, 34)$, ${\textrm{lk}}(34) = C_{20}([{\textbf{47}}, 46, 45, 44, 43, 42, 41, 40, 39, 38, 37, 36, 35, 20, {\textbf{19}}], 18$, ${\textbf{21}}, 22, 11, 10)$, ${\textrm{lk}}(47) = C_{20}([{\textbf{34}}, 19, 20, 35, 36, 37, 38, 39, 40, 41, 42, 43, 44, 45, {\textbf{46}}], 9, {\textbf{10}}, 11$, $22, 21)$ and ${\textrm{lk}}(15)$ $= C_{20}([{\textbf{14}}, 13, 12, 11, 10, 9, 8, 7, 6, 5, 4, 3, 2, 1, {\textbf{0}}], 17, {\textbf{16}}, 33, m, n)$ for some $m, n \in V$. In this case we see that $(m, n) \in \{$(39, 40), (40, 39), (41, 42), (42, 41), (43, 44), (44, 43)$\}$. If $(m, n) = (39, 40)$ then successively considering ${\textrm{lk}}(14)$, ${\textrm{lk}}(15)$, ${\textrm{lk}}(16)$, ${\textrm{lk}}(33)$, ${\textrm{lk}}(39)$, ${\textrm{lk}}(40)$, ${\textrm{lk}}(13)$, ${\textrm{lk}}(29)$, ${\textrm{lk}}(30)$, ${\textrm{lk}}(37)$ and ${\textrm{lk}}(38)$ we see that ${\textrm{lk}}(11)$ and ${\textrm{lk}}(12)$ can not be completed. Proceeding similarly for $(m, n) \in \{$(40, 39), (41, 42), (42, 41), (43, 44)$\}$ it is easy to see that the map does not exist. If $(m, n) = (44, 43)$ then completing successively we get ${\textrm{lk}}(15) = C_{20}([{\textbf{14}}, 13, 12, 11, 10, 9, 8, 7, 6, 5, 4, 3, 2, 1, {\textbf{0}}], 17, {\textbf{16}}, 33, 44, 43)$, ${\textrm{lk}}(16) = C_{20}([{\textbf{33}}, 32, 31, 30, 29, 28, 27, 26, 25, 24, 23, 22, 21, 18, {\textbf{17}}], 0, {\textbf{15}}, 14, 43, 44)$, ${\textrm{lk}}(14) = C_{20}([{\textbf{15}}$, $0, 1, 2, 3, 4, 5, 6, 7, 8, 9, 10, 11, 12, {\textbf{13}}], 42, {\textbf{43}}, 44, 33, 16)$, ${\textrm{lk}}(33) = C_{20}([{\textbf{16}}, 17, 18, 21, 22, 23, 24$, $25, 26, 27, 28, 29, 30, 31, {\textbf{32}}], 45, {\textbf{44}}, 43, 14, 15)$, ${\textrm{lk}}(43) = C_{20}([{\textbf{44}}, 45, 46, 47, 34, 19, 20, 35, 36$, $37, 38, 39, 40, 41, {\textbf{42}}], 13, {\textbf{14}}, 15, 16, 33)$, ${\textrm{lk}}(44) = C_{20}([{\textbf{43}}, 42, 41, 40, 39, 38, 37, 36, 35, 20, 19$, $34, 47, 46, {\textbf{45}}], 32, {\textbf{33}}, 16, 15, 14)$, ${\textrm{lk}}(4) = C_{20}([{\textbf{5}}, 6, 7, 8, 9, 10, 11, 12, 13, 14, 15, 0, 1, 2, {\textbf{3}}], 27$, ${\textbf{26}}, 25, 40, 39)$, ${\textrm{lk}}(5) = C_{20}([{\textbf{4}}, 3, 2, 1, 0, 15, 14, 13, 12, 11, 10, 9, 8, 7, {\textbf{6}}], 38, {\textbf{39}}, 40, 25, 26)$, ${\textrm{lk}}$ $(25) = C_{20}$ $([{\textbf{26}}, 27, 28, 29, 30, 31, 32, 33, 16, 17, 18, 21, 22, 23, {\textbf{24}}], 41, {\textbf{40}}, 39, 5, 4)$, ${\textrm{lk}}(26) = C_{20}([{\textbf{25}}, 24, 23, 22, 21, 18, 17, 16, 33, 32, 31, 30, 29, 28, {\textbf{27}}], 3, {\textbf{4}}, 5, 39, 40)$, ${\textrm{lk}}(39) = C_{20}([{\textbf{40}}$, $41, 42, 43, 44, 45, 46, 47, 34, 19, 20, 35, 36, 37, {\textbf{38}}], 6, {\textbf{5}}, 4, 26, 25)$, ${\textrm{lk}}(40) = C_{20}([{\textbf{39}}, 38, 37, 36$, $35, 20, 19, 34, 47, 46, 45, 44, 43, 42, {\textbf{41}}], 24, {\textbf{25}}, 26, 4, 5)$, ${\textrm{lk}}(6) = C_{20}([{\textbf{7}}, 8, 9, 10, 11, 12, 13, 14$, $15, 0, 1, 2, 3, 4, {\textbf{5}}], 39, {\textbf{38}}, 37, 29, 30)$, ${\textrm{lk}}(7) = C_{20}([{\textbf{6}}, 5, 4, 3, 2, 1, 0, 15, 14, 13, 12, 11, 10, 9, {\textbf{8}}]$, $31, {\textbf{30}}, 29, 37, 38)$, ${\textrm{lk}}(29) = C_{20}([{\textbf{30}}, 31, 32, 33, 16, 17, 18, 21, 22, 23, 24, 25, 26, 27, {\textbf{28}}], 36, {\textbf{37}}$, $38, 6, 7)$, ${\textrm{lk}}(30) = C_{20}([{\textbf{29}}, 28, 27, 26, 25, 24, 23, 22, 21, 18, 17, 16, 33, 32, {\textbf{31}}], 8, {\textbf{7}}, 6, 38, 37)$, ${\textrm{lk}}$ $(37) = C_{20}([{\textbf{38}}, 39, 40, 41, 42, 43, 44, 45, 46, 47, 34, 19, 20, 35, {\textbf{36}}], 28, {\textbf{29}}, 30, 7, 6)$, ${\textrm{lk}}(38) = C_{20}([{\textbf{37}}, 36, 35, 20, 19, 34, 47, 46, 45, 44, 43, 42, 41, 40, {\textbf{39}}], 5, {\textbf{6}}, 7, 30, 29)$, ${\textrm{lk}}(8) = C_{20}([{\textbf{9}}, 10$, $11, 12, 13, 14, 15, 0, 1, 2, 3, 4, 5, 6, {\textbf{7}}], 30$, ${\textbf{31}}, 32, 45, 46)$, ${\textrm{lk}}(9) = C_{20}([{\textbf{8}}, 7, 6, 5, 4, 3, 2, 1, 0, 15$, $14, 13, 12, 11, {\textbf{10}}], 47, {\textbf{46}}, 45, 32, 31)$, ${\textrm{lk}}(31) = C_{20}([{\textbf{32}}, 33, 16, 17, 18, 21, 22, 23, 24, 25, 26, 27$, $28, 29, {\textbf{30}}], 7, {\textbf{8}}, 9, 46, 45)$, ${\textrm{lk}}(32) = C_{20}([{\textbf{31}}, 30, 29, 28, 27, 26, 25, 24, 23, 22, 21, 18, 17, 16, {\textbf{33}}]$, $44, {\textbf{45}}, 46, 9, 8)$, ${\textrm{lk}}(45) =  C_{20}([{\textbf{46}}, 47, 34, 19, 20, 35, 36, 37, 38, 39, 40, 41, 42, 43, {\textbf{44}}], 33, {\textbf{32}}, 31$, $8, 9)$, ${\textrm{lk}}(46)= C_{20}([{\textbf{45}}, 44, 43, 42, 41, 40, 39, 38, 37, 36, 35, 20, 19, 34, {\textbf{47}}], 10, {\textbf{9}}, 8, 31, 32)$, ${\textrm{lk}}$ $(12) = C_{20}([{\textbf{13}}, 14, 15, 0, 1, 2, 3, 4, 5, 6, 7, 8, 9, 10, {\textbf{11}}], 22, {\textbf{23}}, 24, 41, 42)$, ${\textrm{lk}}(13) = C_{20}([{\textbf{12}}, 11$, $10, 9, 8, 7, 6, 5, 4, 3, 2, 1, 0, 15, {\textbf{14}}], 43, {\textbf{42}}, 41, 24, 23)$, ${\textrm{lk}}(23) = C_{20}([{\textbf{24}}, 25, 26, 27, 28, 29, 30, 31$, $32, 33, 16, 17, 18, 21, {\textbf{22}}], 11, {\textbf{12}}, 13, 42, 41)$, ${\textrm{lk}}$ $(24)= C_{20}([{\textbf{23}}, 22, 21, 18, 17, 16, 33, 32, 31, 30$, $29, 28, 27, 26, {\textbf{25}}], 40, {\textbf{41}}, 42, 13, 12)$, ${\textrm{lk}}(41) = C_{20}$ $([{\textbf{42}}, 43, 44, 45, 46, 47, 34, 19, 20, 35, 36, 37$, $38, 39, {\textbf{40}}], 25, {\textbf{24}}, 23, 12, 13)$, ${\textrm{lk}}(42) = C_{20}([{\textbf{41}}, 40, 39, 38, 37, 36, 35, 20, 19, 34, 47, 46, 45, 44$, ${\textbf{43}}], 14, {\textbf{13}}, 12, 23, 24)$. This map is isomorphic to $M_2(4, 6, 16)$ as given in Section \ref{example 4.1}, by the map (0, 7)(1, 6)(2, 5)(3, 4)(8, 15)(9, 14)(10, 13)(11, 12)(16, 31) (17, 30)(18, 29)(19, 45, 41, 37)(20, 46, 42, 38)(21, 28)(22, 27)(23, 26)(24, 25)(32, 33)(34, 44, 40, 36)(35, 47, 43, 39).
\smallskip

\noindent{\bf Subcase 1.3\,:} If $(k, j) = (13, 12)$ then successively we get ${\textrm{lk}}(12)$ = $C_{20}([{\textbf{13}}$, 14, 15, 0, 1, 2, 3, 4, 5, 6, 7, 8, 9, 10, ${\textbf{11}}$], 23, ${\textbf{22}}$, 21, 34, 47),
${\textrm{lk}}(13)$ = $C_{20}([{\textbf{12}}$, 11, 10, 9, 8, 7, 6, 5, 4, 3, 2, 1, 0, 15, ${\textbf{14}}$],46, ${\textbf{47}}$, 34, 21, 22), ${\textrm{lk}}(21)$ = $C_{20}([{\textbf{22}}$, 23, 24, 25, 26, 27, 28, 29, 30, 31, 32, 33, 16, 17, ${\textbf{18}}$], 19, ${\textbf{34}}$, 47, 13, 12), ${\textrm{lk}}(22)$ = $C_{20}([{\textbf{21}}$, 18, 17, 16, 33, 32, 31, 30, 29, 28, 27, 26, 25, 24, ${\textbf{23}}$], 11, ${\textbf{12}}$, 13, 47, 34), ${\textrm{lk}}(34)$ = $C_{20}([{\textbf{47}}$, 46, 45, 44, 43, 42, 41, 40, 39, 38, 37, 36, 35, 20, ${\textbf{19}}$], 18, ${\textbf{21}}$, 22, 12, 13), ${\textrm{lk}}(47)$ = $C_{20}([{\textbf{34}}$, 19, 20, 35, 36, 37, 38, 39, 40, 41, 42, 43, 44, 45, ${\textbf{46}}$], 14, ${\textbf{13}}$, 12, 22, 21), ${\textrm{lk}}(33)$ = $C_{20}([{\textbf{16}}$, 17, 18, 21, 22, 23, 24, 25, 26, 27, 28, 29, 30, 31, ${\textbf{32}}$], 44, ${\textbf{45}}$, 46, 14, 15), ${\textrm{lk}}(14)$ = $C_{20}([{\textbf{15}}$, 0, 1, 2, 3, 4, 5, 6, 7, 8, 9, 10, 11, 12, ${\textbf{13}}$], 47, ${\textbf{46}}$, 45, 33, 16), ${\textrm{lk}}(15)$ = $C_{20}([{\textbf{14}}$, 13, 12, 11, 10, 9, 8, 7, 6, 5, 4, 3, 2, 1, ${\textbf{0}}$], 17, ${\textbf{16}}$, 33, 45, 46), ${\textrm{lk}}(16)$ = $C_{20}([{\textbf{33}}$, 32, 31, 30, 29, 28, 27, 26, 25, 24, 23, 22, 21, 18, ${\textbf{17}}$], 0, ${\textbf{15}}$, 14, 46, 45), ${\textrm{lk}}(45)$ = $C_{20}([{\textbf{46}}$, 47, 34, 19, 20, 35, 36, 37, 38, 39, 40, 41, 42, 43, ${\textbf{44}}$],32, ${\textbf{33}}$, 16, 15, 14),
${\textrm{lk}}(46)$ = $C_{20}([{\textbf{45}}$, 44, 43, 42, 41, 40, 39, 38, 37, 36, 35, 20, 19, 34, ${\textbf{47}}$], 13, ${\textbf{14}}$, 15, 16, 33) and ${\textrm{lk}}(44)$ = $C_{20}([{\textbf{43}}$, 42, 41, 40, 39, 38, 37, 36, 35, 20, 19, 34, 47, 46, ${\textbf{45}}$], 33, ${\textbf{32}}$, 31, $m$, $n)$ for some $m, n \in V$. In this case, $(m, n) \in \{$(6, 7), (7, 6), (8, 9), (9, 8)$\}$. If $(m, n) = (6, 7)$ then successively considering ${\textrm{lk}}(6)$, ${\textrm{lk}}(7)$, ${\textrm{lk}}(31)$, ${\textrm{lk}}(32)$, ${\textrm{lk}}(43)$,
${\textrm{lk}}(44)$, ${\textrm{lk}}(4)$, ${\textrm{lk}}(5)$, we see 25\,29 as an edge and a non-edge both.
If $(m, n) = (9, 8)$ then successively considering ${\textrm{lk}}(8)$, ${\textrm{lk}}(9)$, ${\textrm{lk}}(31)$, ${\textrm{lk}}(32)$, ${\textrm{lk}}(43)$, ${\textrm{lk}}(44)$, ${\textrm{lk}}(10)$, ${\textrm{lk}}(11)$ we see 24\,29 as an edge and a non-edge both. So we have $(m, n) \in \{$(7, 6), (8, 9)$\}$.

If $(m, n) = (7, 6)$ then completing successively we get ${\textrm{lk}}(44)$ = $C_{20}([{\textbf{43}}$, 42, 41, 40, 39, 38, 37, 36, 35, 20, 19, 34, 47, 46, ${\textbf{45}}$], 33, ${\textbf{32}}$, 31, 7, 6), ${\textrm{lk}}(6)$ = $C_{20}([{\textbf{7}}$, 8, 9, 10, 11, 12, 13, 14, 15, 0, 1, 2, 3, 4, ${\textbf{5}}$], 42, ${\textbf{43}}$, 44, 32, 31), ${\textrm{lk}}(7)$ = $C_{20}([{\textbf{6}}$, 5, 4, 3, 2, 1, 0, 15, 14, 13, 12, 11, 10, 9, ${\textbf{8}}$], 30, ${\textbf{31}}$, 32, 44, 43), ${\textrm{lk}}(31)$ = $C_{20}([{\textbf{32}}$, 33, 16, 17, 18, 21, 22, 23, 24, 25, 26, 27, 28, 29, ${\textbf{30}}$], 8, ${\textbf{7}}$, 6, 43, 44), ${\textrm{lk}}(32)$ = $C_{20}([{\textbf{31}}$, 30, 29, 28, 27, 26, 25, 24, 23, 22, 21, 18, 17, 16, ${\textbf{33}}$], 45, ${\textbf{44}}$, 43, 6, 7), ${\textrm{lk}}(43)$ = $C_{20}([{\textbf{44}}$, 45, 46, 47, 34, 19, 20, 35, 36, 37, 38, 39, 40, 41, ${\textbf{42}}$], 5, ${\textbf{6}}$, 7, 31, 32),
${\textrm{lk}}(4)$ = $C_{20}([{\textbf{5}}$, 6, 7, 8, 9, 10, 11, 12, 13, 14, 15, 0, 1, 2, ${\textbf{3}}$], 27, ${\textbf{26}}$, 25, 41, 42), ${\textrm{lk}}(5)$ = $C_{20}([{\textbf{4}}$, 3, 2, 1, 0, 15, 14, 13, 12, 11, 10, 9, 8, 7, ${\textbf{6}}$], 43, ${\textbf{42}}$, 41, 25, 26), ${\textrm{lk}}(25)$ = $C_{20}([{\textbf{26}}$, 27, 28, 29, 30, 31, 32, 33, 16, 17, 18, 21, 22, 23, ${\textbf{24}}$],40, ${\textbf{41}}$, 42, 5, 4), ${\textrm{lk}}(26)$ = $C_{20}([{\textbf{25}}$, 24, 23, 22, 21, 18, 17, 16, 33, 32, 31, 30, 29, 28, ${\textbf{27}}$], 3, ${\textbf{4}}$, 5, 42, 41), ${\textrm{lk}}(41)$ = $C_{20}([{\textbf{42}}$, 43, 44, 45, 46, 47, 34, 19, 20, 35, 36, 37, 38, 39, ${\textbf{40}}$], 24, ${\textbf{25}}$, 26, 4, 5), ${\textrm{lk}}(42)$ = $C_{20}([{\textbf{41}}$, 40, 39, 38, 37, 36, 35, 20, 19, 34, 47, 46, 45, 44, ${\textbf{43}}$], 6, ${\textbf{5}}$, 4, 26, 25), ${\textrm{lk}}(8)$ = $C_{20}([{\textbf{9}}$, 10, 11, 12, 13, 14, 15, 0, 1, 2, 3, 4, 5, 6, ${\textbf{7}}$], 31, ${\textbf{30}}$, 29, 37, 38), ${\textrm{lk}}(9)$ = $C_{20}([{\textbf{8}}$, 7, 6, 5, 4, 3, 2, 1, 0, 15, 14, 13, 12, 11, ${\textbf{10}}$], 39, ${\textbf{38}}$, 37, 29, 30), ${\textrm{lk}}(29)$ = $C_{20}([{\textbf{30}}$, 31, 32, 33, 16, 17, 18, 21, 22, 23, 24, 25, 26, 27, ${\textbf{28}}$], 36, ${\textbf{37}}$, 38, 9, 8), ${\textrm{lk}}(30)$ = $C_{20}([{\textbf{29}}$, 28, 27, 26, 25, 24, 23, 22, 21, 18, 17, 16, 33, 32, ${\textbf{31}}$], 7, ${\textbf{8}}$, 9, 38, 37), ${\textrm{lk}}(37)$ = $C_{20}([{\textbf{38}}$, 39, 40, 41, 42, 43, 44, 45, 46, 47, 34, 19, 20, 35, ${\textbf{36}}$], 28, ${\textbf{29}}$, 30, 8, 9), ${\textrm{lk}}(38)$ = $C_{20}([{\textbf{37}}$, 36, 35, 20, 19, 34, 47, 46, 45, 44, 43, 42, 41, 40, ${\textbf{39}}$], 10, ${\textbf{9}}$, 8, 30, 29), ${\textrm{lk}}(10)$ = $C_{20}([{\textbf{11}}$, 12, 13, 14, 15, 0, 1, 2, 3, 4, 5, 6, 7, 8, ${\textbf{9}}$], 38, ${\textbf{39}}$, 40, 24, 23),
${\textrm{lk}}(11)$ = $C_{20}([{\textbf{10}}$, 9, 8, 7, 6, 5, 4, 3, 2, 1, 0, 15, 14, 13, ${\textbf{12}}$], 22, ${\textbf{23}}$, 24, 40, 39), ${\textrm{lk}}(23)$ = $C_{20}([{\textbf{24}}$, 25, 26, 27, 28, 29, 30, 31, 32, 33, 16, 17, 18, 21, ${\textbf{22}}$], 12, ${\textbf{11}}$, 10, 39, 40), ${\textrm{lk}}(24)$ = $C_{20}([{\textbf{23}}$, 22, 21, 18, 17, 16, 33, 32, 31, 30, 29, 28, 27, 26, ${\textbf{25}}$], 41, ${\textbf{40}}$, 39, 10, 11), ${\textrm{lk}}(39)$ = $C_{20}([{\textbf{40}}$, 41, 42, 43, 44, 45, 46, 47, 34, 19, 20, 35, 36, 37, ${\textbf{38}}$],9, ${\textbf{10}}$, 11, 23, 24), ${\textrm{lk}}(40)$ = $C_{20}([{\textbf{39}}$, 38, 37, 36, 35, 20, 19, 34, 47, 46, 45, 44, 43, 42, ${\textbf{41}}$], 25, ${\textbf{24}}$, 23, 11, 10). This is isomorphic to $M_2(4, 6, 16)$ by the map (0, 1)(2, 15)(3, 14)(4, 13)(5, 12)(6, 11)(7, 10)(8, 9)(16, 35) (17, 20)(18, 19)(21, 34)(22, 47)(23, 46)(24, 45)(25, 44)(26, 43)(27, 42)(28, 41)(29, 40)(30, 39)(31, 38)(32, 37)(33, 36).

If $(m, n) = (8, 9)$ then completing successively we get ${\textrm{lk}}(44)$ = $C_{20}([{\textbf{43}}$, 42, 41, 40, 39, 38, 37, 36, 35, 20, 19, 34, 47, 46, ${\textbf{45}}$], 33, ${\textbf{32}}$, 31, 8, 9), ${\textrm{lk}}(43)$ = $C_{20}([{\textbf{44}}$, 45, 46, 47, 34, 19, 20, 35, 36, 37, 38, 39, 40, 41, ${\textbf{42}}$], 10, ${\textbf{9}}$, 8, 31, 32), ${\textrm{lk}}(31)$ = $C_{20}([{\textbf{32}}$, 33, 16, 17, 18, 21, 22, 23, 24, 25, 26, 27, 28, 29, ${\textbf{30}}$], 7, ${\textbf{8}}$, 9, 43, 44), ${\textrm{lk}}(32)$ = $C_{20}([{\textbf{31}}$, 30, 29, 28, 27, 26, 25, 24, 23, 22, 21, 18, 17, 16, ${\textbf{33}}$],45, ${\textbf{44}}$, 43, 9, 8), ${\textrm{lk}}(8)$ = $C_{20}([{\textbf{9}}$, 10, 11, 12, 13, 14, 15, 0, 1, 2, 3, 4, 5, 6, ${\textbf{7}}$], 30, ${\textbf{31}}$, 32, 44, 43), ${\textrm{lk}}$ $(9)$ = $C_{20}([{\textbf{8}}$, 7, 6, 5, 4, 3, 2, 1, 0, 15, 14, 13, 12, 11, ${\textbf{10}}$], 42, ${\textbf{43}}$, 44, 32, 31),${\textrm{lk}}(4)$ = $C_{20}([{\textbf{5}}$, 6, 7, 8, 9, 10, 11, 12, 13, 14, 15, 0, 1, 2, ${\textbf{3}}$], 27, ${\textbf{26}}$, 25, 40, 39), ${\textrm{lk}}(5)$ = $C_{20}([{\textbf{4}}$, 3, 2, 1, 0, 15, 14, 13, 12, 11, 10, 9, 8, 7, ${\textbf{6}}$], 38, ${\textbf{39}}$, 40, 25, 26), ${\textrm{lk}}(25)$ = $C_{20}([{\textbf{26}}$, 27, 28, 29, 30, 31, 32, 33, 16, 17, 18, 21, 22, 23, ${\textbf{24}}$], 41, ${\textbf{40}}$, 39, 5, 4), ${\textrm{lk}}(26)$ = $C_{20}([{\textbf{25}}$, 24, 23, 22, 21, 18, 17, 16, 33, 32, 31, 30, 29, 28, ${\textbf{27}}$], 3, ${\textbf{4}}$, 5, 39, 40), ${\textrm{lk}}(39)$ = $C_{20}([{\textbf{40}}$, 41, 42, 43, 44, 45, 46, 47, 34, 19, 20, 35, 36, 37, ${\textbf{38}}$], 6, ${\textbf{5}}$, 4, 26, 25), ${\textrm{lk}}$ $(40)$ = $C_{20}([{\textbf{39}}$, 38, 37, 36, 35, 20, 19, 34, 47, 46, 45, 44, 43, 42, ${\textbf{41}}$], 24, ${\textbf{25}}$, 26, 4, 5), ${\textrm{lk}}(6)$ = $C_{20}$ $([{\textbf{7}}$, 8, 9, 10, 11, 12, 13, 14, 15, 0, 1, 2, 3, 4, ${\textbf{5}}$], 39, ${\textbf{38}}$, 37, 29, 30), ${\textrm{lk}}(7)$ = $C_{20}([{\textbf{6}}$, 5, 4, 3, 2, 1, 0, 15, 14, 13, 12, 11, 10, 9, ${\textbf{8}}$], 31, ${\textbf{30}}$, 29, 37, 38), ${\textrm{lk}}(29)$ = $C_{20}([{\textbf{30}}$, 31, 32, 33, 16, 17, 18, 21, 22, 23, 24, 25, 26, 27, ${\textbf{28}}$], 36, ${\textbf{37}}$, 38, 6, 7), ${\textrm{lk}}(30)$ = $C_{20}([{\textbf{29}}$, 28, 27, 26, 25, 24, 23, 22, 21, 18, 17, 16, 33, 32, ${\textbf{31}}$], 8, ${\textbf{7}}$, 6, 38, 37), ${\textrm{lk}}(37)$ = $C_{20}([{\textbf{38}}$, 39, 40, 41, 42, 43, 44, 45, 46, 47, 34, 19, 20, 35, ${\textbf{36}}$], 28, ${\textbf{29}}$, 30, 7, 6), ${\textrm{lk}}(38)$ = $C_{20}([{\textbf{37}}$, 36, 35, 20, 19, 34, 47, 46, 45, 44, 43, 42, 41, 40, ${\textbf{39}}$], 5, ${\textbf{6}}$, 7, 30, 29), ${\textrm{lk}}$ $(10)$ = $C_{20}([{\textbf{11}}$, 12, 13, 14, 15, 0, 1, 2, 3, 4, 5, 6, 7, 8, ${\textbf{9}}$], 43, ${\textbf{42}}$, 41, 24, 23), ${\textrm{lk}}(11)$ = $C_{20}([{\textbf{10}}$, 9, 8, 7, 6, 5, 4, 3, 2, 1, 0, 15, 14, 13, ${\textbf{12}}$], 22, ${\textbf{23}}$, 24, 41, 42), ${\textrm{lk}}(23)$ = $C_{20}([{\textbf{24}}$, 25, 26, 27, 28, 29, 30, 31, 32, 33, 16, 17, 18, 21, ${\textbf{22}}$], 12, ${\textbf{11}}$, 10, 42, 41), ${\textrm{lk}}(24)$ = $C_{20}([{\textbf{23}}$, 22, 21, 18, 17, 16, 33, 32, 31, 30, 29, 28, 27, 26, ${\textbf{25}}$], 40, ${\textbf{41}}$, 42, 10, 11), ${\textrm{lk}}(41)$ = $C_{20}([{\textbf{42}}$, 43, 44, 45, 46, 47, 34, 19, 20, 35, 36, 37, 38, 39, ${\textbf{40}}$], 25, ${\textbf{24}}$, 23, 11, 10), ${\textrm{lk}}(42)$ = $C_{20}([{\textbf{41}}$, 40, 39, 38, 37, 36, 35, 20, 19, 34, 47, 46, 45, 44, ${\textbf{43}}$], 9, ${\textbf{10}}$, 11, 23, 24). This is isomorphic to $M_1(4, 6, 16)$ by the map (0, 1)(2, 15)(3, 14)(4, 13)(5, 12)(6, 11)(7, 10)(8, 9)(16, 35)(17, 20)(18, 19)(21, 34)(22, 47)(23, 46)(24, 45)(25, 44)(26, 43)(27, 42)(28, 41)(29, 40)(30, 39)(31, 38)(32, 37)(33, 36). This completes the search for
$(c, d) = (27, 28)$.

\noindent{\bf Case 2\,:} If $(c, d) = (30, 29)$ then constructing successively we get ${\textrm{lk}}(2)$ = $C_{20}([{\textbf{3}}$, 4, 5, 6, 7, 8, 9, 10, 11, 12, 13, 14, 15, 0, ${\textbf{1}}$], 20, ${\textbf{35}}$, 36, 29, 30), ${\textrm{lk}}(3)$ = $C_{20}([{\textbf{2}}$, 1, 0, 15, 14, 13, 12, 11, 10, 9, 8, 7, 6, 5, ${\textbf{4}}$], 31, ${\textbf{30}}$, 29, 36, 35), ${\textrm{lk}}(29)$ = $C_{20}([{\textbf{30}}$, 31, 32, 33, 16, 17, 18, 21, 22, 23, 24, 25, 26, 27, ${\textbf{28}}$], 37, ${\textbf{36}}$, 35, 2, 3), ${\textrm{lk}}(30)$ = $C_{20}([{\textbf{29}}$, 28, 27, 26, 25, 24, 23, 22, 21, 18, 17, 16, 33, 32, ${\textbf{31}}$], 4, ${\textbf{3}}$, 2, 35, 36), ${\textrm{lk}}(35)$ = $C_{20}([{\textbf{36}}$, 37, 38, 39, 40, 41, 42, 43, 44, 45, 46, 47, 34, 19, ${\textbf{20}}$], 1, ${\textbf{2}}$, 3, 30, 29), ${\textrm{lk}}(36)$ = $C_{20}([{\textbf{35}}$, 20, 19, 34, 47, 46, 45, 44, 43, 42, 41, 40, 39, 38, ${\textbf{37}}$], 28, ${\textbf{29}}$, 30, 3, 2) and ${\textrm{lk}}(31)$ = $C_{20}([{\textbf{32}}$, 33, 16, 17, 18, 21, 22, 23, 24, 25, 26, 27, 28, 29, ${\textbf{30}}$], 3, ${\textbf{4}}$, 5, $k$, $j$) for some $k, j \in V$. In this case, $(k, j)\in\{(39, 40), (40, 39), (41, 42), (42, 41), (43, 44), (44, 43), (45, 46), (46, 45)\}$. If $(k, j) = (40, 39)$ then successively considering ${\textrm{lk}}(4)$, ${\textrm{lk}}(5)$, ${\textrm{lk}}(31)$, ${\textrm{lk}}(32)$, ${\textrm{lk}}(39)$, ${\textrm{lk}}(40)$, ${\textrm{lk}}(33)$, ${\textrm{lk}}(27)$, ${\textrm{lk}}(28)$, ${\textrm{lk}}(37)$, ${\textrm{lk}}(38)$, it is easy to see that ${\textrm{lk}}(16)$ and ${\textrm{lk}}(17)$
can not be completed. Now, proceeding similarly for $(k, j) \in \{$(41, 42), (44, 43), (45, 46)$\}$, we see that the map does not exist. Also, $(42, 41)\cong(46, 45)$ by the map (0, 5)(1, 4)(2, 3)(6, 15)(7, 14)(8, 13)
(9, 12)(10, 11)(16, 47)(17, 46)(18, 45)(19, 32)(20, 31)(21, 44)(22, 43)(23, 42)(24, 41)(25, 40)(26, 39)(27, 38)(28, 37)(29, 36)(30, 35). So, we search the map for $(k, j) \in \{$(39, 40), (42, 41), (43, 44)$\}$.

\noindent{\bf Subcase 2.1\,:} If $(k, j) = (39, 40)$ then constructing successively we get ${\textrm{lk}}(4)$ = $C_{20}([{\textbf{5}}$, 6, 7, 8, 9, 10, 11, 12, 13, 14, 15, 0, 1, 2, ${\textbf{3}}$], 30, ${\textbf{31}}$, 32, 40, 39), ${\textrm{lk}}(5)$ = $C_{20}([{\textbf{4}}$, 3, 2, 1, 0, 15, 14, 13, 12, 11, 10, 9, 8, 7, ${\textbf{6}}$],38, ${\textbf{39}}$, 40, 32, 31), ${\textrm{lk}}(39)$ = $C_{20}([{\textbf{40}}$, 41, 42, 43, 44, 45, 46, 47, 34, 19, 20, 35, 36, 37, ${\textbf{38}}$], 6, ${\textbf{5}}$, 4, 31, 32), ${\textrm{lk}}(40)$ = $C_{20}([{\textbf{39}}$, 38, 37, 36, 35, 20, 19, 34, 47, 46, 45, 44, 43, 42, ${\textbf{41}}$], 33, ${\textbf{32}}$, 31, 4, 5), ${\textrm{lk}}(31)$ = $C_{20}([{\textbf{32}}$, 33, 16, 17, 18, 21, 22, 23, 24, 25, 26, 27, 28, 29, ${\textbf{30}}$], 3, ${\textbf{4}}$, 5, 30, 40), ${\textrm{lk}}(32)$ = $C_{20}([{\textbf{31}}$, 30, 29, 28, 27, 26, 25, 24, 23, 22, 21, 18, 17, 16, ${\textbf{33}}$], 41, ${\textbf{40}}$, 39, 5, 4), ${\textrm{lk}}(6)$ = $C_{20}([{\textbf{7}}$, 8, 9, 10, 11, 12, 13, 14, 15, 0, 1, 2, 3, 4, ${\textbf{5}}$], 39, ${\textbf{38}}$, 37, 28, 27), ${\textrm{lk}}(7)$ = $C_{20}([{\textbf{6}}$, 5, 4, 3, 2, 1, 0, 15, 14, 13, 12, 11, 10, 9, ${\textbf{8}}$], 26, ${\textbf{27}}$, 28, 37, 38), ${\textrm{lk}}(27)$ = $C_{20}([{\textbf{28}}$, 29, 30, 31, 32, 33, 16, 17, 18, 21, 22, 23, 24, 25, ${\textbf{26}}$], 8, ${\textbf{7}}$, 6, 38, 37), ${\textrm{lk}}(28)$ = $C_{20}([{\textbf{27}}$, 26, 25, 24, 23, 22, 21, 18, 17, 16, 33, 32, 31, 30, ${\textbf{29}}$], 36, ${\textbf{37}}$, 38, 6, 7), ${\textrm{lk}}(37)$ = $C_{20}([{\textbf{38}}$, 39, 40, 41, 42, 43, 44, 45, 46, 47, 34, 19, 20, 35, ${\textbf{36}}$], 29, ${\textbf{28}}$, 27, 7, 6), ${\textrm{lk}}$ $(38)$ = $C_{20}([{\textbf{37}}$, 36, 35, 20, 19, 34, 47, 46, 45, 44, 43, 42, 41, 40, ${\textbf{39}}$], 5, ${\textbf{6}}$, 7, 27, 28), ${\textrm{lk}}(14)$ = $C_{20}([{\textbf{15}}$, 0, 1, 2, 3, 4, 5, 6, 7, 8, 9, 10, 11, 12, ${\textbf{13}}$], 43, ${\textbf{42}}$, 41, 33, 16), ${\textrm{lk}}(15)$ = $C_{20}([{\textbf{14}}$, 13, 12, 11, 10, 9, 8, 7, 6, 5, 4, 3, 2, 1, ${\textbf{0}}$], 17, ${\textbf{16}}$, 33, 41, 42), ${\textrm{lk}}(16)$ = $C_{20}([{\textbf{33}}$, 32, 31, 30, 29, 28, 27, 26, 25, 24, 23, 22, 21, 18, ${\textbf{17}}$], 0, ${\textbf{15}}$, 14, 42, 41), ${\textrm{lk}}(33)$ = $C_{20}([{\textbf{16}}$, 17, 18, 21, 22, 23, 24, 25, 26, 27, 28, 29, 30, 31, ${\textbf{32}}$],40, ${\textbf{41}}$, 42, 14, 15), ${\textrm{lk}}(41)$ = $C_{20}([{\textbf{42}}$, 43, 44, 45, 46, 47, 34, 19, 20, 35, 36, 37, 38, 39, ${\textbf{40}}$], 32, ${\textbf{33}}$, 16, 15, 14), ${\textrm{lk}}(42)$ = $C_{20}([{\textbf{41}}$, 40, 39, 38, 37, 36, 35, 20, 19, 34, 47, 46, 45, 44, ${\textbf{43}}$], 13, ${\textbf{14}}$, 15, 16, 33) and ${\textrm{lk}}(21)$ = $C_{20}([{\textbf{22}}$, 23, 24, 25, 26, 27, 28, 29, 30, 31, 32, 33, 16, 17, ${\textbf{18}}$], 19, ${\textbf{34}}$, 47, $m, n)$ for some $m, n \in V$. In this case, $(m, n)\in\{$(10, 11), (11, 10)$\}$. If $(m, n) = (11, 10)$ then successively considering ${\textrm{lk}}(10)$, ${\textrm{lk}}(11)$, ${\textrm{lk}}(21)$, ${\textrm{lk}}(22)$, ${\textrm{lk}}(34)$, ${\textrm{lk}}(47)$, ${\textrm{lk}}(12)$, ${\textrm{lk}}(13)$, we see that ${\textrm{lk}}(44)$ can not be completed. So $(m, n) = (10, 11)$.

Then completing successively we get ${\textrm{lk}}(21)$ = $C_{20}([{\textbf{22}}$, 23, 24, 25, 26, 27, 28, 29, 30, 31, 32, 33, 16, 17, ${\textbf{18}}$], 19, ${\textbf{34}}$, 47, 10, 11), ${\textrm{lk}}(22)$ = $C_{20}([{\textbf{21}}$, 18, 17, 16, 33, 32, 31, 30, 29, 28, 27, 26, 25, 24, ${\textbf{23}}$], 12, ${\textbf{11}}$, 10, 47, 34), ${\textrm{lk}}(10)$ = $C_{20}([{\textbf{11}}$, 12, 13, 14, 15, 0, 1, 2, 3, 4, 5, 6, 7, 8, ${\textbf{9}}$], 46, ${\textbf{47}}$, 34, 21, 22), ${\textrm{lk}}(11)$ = $C_{20}([{\textbf{10}}$, 9, 8, 7, 6, 5, 4, 3, 2, 1, 0, 15, 14, 13, ${\textbf{12}}$], 23, ${\textbf{22}}$, 21, 34, 47), ${\textrm{lk}}(34)$ = $C_{20}([{\textbf{47}}$, 46, 45, 44, 43, 42, 41, 40, 39, 38, 37, 36, 35, 20, ${\textbf{19}}$], 18, ${\textbf{21}}$, 22, 11, 10), ${\textrm{lk}}(47)$ = $C_{20}([{\textbf{34}}$, 19, 20, 35, 36, 37, 38, 39, 40, 41, 42, 43, 44, 45, ${\textbf{46}}$], 9, ${\textbf{10}}$, 11, 22, 21), ${\textrm{lk}}(8)$ = $C_{20}([{\textbf{9}}$, 10, 11, 12, 13, 14, 15, 0, 1, 2, 3, 4, 5, 6, ${\textbf{7}}$], 27, ${\textbf{26}}$, 25, 45, 46), ${\textrm{lk}}(9)$ = $C_{20}([{\textbf{8}}$, 7, 6, 5, 4, 3, 2, 1, 0, 15, 14, 13, 12, 11, ${\textbf{10}}$], 47, ${\textbf{46}}$, 45, 25, 26), ${\textrm{lk}}(25)$ = $C_{20}([{\textbf{26}}$, 27, 28, 29, 30, 31, 32, 33, 16, 17, 18, 21, 22, 23, ${\textbf{24}}$], 44, ${\textbf{45}}$, 46, 9, 8), ${\textrm{lk}}(26)$ = $C_{20}([{\textbf{25}}$, 24, 23, 22, 21, 18, 17, 16, 33, 32, 31, 30, 29, 28, ${\textbf{27}}$], 7, ${\textbf{8}}$, 9, 46, 45), ${\textrm{lk}}(45)$ = $C_{20}([{\textbf{46}}$, 47, 34, 19, 20, 35, 36, 37, 38, 39, 40, 41, 42, 43, ${\textbf{44}}$], 24, ${\textbf{25}}$, 26, 8, 9), ${\textrm{lk}}(46)= C_{20}([{\textbf{45}}$, 44, 43, 42, 41, 40, 39, 38, 37, 36, 35, 20, 19, 34, ${\textbf{47}}$], 10, ${\textbf{9}}$, 8, 26, 25), ~${\textrm{lk}}(12)$ = $C_{20}([{\textbf{13}}$, 14, 15, 0, 1, 2, 3, 4, 5, 6, 7, 8, 9, 10, ${\textbf{11}}$], 22, ${\textbf{23}}$, 24, 44, 43), ${\textrm{lk}}(13)$ = $C_{20}([{\textbf{12}}$, 11, 10, 9, 8, 7, 6, 5, 4, 3, 2, 1, 0, 15, ${\textbf{14}}$], 42, ${\textbf{43}}$, 44, 24, 23), ${\textrm{lk}}(23)$ = $C_{20}([{\textbf{24}}$, 25, 26, 27, 28, 29, 30, 31, 32, 33, 16, 17, 18, 21, ${\textbf{22}}$], 11, ${\textbf{12}}$, 13, 43, 44), ${\textrm{lk}}$ $(24)$ = $C_{20}([{\textbf{23}}$, 22, 21, 18, 17, 16, 33, 32, 31, 30, 29, 28, 27, 26, ${\textbf{25}}$], 45, ${\textbf{44}}$, 43, 13, 12), ${\textrm{lk}}(43)$ = $C_{20}([{\textbf{44}}$, 45, 46, 47, 34, 19, 20, 35, 36, 37, 38, 39, 40, 41, ${\textbf{42}}$], 14, ${\textbf{13}}$, 12, 23, 24), ${\textrm{lk}}(44)$ = $C_{20}([{\textbf{43}}$, 42, 41, 40, 39, 38, 37, 36, 35, 20, 19, 34, 47, 46, ${\textbf{45}}$], 25, ${\textbf{24}}$, 23, 12, 13). This is isomorphic to $M_1(4, 6, 16)$ by the map (0, 7)(1, 6)(2, 5)(3, 4)(8, 15)(9, 14)(10, 13)(11, 12)(16, 26)(17, 27)(18, 28)(19, 45, 41, 37)(20, 46, 42, 38)(21, 29)(22, 30)(23, 31)(24, 32)(25, 33)(34, 44, 40, 36)(35, 47, 43, 39).

\noindent{\bf Subcase 2.2\,:} If $(k, j) = (42, 41)$ then successively we get ${\textrm{lk}}(4)$ = $C_{20}([{\textbf{5}}$, 6, 7, 8, 9, 10, 11, 12, 13, 14, 15, 0, 1, 2, ${\textbf{3}}$], 30, ${\textbf{31}}$, 32, 41, 42), ${\textrm{lk}}(5)$ = $C_{20}([{\textbf{4}}$, 3, 2, 1, 0, 15, 14, 13, 12, 11, 10, 9, 8, 7, ${\textbf{6}}$], 43, ${\textbf{42}}$, 41, 32, 31), ${\textrm{lk}}(41)$ = $C_{20}([{\textbf{42}}$, 43, 44, 45, 46, 47, 34, 19, 20, 35, 36, 37, 38, 39, ${\textbf{40}}$], 33, ${\textbf{32}}$, 31, 4, 5), ${\textrm{lk}}(42)$ = $C_{20}([{\textbf{41}}$, 40, 39, 38, 37, 36, 35, 20, 19, 34, 47, 46, 45, 44, ${\textbf{43}}$], 6, ${\textbf{5}}$, 4, 31, 32), ${\textrm{lk}}(31)$ = $C_{20}([{\textbf{32}}$, 33, 16, 17, 18, 21, 22, 23, 24, 25, 26, 27, 28, 29, ${\textbf{30}}$], 3, ${\textbf{4}}$, 5, 42, 41), ${\textrm{lk}}(32)$ = $C_{20}([{\textbf{31}}$, 30, 29, 28, 27, 26, 25, 24, 23, 22, 21, 18, 17, 16, ${\textbf{33}}$], 40, ${\textbf{41}}$, 42, 5, 4), ${\textrm{lk}}(13)$ = $C_{20}([{\textbf{12}}$, 11, 10, 9, 8, 7, 6, 5, 4, 3, 2, 1, 0, 15, ${\textbf{14}}$], 39, ${\textbf{38}}$, 37, 28, 27), ${\textrm{lk}}(27)$ = $C_{20}([{\textbf{28}}$, 29, 30, 31, 32, 33, 16, 17, 18, 21, 22, 23, 24, 25, ${\textbf{26}}$], 11, ${\textbf{12}}$, 13, 38, 37), ${\textrm{lk}}(28)$ = $C_{20}([{\textbf{27}}$, 26, 25, 24, 23, 22, 21, 18, 17, 16, 33, 32, 31, 30, ${\textbf{29}}$], 36, ${\textbf{37}}$, 38, 13, 12), ${\textrm{lk}}(37)$ = $C_{20}([{\textbf{38}}$, 39, 40, 41, 42, 43, 44, 45, 46, 47, 34, 19, 20, 35, ${\textbf{36}}$], 29, ${\textbf{28}}$, 27, 12, 13), ${\textrm{lk}}(38)$ = $C_{20}([{\textbf{37}}$, 36, 35, 20, 19, 34, 47, 46, 45, 44, 43, 42, 41, 40, ${\textbf{39}}$],14, ${\textbf{13}}$, 12, 27, 28) and ${\textrm{lk}}(21)$ = $C_{20}([{\textbf{22}}$, 23, 24, 25, 26, 27, 28, 29, 30, 31, 32, 33, 16, 17, ${\textbf{18}}$], 19, ${\textbf{34}}$, 47, $m$, $n$) for some $m, n \in V$. In this case, $(m, n)\in\{$(8, 9), (9, 8)$\}$. If $(m, n) = (8, 9)$ then successively considering ${\textrm{lk}}(8)$, ${\textrm{lk}}(9)$, ${\textrm{lk}}(21)$, ${\textrm{lk}}(22)$, ${\textrm{lk}}(34)$, ${\textrm{lk}}(47)$, we see that ${\textrm{lk}}(24)$ and ${\textrm{lk}}(25)$ can not be completed. If $(m, n) = (9, 8)$ then completing successively we get ${\textrm{lk}}(21)$ = $C_{20}([{\textbf{22}}$, 23, 24, 25, 26, 27, 28, 29, 30, 31, 32, 33, 16, 17, ${\textbf{18}}$], 19, ${\textbf{34}}$, 47, 9, 8), ${\textrm{lk}}(22)$ = $C_{20}([{\textbf{21}}$, 18, 17, 16, 33, 32, 31, 30, 29, 28, 27, 26, 25, 24, ${\textbf{23}}$], 7, ${\textbf{8}}$, 9, 47, 34), ${\textrm{lk}}(8)$ = $C_{20}([{\textbf{9}}$, 10, 11, 12, 13, 14, 15, 0, 1, 2, 3, 4, 5, 6, ${\textbf{7}}$], 23, ${\textbf{22}}$, 21, 34, 47), ${\textrm{lk}}(9)$ = $C_{20}([{\textbf{8}}$, 7, 6, 5, 4, 3, 2, 1, 0, 15, 14, 13, 12, 11, ${\textbf{10}}$], 46, ${\textbf{47}}$, 34, 21, 22), ${\textrm{lk}}(34)$ = $C_{20}([{\textbf{47}}$, 46, 45, 44, 43, 42, 41, 40, 39, 38, 37, 36, 35, 20, ${\textbf{19}}$], 18, ${\textbf{21}}$, 22, 8, 9), ${\textrm{lk}}(47)$ = $C_{20}([{\textbf{34}}$, 19, 20, 35, 36, 37, 38, 39, 40, 41, 42, 43, 44, 45, ${\textbf{46}}$], 10, ${\textbf{9}}$, 8, 22, 21), ${\textrm{lk}}(6)$ = $C_{20}([{\textbf{7}}$, 8, 9, 10, 11, 12, 13, 14, 15, 0, 1, 2, 3, 4, ${\textbf{5}}$], 42, ${\textbf{43}}$, 44, 24, 23), ${\textrm{lk}}(7)$ = $C_{20}([{\textbf{6}}$, 5, 4, 3, 2, 1, 0, 15, 14, 13, 12, 11, 10, 9, ${\textbf{8}}$], 22, ${\textbf{23}}$, 24, 44, 43), ${\textrm{lk}}(23)$ = $C_{20}([{\textbf{24}}$, 25, 26, 27, 28, 29, 30, 31, 32, 33, 16, 17, 18, 21, ${\textbf{22}}$], 8, ${\textbf{7}}$, 6, 43, 44), ${\textrm{lk}}(24)$ = $C_{20}([{\textbf{23}}$, 22, 21, 18, 17, 16, 33, 32, 31, 30, 29, 28, 27, 26, ${\textbf{25}}$], 45, ${\textbf{44}}$, 43, 6, 7), ${\textrm{lk}}(43)$ = $C_{20}([{\textbf{44}}$, 45, 46, 47, 34, 19, 20, 35, 36, 37, 38, 39, 40, 41, ${\textbf{42}}$], 5, ${\textbf{6}}$, 7, 23, 24), ${\textrm{lk}}(44)$ = $C_{20}([{\textbf{43}}$, 42, 41, 40, 39, 38, 37, 36, 35, 20, 19, 34, 47, 46, ${\textbf{45}}$], 25, ${\textbf{24}}$, 23, 7, 6), ${\textrm{lk}}(10)$ = $C_{20}([{\textbf{11}}$, 12, 13, 14, 15, 0, 1, 2, 3, 4, 5, 6, 7, 8, ${\textbf{9}}$], 47, ${\textbf{46}}$, 45, 25, 26), ${\textrm{lk}}(11)$ = $C_{20}([{\textbf{10}}$, 9, 8, 7, 6, 5, 4, 3, 2, 1, 0, 15, 14, 13, ${\textbf{12}}$], 27, ${\textbf{26}}$, 25, 45, 46), ${\textrm{lk}}(25)$ = $C_{20}([{\textbf{26}}$, 27, 28, 29, 30, 31, 32, 33, 16, 17, 18, 21, 22, 23, ${\textbf{24}}$], 44, ${\textbf{45}}$, 46, 10, 11), ${\textrm{lk}}(26)$ = $C_{20}([{\textbf{25}}$, 24, 23, 22, 21, 18, 17, 16, 33, 32, 31, 30, 29, 28, ${\textbf{27}}$], 12, ${\textbf{11}}$, 10, 46, 45), ${\textrm{lk}}(45)$ = $C_{20}([{\textbf{46}}$, 47, 34, 19, 20, 35, 36, 37, 38, 39, 40, 41, 42, 43, ${\textbf{44}}$], 24, ${\textbf{25}}$, 26, 11, 10), ${\textrm{lk}}(46)$ = $C_{20}([{\textbf{45}}$, 44, 43, 42, 41, 40, 39, 38, 37, 36, 35, 20, 19, 34, ${\textbf{47}}$], 9, ${\textbf{10}}$, 11, 26, 25), ${\textrm{lk}}(14)$ = $C_{20}([{\textbf{15}}$, 0, 1, 2, 3, 4, 5, 6, 7, 8, 9, 10, 11, 12, ${\textbf{13}}$], 38, ${\textbf{39}}$, 40, 33, 16), ${\textrm{lk}}(15)$ = $C_{20}([{\textbf{14}}$, 13, 12, 11, 10, 9, 8, 7, 6, 5, 4, 3, 2, 1, ${\textbf{0}}$], 17, ${\textbf{16}}$, 33, 40, 39), ${\textrm{lk}}(16)$ = $C_{20}([{\textbf{33}}$, 32, 31, 30, 29, 28, 27, 26, 25, 24, 23, 22, 21, 18, ${\textbf{17}}$], 0, ${\textbf{15}}$, 14, 39, 40), ${\textrm{lk}}(33)$ = $C_{20}([{\textbf{16}}$, 17, 18, 21, 22, 23, 24, 25, 26, 27, 28, 29, 30, 31, ${\textbf{32}}$], 41, ${\textbf{40}}$, 39, 14, 15), ${\textrm{lk}}(39)$ = $C_{20}([{\textbf{40}}$, 41, 42, 43, 44, 45, 46, 47, 34, 19, 20, 35, 36, 37, ${\textbf{38}}$], 13, ${\textbf{14}}$, 15, 16, 33). This is isomorphic to $M_2(4, 6, 16)$ by the map (0, 11)(1, 10)(2, 9)(3, 8)(4, 7)(5, 6)(12, 15)(13, 14)(16, 27)(17, 26)(18, 25)(19, 37, 41, 45)(20, 38, 42, 46)(21, 24)(22, 23)(28, 33)(29, 32)(30, 31)(34, 36, 40, 44)(35, 39, 43, 47).

\noindent{\bf Subcase 2.3\,:} If $(k, j) = (43, 44)$ then constructing successively we get ${\textrm{lk}}(4)$ = $C_{20}([{\textbf{5}}$, 6, 7, 8, 9, 10, 11, 12, 13, 14, 15, 0, 1, 2, ${\textbf{3}}$], 30, ${\textbf{31}}$, 32, 44, 43), ${\textrm{lk}}(5)$ = $C_{20}([{\textbf{4}}$, 3, 2, 1, 0, 15, 14, 13, 12, 11, 10, 9, 8, 7, ${\textbf{6}}$],42, ${\textbf{43}}$, 44, 32, 31), ${\textrm{lk}}(43)$ = $C_{20}([{\textbf{44}}$, 45, 46, 47, 34, 19, 20, 35, 36, 37, 38, 39, 40, 41, ${\textbf{42}}$], 6, ${\textbf{5}}$, 4, 31, 32), ${\textrm{lk}}(44)$ = $C_{20}([{\textbf{43}}$, 42, 41, 40, 39, 38, 37, 36, 35, 20, 19, 34, 47, 46, ${\textbf{45}}$], 33, ${\textbf{32}}$, 31, 4, 5), ${\textrm{lk}}(31)$ = $C_{20}([{\textbf{32}}$, 33, 16, 17, 18, 21, 22, 23, 24, 25, 26, 27, 28, 29, ${\textbf{30}}$], 3, ${\textbf{4}}$, 5, 43, 44), ${\textrm{lk}}(32)$ = $C_{20}([{\textbf{31}}$, 30, 29, 28, 27, 26, 25, 24, 23, 22, 21, 18, 17, 16, ${\textbf{33}}$], 45, ${\textbf{44}}$, 43, 5, 4), ${\textrm{lk}}(14)$ = $C_{20}([{\textbf{15}}$, 0, 1, 2, 3, 4, 5, 6, 7, 8, 9, 10, 11, 12, ${\textbf{13}}$], 47, ${\textbf{46}}$, 45, 33, 16), ${\textrm{lk}}(15)$ = $C_{20}([{\textbf{14}}$, 13, 12, 11, 10, 9, 8, 7, 6, 5, 4, 3, 2, 1, ${\textbf{0}}$], 17, ${\textbf{16}}$, 33, 45, 46), ${\textrm{lk}}(16)$ = $C_{20}([{\textbf{33}}$, 32, 31, 30, 29, 28, 27, 26, 25, 24, 23, 22, 21, 18, ${\textbf{17}}$], 0, ${\textbf{15}}$, 14, 46, 45), ${\textrm{lk}}(33)$ = $C_{20}([{\textbf{16}}$, 17, 18, 21, 22, 23, 24, 25, 26, 27, 28, 29, 30, 31, ${\textbf{32}}$], 44, ${\textbf{45}}$, 46, 14, 15), ${\textrm{lk}}(45)$ = $C_{20}([{\textbf{46}}$, 47, 34, 19, 20, 35, 36, 37, 38, 39, 40, 41, 42, 43, ${\textbf{44}}$], 32, ${\textbf{33}}$, 16, 15, 14), ${\textrm{lk}}(46)$ = $C_{20}([{\textbf{45}}$, 44, 43, 42, 41, 40, 39, 38, 37, 36, 35, 20, 19, 34, ${\textbf{47}}$], 13, ${\textbf{14}}$, 15, 16, 33), ${\textrm{lk}}(12)$ = $C_{20}([{\textbf{13}}$, 14, 15, 0, 1, 2, 3, 4, 5, 6, 7, 8, 9, 10, ${\textbf{11}}$], 23, ${\textbf{22}}$, 21, 34, 47), ${\textrm{lk}}(13)$ = $C_{20}([{\textbf{12}}$, 11, 10, 9, 8, 7, 6, 5, 4, 3, 2, 1, 0, 15, ${\textbf{14}}$], 46, ${\textbf{47}}$, 34, 21, 22), ${\textrm{lk}}(47)$ = $C_{20}([{\textbf{34}}$, 19, 20, 35, 36, 37, 38, 39, 40, 41, 42, 43, 44, 45, ${\textbf{46}}$], 14, ${\textbf{13}}$, 12, 22, 21), ${\textrm{lk}}(34)$ = $C_{20}([{\textbf{47}}$, 46, 45, 44, 43, 42, 41, 40, 39, 38, 37, 36, 35, 20, ${\textbf{19}}$], 18, ${\textbf{21}}$, 22, 12, 13), ${\textrm{lk}}(21)$ = $C_{20}([{\textbf{22}}$, 23, 24, 25, 26, 27, 28, 29, 30, 31, 32, 33, 16, 17, ${\textbf{18}}$], 19, ${\textbf{34}}$, 47, 13, 12), ${\textrm{lk}}(22)$ = $C_{20}([{\textbf{21}}$, 18, 17, 16, 33, 32, 31, 30, 29, 28, 27, 26, 25, 24, ${\textbf{23}}$], 11, ${\textbf{12}}$, 13, 47, 34) and ${\textrm{lk}}(37)$ = $C_{20}([{\textbf{38}}$, 39, 40, 41, 42, 43, 44, 45, 46, 47, 34, 19, 20, 35, ${\textbf{36}}$], 29, ${\textbf{28}}$, 27, $m$, $n$) for some $m, n \in V$. Observe that $(m, n)\in\{$(8, 9), (9, 8)$\}$. In case $(m, n) = (9, 8)$, completing ${\textrm{lk}}(8)$, ${\textrm{lk}}(9)$,
${\textrm{lk}}(27)$, ${\textrm{lk}}(28)$, ${\textrm{lk}}(37)$, ${\textrm{lk}}(38)$, we see easily that ${\textrm{lk}}(24)$ and ${\textrm{lk}}(25)$ can not be completed. On the other hand when $(m, n) = (8, 9)$ then completing successively we get ${\textrm{lk}}(37)$ = $C_{20}([{\textbf{38}}$, 39, 40, 41, 42, 43, 44, 45, 46, 47, 34, 19, 20, 35, ${\textbf{36}}$], 29, ${\textbf{28}}$, 27, 8, 9), ${\textrm{lk}}(38)$ = $C_{20}([{\textbf{37}}$, 36, 35, 20, 19, 34, 47, 46, 45, 44, 43, 42, 41, 40, ${\textbf{39}}$], 10, ${\textbf{9}}$, 8, 27, 28), ${\textrm{lk}}(8)$ = $C_{20}([{\textbf{9}}$, 10, 11, 12, 13, 14, 15, 0, 1, 2, 3, 4, 5, 6, ${\textbf{7}}$], 26, ${\textbf{27}}$, 28, 37, 38), ${\textrm{lk}}(9)$ = $C_{20}([{\textbf{8}}$, 7, 6, 5, 4, 3, 2, 1, 0, 15, 14, 13, 12, 11, ${\textbf{10}}$], 39, ${\textbf{38}}$, 37, 28, 27), ${\textrm{lk}}(27)$ = $C_{20}([{\textbf{28}}$, 29, 30, 31, 32, 33, 16, 17, 18, 21, 22, 23, 24, 25, ${\textbf{26}}$], 7, ${\textbf{8}}$, 9, 38, 37), ${\textrm{lk}}(28)$ = $C_{20}([{\textbf{27}}$, 26, 25, 24, 23, 22, 21, 18, 17, 16, 33, 32, 31, 30, ${\textbf{29}}$], 36, ${\textbf{37}}$, 38, 9, 8), ${\textrm{lk}}(6)$ = $C_{20}([{\textbf{7}}$, 8, 9, 10, 11, 12, 13, 14, 15, 0, 1, 2, 3, 4, ${\textbf{5}}$], 43, ${\textbf{42}}$, 41, 25, 26), ${\textrm{lk}}(7)$ = $C_{20}([{\textbf{6}}$, 5, 4, 3, 2, 1, 0, 15, 14, 13, 12, 11, 10, 9, ${\textbf{8}}$], 27, ${\textbf{26}}$, 25, 41, 42), ${\textrm{lk}}(41)$ = $C_{20}([{\textbf{42}}$, 43, 44, 45, 46, 47, 34, 19, 20, 35, 36, 37, 38, 39, ${\textbf{40}}$], 24, ${\textbf{25}}$, 26, 7, 6), ${\textrm{lk}}(42)$ = $C_{20}([{\textbf{41}}$, 40, 39, 38, 37, 36, 35, 20, 19, 34, 47, 46, 45, 44, ${\textbf{43}}$], 5, ${\textbf{6}}$, 7, 26, 25), ${\textrm{lk}}(25)$ = $C_{20}([{\textbf{26}}$, 27, 28, 29, 30, 31, 32, 33, 16, 17, 18, 21, 22, 23, ${\textbf{24}}$], 40, ${\textbf{41}}$, 42, 6, 7), ${\textrm{lk}}(26)$ = $C_{20}([{\textbf{25}}$, 24, 23, 22, 21, 18, 17, 16, 33, 32, 31, 30, 29, 28, ${\textbf{27}}$], 8, ${\textbf{7}}$, 6, 42, 41), ${\textrm{lk}}(10)$ = $C_{20}([{\textbf{11}}$, 12, 13, 14, 15, 0, 1, 2, 3, 4, 5, 6, 7, 8, ${\textbf{9}}$], 38, ${\textbf{39}}$, 40, 24, 23), ${\textrm{lk}}(11)$ = $C_{20}([{\textbf{10}}$, 9, 8, 7, 6, 5, 4, 3, 2, 1, 0, 15, 14, 13, ${\textbf{12}}$], 22, ${\textbf{23}}$, 24, 40, 39), ${\textrm{lk}}(23)$ = $C_{20}([{\textbf{24}}$, 25, 26, 27, 28, 29, 30, 31, 32, 33, 16, 17, 18, 21, ${\textbf{22}}$], 12, ${\textbf{11}}$, 10, 39, 40), ${\textrm{lk}}(24)$ = $C_{20}([{\textbf{23}}$, 22, 21, 18, 17, 16, 33, 32, 31, 30, 29, 28, 27, 26, ${\textbf{25}}$], 41, ${\textbf{40}}$, 39, 10, 11), ${\textrm{lk}}(39)$ = $C_{20}([{\textbf{40}}$, 41, 42, 43, 44, 45, 46, 47, 34, 19, 20, 35, 36, 37, ${\textbf{38}}$],9, ${\textbf{10}}$, 11, 23, 24), ${\textrm{lk}}(40)$ = $C_{20}([{\textbf{39}}$, 38, 37, 36, 35, 20, 19, 34, 47, 46, 45, 44, 43, 42, ${\textbf{41}}$], 25, ${\textbf{24}}$, 23, 11, 10). This is isomorphic to $M_1(4, 6, 16)$ by the map (0, 11)(1, 10)(2, 9)(3, 8)(4, 7)(5, 6)(12, 15)(13, 14)(16, 30, 26, 22)(17, 31, 27, 23)(18, 32, 28, 24)(19, 40)(20, 39)(21, 33, 29, 25)(34, 41)(35, 38)(36, 37)(42, 47)(43, 46)(44, 45). This completes the search for $(c, d) = (30, 29)$.

\smallskip

\noindent{\bf Case 3\,:} If $(c, d) = (31, 32)$ then successively we get ${\textrm{lk}}(3)$ = $C_{20}([{\textbf{2}}$, 1, 0, 15, 14, 13, 12, 11, 10, 9, 8, 7, 6, 5, ${\textbf{4}}$], 30, ${\textbf{31}}$, 32, 36, 35), ${\textrm{lk}}(35)$ = $C_{20}([{\textbf{36}}$, 37, 38, 39, 40, 41, 42, 43, 44, 45, 46, 47, 34, 19, ${\textbf{20}}$], 1, ${\textbf{2}}$, 3, 31, 32), ${\textrm{lk}}(36)$ = $C_{20}([{\textbf{35}}$, 20, 19, 34, 47, 46, 45, 44, 43, 42, 41, 40, 39, 38, ${\textbf{37}}$], 33, ${\textbf{32}}$, 31, 3, 2), ${\textrm{lk}}(31)$ = $C_{20}([{\textbf{32}}$, 33, 16, 17, 18, 21, 22, 23, 24, 25, 26, 27, 28, 29, ${\textbf{30}}$], 4, ${\textbf{3}}$, 2, 35, 36), ${\textrm{lk}}(32)$ = $C_{20}([{\textbf{31}}$, 30, 29, 28, 27, 26, 25, 24, 23, 22, 21, 18, 17, 16, ${\textbf{33}}$], 37, ${\textbf{36}}$, 35, 2, 3), ${\textrm{lk}}(14)$ = $C_{20}([{\textbf{15}}$, 0, 1, 2, 3, 4, 5, 6, 7, 8, 9, 10, 11, 12, ${\textbf{13}}$], 39, ${\textbf{38}}$, 37, 33, 16), ${\textrm{lk}}(15)$ = $C_{20}([{\textbf{14}}$, 13, 12, 11, 10, 9, 8, 7, 6, 5, 4, 3, 2, 1, ${\textbf{0}}$], 17, ${\textbf{16}}$, 33, 37, 38), ${\textrm{lk}}(16)$ = $C_{20}([{\textbf{33}}$, 32, 31, 30, 29, 28, 27, 26, 25, 24, 23, 22, 21, 18, ${\textbf{17}}$], 0, ${\textbf{15}}$, 14, 38, 37), ${\textrm{lk}}(33)$ = $C_{20}([{\textbf{16}}$, 17, 18, 21, 22, 23, 24, 25, 26, 27, 28, 29, 30, 31, ${\textbf{32}}$], 36, ${\textbf{37}}$, 38, 14, 15), ${\textrm{lk}}(37)$ = $C_{20}([{\textbf{38}}$, 39, 40, 41, 42, 43, 44, 45, 46, 47, 34, 19, 20, 35, ${\textbf{36}}$], 32, ${\textbf{33}}$, 16, 15, 14), ${\textrm{lk}}(38)$ = $C_{20}([{\textbf{37}}$, 36, 35, 20, 19, 34, 47, 46, 45, 44, 43, 42, 41, 40, ${\textbf{39}}$], 13, ${\textbf{14}}$, 15, 16, 33) and ${\textrm{lk}}(21)$ = $C_{20}([{\textbf{22}}$, 23, 24, 25, 26, 27, 28, 29, 30, 31, 32, 33, 16, 17, ${\textbf{18}}$], 19, ${\textbf{34}}$, 47, $k$, $j$) for some $k, j \in V$. Then we see that $(k, j) \in \{$(6, 7), (8, 9), (9, 8), (10, 11)$\}$. In case $(k, j) = (8, 9)$, considering ${\textrm{lk}}(8)$, ${\textrm{lk}}(9)$, ${\textrm{lk}}(22)$, ${\textrm{lk}}(21)$, ${\textrm{lk}}(34)$ and ${\textrm{lk}}(47)$ successively we see that ${\textrm{lk}}(7)$ can not be completed. Also, $(6, 7)\cong (10, 11)$ by the map (0, 1)(2, 15)(3, 14)(4, 13)(5, 12)(6, 11)(7, 10)(8, 9)(16, 35)(17, 20)(18, 19)(21, 34)(22, 47)(23, 46)(24, 45)(25, 44)(26, 43)(27, 42)(28, 41)(29, 40)(30, 39)(31, 38)(32, 37)(33, 36). So we search for $(k, j)\in \{$(6, 7), (9, 8)$\}$.

\noindent{\bf Subcase 3.1\,:} If $(k, j) = (6, 7)$ then constructing successively we get ${\textrm{lk}}(6)$ = $C_{20}([{\textbf{7}}$, 8, 9, 10, 11, 12, 13, 14, 15, 0, 1, 2, 3, 4, ${\textbf{5}}$], 46, ${\textbf{47}}$, 34, 21, 22), ${\textrm{lk}}(7)$ = $C_{20}([{\textbf{6}}$, 5, 4, 3, 2, 1, 0, 15, 14, 13, 12, 11, 10, 9, ${\textbf{8}}$],23, ${\textbf{22}}$, 21, 34, 47), ${\textrm{lk}}(21)$ = $C_{20}([{\textbf{22}}$, 23, 24, 25, 26, 27, 28, 29, 30, 31, 32, 33, 16, 17, ${\textbf{18}}$], 19, ${\textbf{34}}$, 47, 6, 7), ${\textrm{lk}}(22)$ = $C_{20}([{\textbf{21}}$, 18, 17, 16, 33, 32, 31, 30, 29, 28, 27, 26, 25, 24, ${\textbf{23}}$], 8, ${\textbf{7}}$, 6, 47, 34), ${\textrm{lk}}(34)$ = $C_{20}([{\textbf{47}}$, 46, 45, 44, 43, 42, 41, 40, 39, 38, 37, 36, 35, 20, ${\textbf{19}}$],18, ${\textbf{21}}$, 22, 7, 6), ${\textrm{lk}}(47)$ = $C_{20}([{\textbf{34}}$, 19, 20, 35, 36, 37, 38, 39, 40, 41, 42, 43, 44, 45, ${\textbf{46}}$], 5, ${\textbf{6}}$, 7, 22, 21), ${\textrm{lk}}(4)$ = $C_{20}([{\textbf{5}}$, 6, 7, 8, 9, 10, 11, 12, 13, 14, 15, 0, 1, 2, ${\textbf{3}}$], 31, ${\textbf{30}}$, 29, 45, 46), ${\textrm{lk}}(5)$ = $C_{20}([{\textbf{4}}$, 3, 2, 1, 0, 15, 14, 13, 12, 11, 10, 9, 8, 7, ${\textbf{6}}$], 47, ${\textbf{46}}$, 45, 29, 30),
${\textrm{lk}}(45)$ = $C_{20}([{\textbf{46}}$, 47, 34, 19, 20, 35, 36, 37, 38, 39, 40, 41, 42, 43, ${\textbf{44}}$], 28, ${\textbf{29}}$, 30, 4, 5), ${\textrm{lk}}(46)$ = $C_{20}([{\textbf{45}}$, 44, 43, 42, 41, 40, 39, 38, 37, 36, 35, 20, 19, 34, ${\textbf{47}}$], 6, ${\textbf{5}}$, 4, 30, 29), ${\textrm{lk}}(29)$ = $C_{20}([{\textbf{30}}$, 31, 32, 33, 16, 17, 18, 21, 22, 23, 24, 25, 26, 27, ${\textbf{28}}$], 44, ${\textbf{45}}$, 46, 5, 4), ${\textrm{lk}}(30)$ = $C_{20}([{\textbf{29}}$, 28, 27, 26, 25, 24, 23, 22, 21, 18, 17, 16, 33, 32, ${\textbf{31}}$], 3, ${\textbf{4}}$, 5, 46, 45) and ${\textrm{lk}}(13)$ = $C_{20}([{\textbf{12}}$, 11, 10, 9, 8, 7, 6, 5, 4, 3, 2, 1, 0, 15, ${\textbf{14}}$], 38, ${\textbf{39}}$, 40, $m$, $n$) for some $m, n \in V$. Observe that $(m, n) \in \{$(25, 26), (26, 25)$\}$. In case $(m, n) = (26, 25)$, completing ${\textrm{lk}}(12)$, ${\textrm{lk}}(13)$, ${\textrm{lk}}(25)$, ${\textrm{lk}}(26)$, ${\textrm{lk}}(39)$, ${\textrm{lk}}(40)$, ${\textrm{lk}}(23)$, ${\textrm{lk}}(24)$, it is easy to see that ${\textrm{lk}}(9)$ and ${\textrm{lk}}(10)$ can not be completed. On the other hand when $(m, n) = (25, 26)$ then completing successively we get ${\textrm{lk}}(12)$ = $C_{20}([{\textbf{13}}$, 14, 15, 0, 1, 2, 3, 4, 5, 6, 7, 8, 9, 10, ${\textbf{11}}$], 27, ${\textbf{26}}$, 25, 40, 39), ${\textrm{lk}}(13)$ = $C_{20}([{\textbf{12}}$, 11, 10, 9, 8, 7, 6, 5, 4, 3, 2, 1, 0, 15, ${\textbf{14}}$], 38, ${\textbf{39}}$, 40, 25, 26), ${\textrm{lk}}(25)$ = $C_{20}([{\textbf{26}}$, 27, 28, 29, 30, 31, 32, 33, 16, 17, 18, 21, 22, 23, ${\textbf{24}}$], 41, ${\textbf{40}}$, 39, 13, 12), ${\textrm{lk}}(26)$ = $C_{20}([{\textbf{25}}$, 24, 23, 22, 21, 18, 17, 16, 33, 32, 31, 30, 29, 28, ${\textbf{27}}$], 11, ${\textbf{12}}$, 13, 39, 40), ${\textrm{lk}}(39)$ = $C_{20}([{\textbf{40}}$, 41, 42, 43, 44, 45, 46, 47, 34, 19, 20, 35, 36, 37, ${\textbf{38}}$],14, ${\textbf{13}}$, 12, 26, 25), ${\textrm{lk}}(40)$ = $C_{20}([{\textbf{39}}$, 38, 37, 36, 35, 20, 19, 34, 47, 46, 45, 44, 43, 42, ${\textbf{41}}$], 24, ${\textbf{25}}$, 26, 12, 13), ${\textrm{lk}}(8)$ = $C_{20}([{\textbf{9}}$, 10, 11, 12, 13, 14, 15, 0, 1, 2, 3, 4, 5, 6, ${\textbf{7}}$], 22, ${\textbf{23}}$, 24, 41, 42), ${\textrm{lk}}(9)$ = $C_{20}([{\textbf{8}}$, 7, 6, 5, 4, 3, 2, 1, 0, 15, 14, 13, 12, 11, ${\textbf{10}}$], 43, ${\textbf{42}}$, 41, 24, 23), ${\textrm{lk}}(23)$ = $C_{20}([{\textbf{24}}$, 25, 26, 27, 28, 29, 30, 31, 32, 33, 16, 17, 18, 21, ${\textbf{22}}$], 7, ${\textbf{8}}$, 9, 42, 41), ${\textrm{lk}}(24)$ = $C_{20}([{\textbf{23}}$, 22, 21, 18, 17, 16, 33, 32, 31, 30, 29, 28, 27, 26, ${\textbf{25}}$], 40, ${\textbf{41}}$, 42, 9, 8), ${\textrm{lk}}(41)$ = $C_{20}([{\textbf{42}}$, 43, 44, 45, 46, 47, 34, 19, 20, 35, 36, 37, 38, 39, ${\textbf{40}}$], 25, ${\textbf{24}}$, 23, 8, 9), ${\textrm{lk}}(42)$ = $C_{20}([{\textbf{41}}$, 40, 39, 38, 37, 36, 35, 20, 19, 34, 47, 46, 45, 44, ${\textbf{43}}$], 10, ${\textbf{9}}$, 8, 23, 24), ${\textrm{lk}}(10)$ = $C_{20}([{\textbf{11}}$, 12, 13, 14, 15, 0, 1, 2, 3, 4, 5, 6, 7, 8, ${\textbf{9}}$],42, ${\textbf{43}}$, 44, 28, 27), ${\textrm{lk}}(11)$ = $C_{20}([{\textbf{10}}$, 9, 8, 7, 6, 5, 4, 3, 2, 1, 0, 15, 14, 13, ${\textbf{12}}$], 26, ${\textbf{27}}$, 28, 44, 43), ${\textrm{lk}}(27)$ = $C_{20}([{\textbf{28}}$, 29, 30, 31, 32, 33, 16, 17, 18, 21, 22, 23, 24, 25, ${\textbf{26}}$], 12, ${\textbf{11}}$, 10, 43, 44), ${\textrm{lk}}(28)$ = $C_{20}([{\textbf{27}}$, 26, 25, 24, 23, 22, 21, 18, 17, 16, 33, 32, 31, 30, ${\textbf{29}}$], 45, ${\textbf{44}}$, 43, 10, 11), ${\textrm{lk}}(43)$ = $C_{20}([{\textbf{44}}$, 45, 46, 47, 34, 19, 20, 35, 36, 37, 38, 39, 40, 41, ${\textbf{42}}$], 9, ${\textbf{10}}$, 11, 27, 28), ${\textrm{lk}}(44)$ = $C_{20}([{\textbf{43}}$, 42, 41, 40, 39, 38, 37, 36, 35, 20, 19, 34, 47, 46, ${\textbf{45}}$], 29, ${\textbf{28}}$, 27, 11, 10). This is isomorphic to $M_1(4, 6, 16)$ by the map (0, 13)(1, 12)(2, 11)(3, 10)(4, 9)(5, 8)(6, 7)(14, 15)(16, 42, 22, 46, 26, 20, 30, 38)(17, 43, 23, 47, 27, 35, 31, 39)(18, 44, 24, 34, 28, 36, 32, 40)(19, 29, 37, 33, 41, 21, 45, 25).

\noindent{\bf Subcase 3.2\,:} If $(k, j) = (9, 8)$ then constructing successively we get ${\textrm{lk}}(8)$ = $C_{20}([{\textbf{9}}$, 10, 11, 12, 13, 14, 15, 0, 1, 2, 3, 4, 5, 6, ${\textbf{7}}$], 23, ${\textbf{22}}$, 21, 34, 47), ${\textrm{lk}}(9)$ = $C_{20}([{\textbf{8}}$, 7, 6, 5, 4, 3, 2, 1, 0, 15, 14, 13, 12, 11, ${\textbf{10}}$],46, ${\textbf{47}}$, 34, 21, 22), ${\textrm{lk}}(22)$ = $C_{20}([{\textbf{21}}$, 18, 17, 16, 33, 32, 31, 30, 29, 28, 27, 26, 25, 24, ${\textbf{23}}$], 7, ${\textbf{8}}$, 9, 47, 34), ${\textrm{lk}}(21)$ = $C_{20}([{\textbf{22}}$, 23, 24, 25, 26, 27, 28, 29, 30, 31, 32, 33, 16, 17, ${\textbf{18}}$], 19, ${\textbf{34}}$, 47, 9, 8), ${\textrm{lk}}(34)$ = $C_{20}([{\textbf{47}}$, 46, 45, 44, 43, 42, 41, 40, 39, 38, 37, 36, 35, 20, ${\textbf{19}}$], 18, ${\textbf{21}}$, 22, 8, 9), ${\textrm{lk}}(47)$ = $C_{20}$ $([{\textbf{34}}$, 19, 20, 35, 36, 37, 38, 39, 40, 41, 42, 43, 44, 45, ${\textbf{46}}$], 10, ${\textbf{9}}$, 8, 22, 21) and ${\textrm{lk}}(10)$ = $C_{20}([{\textbf{11}}$, 12, 13, 14, 15, 0, 1, 2, 3, 4, 5, 6, 7, 8, ${\textbf{9}}$], 47, ${\textbf{46}}$, 45, $m$, $n$) for some $m, n \in V$. It is easy to see that $(m, n)\in\{(25, 26), (28, 27)\}$.

If $(m, n) = (25, 26)$ then completing successively we get ${\textrm{lk}}(10)$ = $C_{20}([{\textbf{11}}$, 12, 13, 14, 15, 0, 1, 2, 3, 4, 5, 6, 7, 8, ${\textbf{9}}$], 47, ${\textbf{46}}$, 45, 25, 26), ${\textrm{lk}}(11)$ = $C_{20}([{\textbf{10}}$, 9, 8, 7, 6, 5, 4, 3, 2, 1, 0, 15, 14, 13, ${\textbf{12}}$],27, ${\textbf{26}}$, 25, 45, 46),
${\textrm{lk}}(25)$ = $C_{20}([{\textbf{26}}$, 27, 28, 29, 30, 31, 32, 33, 16, 17, 18, 21, 22, 23, ${\textbf{24}}$], 44, ${\textbf{45}}$, 46, 10, 11), ${\textrm{lk}}(26)$ = $C_{20}([{\textbf{25}}$, 24, 23, 22, 21, 18, 17, 16, 33, 32, 31, 30, 29, 28, ${\textbf{27}}$], 12, ${\textbf{11}}$, 10, 46, 45), ${\textrm{lk}}(45)$ = $C_{20}([{\textbf{46}}$, 47, 34, 19, 20, 35, 36, 37, 38, 39, 40, 41, 42, 43, ${\textbf{44}}$], 24, ${\textbf{25}}$, 26, 11, 10), ${\textrm{lk}}(46)$ = $C_{20}([{\textbf{45}}$, 44, 43, 42, 41, 40, 39, 38, 37, 36, 35, 20, 19, 34, ${\textbf{47}}$], 9, ${\textbf{10}}$, 11, 26, 25), ${\textrm{lk}}(4)$ = $C_{20}$ $([{\textbf{5}}$, 6, 7, 8, 9, 10, 11, 12, 13, 14, 15, 0, 1, 2, ${\textbf{3}}$], 31, ${\textbf{30}}$, 29, 41, 42), ${\textrm{lk}}(5)$ = $C_{20}([{\textbf{4}}$, 3, 2, 1, 0, 15, 14, 13, 12, 11, 10, 9, 8, 7, ${\textbf{6}}$], 43, ${\textbf{42}}$, 41, 29, 30), ${\textrm{lk}}(29)$ = $C_{20}([{\textbf{30}}$, 31, 32, 33, 16, 17, 18, 21, 22, 23, 24, 25, 26, 27, ${\textbf{28}}$], 40, ${\textbf{41}}$, 42, 5, 4), ${\textrm{lk}}(30)$ = $C_{20}([{\textbf{29}}$, 28, 27, 26, 25, 24, 23, 22, 21, 18, 17, 16, 33, 32, ${\textbf{31}}$], 3, ${\textbf{4}}$, 5, 42, 41), ${\textrm{lk}}(41)$ = $C_{20}([{\textbf{42}}$, 43, 44, 45, 46, 47, 34, 19, 20, 35, 36, 37, 38, 39, ${\textbf{40}}$], 28, ${\textbf{29}}$, 30, 4, 5), ${\textrm{lk}}(42)$ = $C_{20}([{\textbf{41}}$, 40, 39, 38, 37, 36, 35, 20, 19, 34, 47, 46, 45, 44, ${\textbf{43}}$], 6, ${\textbf{5}}$, 4, 30, 29), ${\textrm{lk}}(6)$ = $C_{20}([{\textbf{7}}$, 8, 9, 10, 11, 12, 13, 14, 15, 0, 1, 2, 3, 4, ${\textbf{5}}$], 42, ${\textbf{43}}$, 44, 24, 23), ${\textrm{lk}}(7)$ = $C_{20}([{\textbf{6}}$, 5, 4, 3, 2, 1, 0, 15, 14, 13, 12, 11, 10, 9, ${\textbf{8}}$], 22, ${\textbf{23}}$, 24, 44, 43), ${\textrm{lk}}(23)$ = $C_{20}([{\textbf{24}}$, 25, 26, 27, 28, 29, 30, 31, 32, 33, 16, 17, 18, 21, ${\textbf{22}}$], 8, ${\textbf{7}}$, 6, 43, 44), ${\textrm{lk}}(24)$ = $C_{20}([{\textbf{23}}$, 22, 21, 18, 17, 16, 33, 32, 31, 30, 29, 28, 27, 26, ${\textbf{25}}$], 45, ${\textbf{44}}$, 43, 6, 7), ${\textrm{lk}}(43)$ = $C_{20}([{\textbf{44}}$, 45, 46, 47, 34, 19, 20, 35, 36, 37, 38, 39, 40, 41, ${\textbf{42}}$], 5, ${\textbf{6}}$, 7, 23, 24), ${\textrm{lk}}(44)$ = $C_{20}([{\textbf{43}}$, 42, 41, 40, 39, 38, 37, 36, 35, 20, 19, 34, 47, 46, ${\textbf{45}}$], 25, ${\textbf{24}}$, 23, 7, 6), ${\textrm{lk}}(12)$ = $C_{20}([{\textbf{13}}$, 14, 15, 0, 1, 2, 3, 4, 5, 6, 7, 8, 9, 10, ${\textbf{11}}$], 26, ${\textbf{27}}$, 28, 40, 39), ${\textrm{lk}}(13)$ = $C_{20}$ $([{\textbf{12}}$, 11, 10, 9, 8, 7, 6, 5, 4, 3, 2, 1, 0, 15, ${\textbf{14}}$], 38, ${\textbf{39}}$, 40, 28, 27), ${\textrm{lk}}(27)$ = $C_{20}([{\textbf{28}}$, 29, 30, 31, 32, 33, 16, 17, 18, 21, 22, 23, 24, 25, ${\textbf{26}}$], 11, ${\textbf{12}}$, 13, 39, 40), ${\textrm{lk}}(28)$ = $C_{20}([{\textbf{27}}$, 26, 25, 24, 23, 22, 21, 18, 17, 16, 33, 32, 31, 30, ${\textbf{29}}$], 41, ${\textbf{40}}$, 39, 13, 12), ${\textrm{lk}}(39)$ = $C_{20}([{\textbf{40}}$, 41, 42, 43, 44, 45, 46, 47, 34, 19, 20, 35, 36, 37, ${\textbf{38}}$], 14, ${\textbf{13}}$, 12, 27, 28), ${\textrm{lk}}(40)$ = $C_{20}([{\textbf{39}}$, 38, 37, 36, 35, 20, 19, 34, 47, 46, 45, 44, 43, 42, ${\textbf{41}}$], 29, ${\textbf{28}}$, 27, 12, 13). This is isomorphic to $M_1(4, 6, 16)$ by the map (0, 20, 18)(1, 19, 17)(2, 34, 16)(3, 47, 33)(4, 46, 32)(5, 45, 31)(6, 44, 30)(7, 43, 29)(8, 42, 28)(9, 41, 27)(10, 40, 26)(11, 39, 25)(12, 38, 24)(13, 37, 23)(14, 36, 22)(15, 35, 21).

On the other hand when, $(m, n) = (28, 27)$ then completing successively we get ${\textrm{lk}}(10)$ = $C_{20}([{\textbf{11}}$, 12, 13, 14, 15, 0, 1, 2, 3, 4, 5, 6, 7, 8, ${\textbf{9}}$], 47, ${\textbf{46}}$, 45, 28, 27), ${\textrm{lk}}(11)$ = $C_{20}([{\textbf{10}}$, 9, 8, 7, 6, 5, 4, 3, 2, 1, 0, 15, 14, 13, ${\textbf{12}}$], 26, ${\textbf{27}}$, 28, 45, 46), ${\textrm{lk}}(27)$ = $C_{20}([{\textbf{28}}$, 29, 30, 31, 32, 33, 16, 17, 18, 21, 22, 23, 24, 25, ${\textbf{26}}$], 12, ${\textbf{11}}$, 10, 46, 45), ${\textrm{lk}}(28)$ = $C_{20}([{\textbf{27}}$, 26, 25, 24, 23, 22, 21, 18, 17, 16, 33, 32, 31, 30, ${\textbf{29}}$], 44, ${\textbf{45}}$, 46, 10, 11), ${\textrm{lk}}(45)$ = $C_{20}([{\textbf{46}}$, 47, 34, 19, 20, 35, 36, 37, 38, 39, 40, 41, 42, 43, ${\textbf{44}}$], 29, ${\textbf{28}}$, 27, 11, 10), ${\textrm{lk}}(46)$ = $C_{20}([{\textbf{45}}$, 44, 43, 42, 41, 40, 39, 38, 37, 36, 35, 20, 19, 34, ${\textbf{47}}$], 9, ${\textbf{10}}$, 11, 27, 28), ${\textrm{lk}}(4)$ = $C_{20}([{\textbf{5}}$, 6, 7, 8, 9, 10, 11, 12, 13, 14, 15, 0, 1, 2, ${\textbf{3}}$], 31, ${\textbf{30}}$, 29, 44, 43), ${\textrm{lk}}(5)$ $= C_{20}([{\textbf{4}}$, 3, 2, 1, 0, 15, 14, 13, 12, 11, 10, 9, 8, 7, ${\textbf{6}}$], 42, ${\textbf{43}}$, 44, 29, 30), ${\textrm{lk}}(29)$ = $C_{20}([{\textbf{30}}$, 31, 32, 33, 16, 17, 18, 21, 22, 23, 24, 25, 26, 27, ${\textbf{28}}$], 45, ${\textbf{44}}$, 43, 5, 4), ${\textrm{lk}}(30)$ = $C_{20}([{\textbf{29}}$, 28, 27, 26, 25, 24, 23, 22, 21, 18, 17, 16, 33, 32, ${\textbf{31}}$], 3, ${\textbf{4}}$, 5, 43, 44), ${\textrm{lk}}(43)$ = $C_{20}([{\textbf{44}}$, 45, 46, 47, 34, 19, 20, 35, 36, 37, 38, 39, 40, 41, ${\textbf{42}}$], 6, ${\textbf{5}}$, 4, 30, 29), ${\textrm{lk}}(44)$ = $C_{20}([{\textbf{43}}$, 42, 41, 40, 39, 38, 37, 36, 35, 20, 19, 34, 47, 46, ${\textbf{45}}$], 28, ${\textbf{29}}$, 30, 4, 5), ${\textrm{lk}}(6)$ = $C_{20}([{\textbf{7}}$, 8, 9, 10, 11, 12, 13, 14, 15, 0, 1, 2, 3, 4, ${\textbf{5}}$], 43, ${\textbf{42}}$, 41, 24, 23), ${\textrm{lk}}(7)$ = $C_{20}([{\textbf{6}}$, 5, 4, 3, 2, 1, 0, 15, 14, 13, 12, 11, 10, 9, ${\textbf{8}}$], 22, ${\textbf{23}}$, 24, 41, 42), ${\textrm{lk}}(23)$ = $C_{20}([{\textbf{24}}$, 25, 26, 27, 28, 29, 30, 31, 32, 33, 16, 17, 18, 21, ${\textbf{22}}$], 8, ${\textbf{7}}$, 6, 42, 41), ${\textrm{lk}}(24)$ = $C_{20}([{\textbf{23}}$, 22, 21, 18, 17, 16, 33, 32, 31, 30, 29, 28, 27, 26, ${\textbf{25}}$], 40, ${\textbf{41}}$, 42, 6, 7), ${\textrm{lk}}(41)$ = $C_{20}([{\textbf{42}}$, 43, 44, 45, 46, 47, 34, 19, 20, 35, 36, 37, 38, 39, ${\textbf{40}}$], 25, ${\textbf{24}}$, 23, 7, 6), ${\textrm{lk}}(42)$ = $C_{20}([{\textbf{41}}$, 40, 39, 38, 37, 36, 35, 20, 19, 34, 47, 46, 45, 44, ${\textbf{43}}$], 5, ${\textbf{6}}$, 7, 23, 24), ${\textrm{lk}}(12)$ = $C_{20}([{\textbf{13}}$, 14, 15, 0, 1, 2, 3, 4, 5, 6, 7, 8, 9, 10, ${\textbf{11}}$], 27, ${\textbf{26}}$, 25, 40, 39), ${\textrm{lk}}(13)$ = $C_{20}([{\textbf{12}}$, 11, 10, 9, 8, 7, 6, 5, 4, 3, 2, 1, 0, 15, ${\textbf{14}}$], 38, ${\textbf{39}}$, 40, 25, 26), ${\textrm{lk}}(25)$ = $C_{20}$ $([{\textbf{26}}$, 27, 28, 29, 30, 31, 32, 33, 16, 17, 18, 21, 22, 23, ${\textbf{24}}$], 41, ${\textbf{40}}$, 39, 13, 12), ${\textrm{lk}}(26)$ = $C_{20}([{\textbf{25}}$, 24, 23, 22, 21, 18, 17, 16, 33, 32, 31, 30, 29, 28, ${\textbf{27}}$], 11, ${\textbf{12}}$, 13, 39, 40), ${\textrm{lk}}(39)$ = $C_{20}([{\textbf{40}}$, 41, 42, 43, 44, 45, 46, 47, 34, 19, 20, 35, 36, 37, ${\textbf{38}}$],14, ${\textbf{13}}$, 12, 26, 25), ${\textrm{lk}}(40)$ = $C_{20}([{\textbf{39}}$, 38, 37, 36, 35, 20, 19, 34, 47, 46, 45, 44, 43, 42, ${\textbf{41}}$], 24, ${\textbf{25}}$, 26, 12, 13). This is isomorphic to $M_2(4, 6, 16)$ by the map (0, 3)(1, 2)(4, 15)(5, 14)(6, 13)(7, 12)(8, 11)(9, 10)(16, 22, 26, 30)(17, 23, 27, 31)(18, 24, 28, 32)(19, 36)(20, 35)(21, 25, 29, 33)(34, 37)(38, 47)(39, 46)(40, 45)(41, 44)(42, 43). This completes the
search for $(c, d) = (31, 32)$ and thus the Lemma \ref{lem2} is proved. $\hfill\Box$

\smallskip

\noindent {\bf Proof of Lemma \ref{lem3}\,:} Let $N$ be a SEM of type $(6^2, 8)$ on the surface of Euler characteristic $-1$. The notation ${\textrm{lk}}(i)$ = $C_{14}([{\textbf{\textit{i}}}_{\textbf{1}}$, $i_2, i_3, i_4, i_5, i_6$, ${\textbf{\textit{i}}}_{\textbf{7}}]$, $i_8, i_9, i_{10}$, ${\textbf{\textit{i}}}_{\textbf{11}}$, $i_{12}, i_{13}, i_{14})$ for the link of $i$ will mean that $[i, i_7, i_{8}, i_{9}, i_{10}, i_{11}],[i, i_1, i_{14}, i_{13}, i_{12}, i_{11}]$ form hexagonal faces and $[i, i_{1}, i_{2}, i_{3}, i_4, i_5, i_6, i_7]$ forms octagonal face. If $|V|$, $E(N)$, $H(N)$ and $O(N)$ denote number vertices, number of edges, number of hexagonal faces and number of octagonal faces in the map $N$, respectively, then $E(N) = \frac{3|V|}{2}$, $H(N) = \frac{2|V|}{6}$ and $O(N) = \frac{|V|}{8}$. Using Euler's equation we see that if the map exists then $|V| = 24$. Let $V = V(M) = \{$0, 1, \ldots, 23$\}$. Now, we prove the proposition by exhaustive search for all $N$.


Assume that, ${\textrm{lk}}(0)$ = $C_{14}([{\textbf{1}}$, 2, 3, 4, 5, 6, ${\textbf{7}}$], 8, 9, 10, ${\textbf{11}}$, 12, 13, 14). This implies ${\textrm{lk}}(11)$ = $C_{14}([{\textbf{12}}$, 19, 18, 17, 16, 15, ${\textbf{10}}$], 9, 8, 7, ${\textbf{0}}$, 1, 14, 13) and ${\textrm{lk}}(8)$ = $C_{14}([{\textbf{9}}$, $f, e, d, c, b$, ${\textit {\textbf a}}$], $g, h, 6$, ${\textbf 7}$, 0, 11, 10) for some $a, b, c, d, e, f, g, h \in V$. Then we get the partial picture of the map as shown in Figure II. Let $V(O_i)$, for $i = 1, 2, 3$, denote the vertex set of octagonal face $O_i$ then we see that $V(O_1$) = $\{0, 1, 2, 3, 4, 5, 6, 7\}$, $V(O_2)$ = $\{10, 11, 12, 15, 16, 17, 18, 19\}$ and $V(O_3)$ = $\{8, 9, 13, 14, 20, 21, 22, 23\}$. In this case we observe that $a \in \{13, 14, 20\}$. If $a = 14$ then completing successively we get $b = 13$, $c = 20$, $d = 21$, $e = 22$ and $f = 23$. This implies $g = 1$. This contradicts the fact that $g \in V(O_2)$, as $1 \in O_1$. So $a \neq 14$. So, $a = 13$ or $20$.

\noindent{\bf Case 1\,:} If $a = 13$ then successively we get $b = 14$, $c = 20$, $d = 21$, $e = 22$, $f = 23$, $g = 12$ and $h = 19$. This implies ${\textrm{lk}}(8)$ = $C_{14}([{\textbf {9}}$, 23, 22, 21, 20, 14, ${\textbf {13}}$], 12, 19, 6, ${\textbf 7}$, 0, 11, 10), ${\textrm{lk}}(7)$ = $C_{14}([{\textbf {0}}$, 1, 2, 3, 4, 5, ${\textbf {6}}$], 19, 12, 13, ${\textbf 8}$, 9, 10, 11), ${\textrm{lk}}(12)$ = $C_{14}([{\textbf {19}}$, 18, 17, 16, 15, 10, ${\textbf {11}}$], 0, 1, 14, ${\textbf{13}}$, 8, 7, 6), ${\textrm{lk}}(13)$ = $C_{14}([{\textbf {8}}$, 9, 23, 22, 21, 20, ${\textbf {14}}$], 1, 0, 11, ${\textbf 12}$, 19, 6, 7) and ${\textrm{lk}}(19)$ = $C_{14}([{\textbf {18}}$, 17, 16, 15, 10, 11, ${\textbf {12}}$], 13, 8, 7, ${\textbf 6}$, 5, $j$, $i$) for some $i, j \in V(O_3)$. Then we see that $(j, i) \in \{$(20, 21), (21, 20), (21, 22), (23, 22)$\}$. Observe that $(23, 22) \cong (21, 20)$ by the map (0, 11)(1, 10)(2, 15)(3, 16)(4, 17)(5, 18)(6, 19)(7, 12)(8, 13)(9, 14)(20, 23)(21, 22), so we search for $(j, i) \in \{(20, 21), (21, 20), (21, 22)\}$.

\vspace{1.5in}

\begin{center}
\begin{picture}(0, 0)(0, 6)

\unitlength=6mm

\drawpolygon(2, .8)(.8, 2)(-.8, 2)(-2, .8)(-2, -.8)(-.8, -2)(.8, -2)(2, -.8)
\drawpolygon(2, .8)(3, 1.2)(4, .6)(4, -.6)(3, -1.2)(2, -.8)
\drawpolygon(-2, .8)(-3, 1.2)(-4, .6)(-4, -.6)(-3, -1.2)(-2, -.8)

\drawpolygon(.8, 2)(1.2, 3)(.6, 4)(-.6, 4)(-1.2, 3)(-.8, 2)
\drawpolygon(.8, -2)(1.2, -3)(.6, -4)(-.6, -4)(-1.2, -3)(-.8, -2)

\drawpolygon(1.2, 3)(.8, 2)(2, .8)(3, 1.2)(3.1, 2.2)(2.3, 3)
\drawpolygon(-1.2, 3)(-.8, 2)(-2, .8)(-3, 1.2)(-3.1, 2.2)(-2.3, 3)
\drawpolygon(1.2, -3)(.8, -2)(2, -.8)(3, -1.2)(3.1, -2.2)(2.3, -3)
\drawpolygon(-1.2, -3)(-.8, -2)(-2, -.8)(-3, -1.2)(-3.1, -2.2)(-2.3, -3)

\drawpolygon(.6, 4)(1.2, 3)(2.3, 3)(3, 3.5)(3.1, 4.4)(2.5, 5.2)(1.4, 5.3)(.7, 4.8)
\drawpolygon(-.6, 4)(-1.2, 3)(-2.3, 3)(-3, 3.5)(-3.1, 4.4)(-2.5, 5.2)(-1.4, 5.3)(-.7, 4.8)
\drawpolygon(-.6, -4)(-1.2, -3)(-2.3, -3)(-3, -3.5)(-3.1, -4.4)(-2.5, -5.2)(-1.4, -5.3)(-.7, -4.8)
\drawpolygon(.6, -4)(1.2, -3)(2.3, -3)(3, -3.5)(3.1, -4.4)(2.5, -5.2)(1.4, -5.3)(.7, -4.8)

\drawpolygon(5, 2.8)(3.8, 2.9)(3.1, 2.2)(3, 1.2)(4, .6)(4.8, .5)(5.6, 1)(5.8, 1.9)
\drawpolygon(-5, 2.8)(-3.8, 2.9)(-3.1, 2.2)(-3, 1.2)(-4, .6)(-4.8, .5)(-5.6, 1)(-5.8, 1.9)
\drawpolygon(5, -2.8)(3.8, -2.9)(3.1, -2.2)(3, -1.2)(4, -.6)(4.8, -.5)(5.6, -1)(5.8, -1.9)
\drawpolygon(-5, -2.8)(-3.8, -2.9)(-3.1, -2.2)(-3, -1.2)(-4, -.6)(-4.8, -.5)(-5.6, -1)(-5.8, -1.9)

\drawpolygon(-.7, 4.8)(-.6, 4)(.6, 4)(.7, 4.8)(.4, 5.4)(-.4, 5.4)

\drawpolygon(-.7, -4.8)(-.6, -4)(.6, -4)(.7, -4.8)(.4, -5.4)(-.4, -5.4)

\drawpolygon(4.8, .5)(4, .6)(4, -.6)(4.8, -.5)(5.2, -.3)(5.2, .3)
\drawpolygon(-4.8, .5)(-4, .6)(-4, -.6)(-4.8, -.5)(-5.2, -.3)(-5.2, .3)

\drawpolygon(3, 3.5)(2.3, 3)(3.1, 2.2)(3.8, 2.9)(4.1, 3.4)(3.7, 3.8)
\drawpolygon(-3, 3.5)(-2.3, 3)(-3.1, 2.2)(-3.8, 2.9)(-4.1, 3.4)(-3.7, 3.8)
\drawpolygon(3, -3.5)(2.3, -3)(3.1, -2.2)(3.8, -2.9)(4.1, -3.4)(3.7, -3.8)
\drawpolygon(-3, -3.5)(-2.3, -3)(-3.1, -2.2)(-3.8, -2.9)(-4.1, -3.4)(-3.7, -3.8)

\put(1.3, 3.1){\scriptsize 8}
\put(-1.6, 3.1){\scriptsize 11}
\put(-.8, 1.6){\scriptsize 0}
\put(-1.9, .6){\scriptsize 1}
\put(-1.9, -.8){\scriptsize 2}
\put(-.8, -1.8){\scriptsize 3}
\put(.6, -1.8){\scriptsize 4}
\put(1.7, -.8){\scriptsize 5}
\put(1.7, .6){\scriptsize 6}
\put(.6, 1.6){\scriptsize 7}
\put(3.1, 1.2){\scriptsize $h$}
\put(3.2, 2.1){\scriptsize $g$}
\put(2.5, 2.85){\scriptsize $a$}
\put(2.7, 3.5){\scriptsize $b$}
\put(-3.6, 1.9){\scriptsize 13}
\put(-2.9, 2.85){\scriptsize 12}
\put(-3, 3.6){\scriptsize 19}

\put(-3.5, 1.2){\scriptsize 14}

\put(-.6, 4.1){\scriptsize 10}
\put(.3, 4.1){\scriptsize 9}

\put(-.55, 4.7){\scriptsize 15}
\put(.7, 5.1){\scriptsize $f$}

\put(-1.6, 5.4){\scriptsize 16}
\put(1.25, 5.4){\scriptsize $e$}

\put(-2.7, 5.3){\scriptsize 17}
\put(2.45, 5.3){\scriptsize $d$}

\put(-3.6, 4.3){\scriptsize 18}
\put(3.2, 4.3){\scriptsize $c$}

\put(-.2, 0){\scriptsize {\large $O_1$}}

\put(4, 1.5){\scriptsize {\large $O_2$}}
\put(4, -1.5){\scriptsize {\large $O_3$}}
\put(-4.5, 1.5){\scriptsize {\large $O_3$}}
\put(-4.3, -1.5){\scriptsize {\large $O_2$}}

\put(1.5, 4){\scriptsize {\large $O_3$}}
\put(-1.8, 4){\scriptsize {\large $O_2$}}
\put(1.5, -4){\scriptsize {\large $O_2$}}
\put(-1.8, -4){\scriptsize {\large $O_3$}}

\put(-5.5, -7){\mbox{Figure II: Semi-equivelar map N of type $(6^2, 8)$}}

\end{picture}

\end{center}

\vspace{5.5cm}

If $(j, i) = (20, 21)$ then ${\textrm{lk}}(19)$ = $C_{14}([{\textbf {18}}$, 17, 16, 15, 10, 11, ${\textbf {12}}$], 13, 8, 7, ${\textbf 6}$, 5, 20, 21) and ${\textrm{lk}}(5)$ = $C_{14}([{\textbf {4}}$, 3, 2, 1, 0, 7, ${\textbf {6}}$], 19, 18, 21, ${\textbf{20}}$, 14, $k, l)$ for some $k, l \in V(O_2)$. This implies $\deg(14) > 3$, a contradiction. So $(j, i) \neq (20, 21)$.

\noindent {\bf Subcase 2.1\,:} If $(j, i) = (21, 20)$ then ${\textrm{lk}}(19)$ = $C_{14}([{\textbf{18}}$, 17, 16, 15, 10, 11, ${\textbf{12}}$], 13, 8, 7, ${\textbf 6}$, 5, 21, 20), ${\textrm{ lk}}(6)$ = $C_{14}([{\textbf 7}$, 0, 1, 2, 3, 4, ${\textbf 5}$], 21, 20, 18, ${\textbf{19}}$, 12, 13, 8) and ${\textrm{lk}}(5)$ = $C_{14}([{\textbf{4}}$, 3, 2, 1, 0, 7, ${\textbf{6}}$], 19, 18, 20, ${\textbf{21}}$, 22, $k, l)$ for some $k, l \in V(O_2)$. Observe that $(l, k)\in\{$(15, 16), (16, 15), (16, 17), (17, 16)$\}$. If $(l, k) = (17, 16)$ then successively considering ${\textrm{lk}}(5)$ and ${\textrm{lk}}(4)$ we get $\deg(14) > 3$. A contradiction. If $(l, k) = (16, 15)$ then $C_{13}(4, 5, 21, 22, 23, 9, 10, 11, 12, 19, 18, 17, 16)\subseteq{\textrm{lk}}(15)$. A contradiction. If $(l, k) = (16, 17)$ then considering ${\textrm{lk}}(5)$ and we see that ${\textrm{lk}}(15)$ and ${\textrm{lk}}(16)$ can not be completed. If $(l, k) = (15, 16)$ then successively we get ${\textrm{lk}}(5)$ = $C_{14}([{\textbf {4}}$, 3, 2, 1, 0, 7, ${\textbf {6}}$], 19, 18, 20, ${\textbf{21}}$, 22, 16, 15), ${\textrm{lk}}(15)$ = $C_{14}([{\textbf{16}}$, 17, 18, 19, 12, 11, ${\textbf 10}$], 9, 23, 3, ${\textbf 4}$, 5, 21, 22), ${\textrm{lk}}(10)$ = $C_{14}([{\textbf{15}}$, 16, 17, 18, 19, 12, ${\textbf{11}}$], 0, 7, 8, ${\textbf 9}$, 23, 3, 4), ${\textrm{lk}}(9)$ = $C_{14}([{\textbf 8}$, 13, 14, 20, 21, 22, ${\textbf{23}}$], 3, 4, 15, ${\textbf{10}}$, 11, 0, 7) and ${\textrm{lk}}(23)$ = $C_{14}([{\textbf{22}}$, 21, 20, 14, 13, 8, ${\textbf 9}$], 10, 15, 4, ${\textbf 3}$, 2, $n, m)$ for some $m, n \in V(O_2)$. Observe that $n = 17$ and $m = 16$. This implies ${\textrm{lk}}(23)$ = $C_{14}([{\textbf{22}}$, 21, 20, 14, 13, 8, ${\textbf 9}$], 10, 15, 4, ${\textbf 3}$, 2, 17, 16), completing successively we get ${\textrm{lk}}(22)$ = $C_{14}([{\textbf{23}}$, 9, 8, 13, 14, 20, ${\textbf{21}}$], 5, 4, 16, ${\textbf{15}}$, 17, 2, 3), ${\textrm{lk}}(16)$ = $C_{14}([{\textbf{15}}$, 10, 11, 12, 19, 18, ${\textbf{17}}$], 2, 3, 23, ${\textbf{22}}$, 21, 5, 4), ${\textrm{ lk}}(17)$ = $C_{14}([{\textbf{16}}$, 15, 10, 11, 12, 19, ${\textbf{18}}$], 20, 14, 1, ${\textbf 2}$, 3, 23, 22), ${\textrm{lk}}(18)$ = $C_{14}([{\textbf{17}}$, 16, 15, 10, 11, 12, ${\textbf{19}}$], 6, 5, 21, ${\textbf{20}}$, 14, 1, 2), ${\textrm{lk}}(1)$ = $C_{14}([{\textbf 2}$, 3, 4, 5, 6, 7, ${\textbf 0}$], 11, 12, 13, ${\textbf{14}}$, 20, 18, 17), ${\textrm{lk}}(2)$ = $C_{14}([{\textbf 3}$, 4, 5, 6, 7, 0, ${\textbf 1}$], 14, 20, 18, ${\textbf{17}}$, 16, 22, 23), ${\textrm{lk}}(3)$ = $C_{14}([{\textbf 4}$, 5, 6, 7, 0, 1, ${\textbf 2}$], 17, 16, 22, ${\textbf{23}}$, 9, 10, 15), ${\textrm{lk}}(4)$ = $C_{14}([{\textbf 5}$, 6, 7, 0, 1, 2, ${\textbf 3}$], 23, 9, 10, ${\textbf{15}}$, 16, 22, 21), ${\textrm{lk}}(14)$ = $C_{14}([{\textbf{20}}$, 21, 22, 23, 9, 8, ${\textbf{13}}$], 12, 11, 0, ${\textbf 1}$, 2, 17, 18) and ${\textrm{lk}}(12)$ = $C_{14}([{\textbf{11}}$, 10, 15, 16, 17, 18, ${\textbf{19}}$],6, 7, 8, ${\textbf{13}}$, 14, 1, 0). This is $N_1(6^2, 8)$ as given in Section \ref{example 4.1}.

\noindent {\bf Subcase 2.2\,:} If $(j, i) = (21, 22)$ then ${\textrm{lk}}(6)$ = $C_{14}([{\textbf 5}$, 4, 3, 2, 1, 0, ${\textbf 7}$], 8, 13, 12, ${\textbf{19}}$, 18, 22, 21). This implies ${\textrm{lk}}(5)$ = $C_{14}([{\textbf 4}$, 3, 2, 1, 0, 7, ${\textbf 6}$], 19, 18, 22, ${\textbf{21}}$, 20, $l, k$) for some $k, l \in V$. Then we see that $(k, l) \in \{$(15, 16), (16, 17), (17, 16)$\}$. In case $(k, l) = (16, 17)$, completing ${\textrm{lk}}(5)$ and ${\textrm{lk}}(20)$ we see that ${\textrm{lk}}(15)$ can not be completed. So we search for $(k, l) \in \{$(15, 16), (17, 16)$\}$.

\noindent {\bf Subcase 2.2.1\,:} If $(k, l) = (15, 16)$ then ${\textrm{lk}}(5)$ = $C_{14}([{\textbf{4}}$, 3, 2, 1, 0, 7, ${\textbf{6}}$], 19, 18, 22, ${\textbf{21}}$, 20, 16, 15), completing successively we get ${\textrm{lk}}(21)$ = $C_{14}([{\textbf{20}}$, 14, 13, 8, 9, 23, ${\textbf{22}}$], 18, 19, 6, ${\textbf{5}}$, 4, 15, 16), ${\textrm{lk}}(2)$ = $C_{14}([{\textbf{3}}$, 4, 5, 6, 7, 0, ${\textbf{1}}$], 14, 20, 16, ${\textbf{17}}$, 18, 22, 23), ${\textrm{lk}}(17)$ = $C_{14}([{\textbf{16}}$, 15, 10, 11, 12, 19, ${\textbf{18}}$],22, 23, 3, ${\textbf{2}}$, 1, 14, 20), ${\textrm{lk}}(3)$ = $C_{14}([{\textbf{4}}$, 5, 6, 7, 0, 1, ${\textbf{2}}$], 17, 18, 22, ${\textbf{23}}$, 9, 10, 15), ${\textrm{lk}}(23)$ = $C_{14}([{\textbf{9}}$, 8, 13, 14, 20, 21, ${\textbf{22}}$], 18, 17, 2, ${\textbf{3}}$, 4, 15, 10), ${\textrm{lk}}(4)$ = $C_{14}([{\textbf{5}}$, 6, 7, 0, 1, 2, ${\textbf{3}}$], 23, 9, 10, ${\textbf{15}}$, 16, 20, 21), ${\textrm{lk}}(15)$ = $C_{14}([{\textbf{16}}$, 17, 18, 19, 12, 11, ${\textbf{10}}$], 9, 23, 3, ${\textbf{4}}$, 5, 21, 20), ${\textrm{lk}}(5)$ = $C_{14}([{\textbf{4}}$, 3, 2, 1, 0, 7, ${\textbf{6}}$], 19, 18, 22, ${\textbf{21}}$, 20, 16, 15), ${\textrm{lk}}(9)$ = $C_{14}([{\textbf{8}}$, 13, 14, 20, 21, 22, ${\textbf{23}}$], 3, 4, 15, ${\textbf{10}}$, 11, 0, 7), ${\textrm{lk}}(10)$ = $C_{14}([{\textbf{15}}$, 16, 17, 18, 19, 12, ${\textbf{11}}$], 0, 7, 8, ${\textbf{9}}$, 23, 3, 4), ${\textrm{lk}}(14)$ = $C_{14}([{\textbf{20}}$, 21, 22, 23, 9, 8, ${\textbf{13}}$], 12, 11, 0, ${\textbf{1}}$, 2, 17, 16), ${\textrm{lk}}(16)$ = $C_{14}([{\textbf{15}}$, 10, 11, 12, 19, 18, ${\textbf{17}}$], 2, 1, 14, ${\textbf{20}}$, 21, 5, 4), ${\textrm{lk}}(20)$ = $C_{14}([{\textbf{21}}$, 22, 23, 9, 8, 13, ${\textbf{14}}$], 1, 2, 17, ${\textbf{16}}$, 15, 4, 5), ${\textrm{lk}}(18)$ = $C_{14}([{\textbf{17}}$, 16, 15, 10, 11, 12, ${\textbf{19}}$], 6, 5, 21, ${\textbf{22}}$, 23, 3, 2), ${\textrm{lk}}(22)$ = $C_{14}([{\textbf{23}}$, 9, 8, 13, 14, 20, ${\textbf{21}}$], 5, 6, 19, ${\textbf{18}}$, 17, 2, 3). This is isomorphic to $N_2(6^2, 8)$ as given in Section \ref{example 4.1}, by the map (0, 5, 2, 7, 4, 1, 6, 3)(8, 19, 23, 11, 21, 15, 14, 17)(9, 12, 22, 10, 13, 18)(16, 20).

\noindent {\bf Subcase 2.2.2\,:} If $(k, l) = (17, 16)$ then ${\textrm{lk}}(5)$ = $C_{14}([{\textbf{4}}$, 3, 2, 1, 0, 7, ${\textbf{6}}$], 19, 18, 22, ${\textbf{21}}$, 20, 16, 17), ${\textrm{lk}}(21)$ = $C_{14}([{\textbf{20}}$, 14, 13, 8, 9, 23, ${\textbf{22}}$], 18, 19, 6, ${\textbf{5}}$, 4, 17, 16), completing successively we get
${\textrm{lk}}(1)$ = $C_{14}([{\textbf{0}}$, 7, 6, 5, 4, 3, ${\textbf{2}}$], 17, 16, 20, ${\textbf{14}}$, 13, 12, 11), ${\textrm{lk}}(14)$ = $C_{14}([{\textbf{20}}$, 21, 22, 23, 9, 8, ${\textbf{13}}$], 12, 11, 0, ${\textbf{1}}$, 2, 15, 16), ${\textrm{lk}}(2)$ = $C_{14}([{\textbf{3}}$, 4, 5, 6, 7, 0, ${\textbf{1}}$], 14, 20, 16, ${\textbf{15}}$, 10, 9, 23), ${\textrm{lk}}(15)$ = $C_{14}$ $([{\textbf{16}}$, 17, 18, 19, 12, 11, ${\textbf{10}}$], 9, 23, 3, ${\textbf{2}}$, 1, 14, 20), ${\textrm{lk}}(3)$ = $C_{14}([{\textbf{4}}$, 5, 6, 7, 0, 1, ${\textbf{2}}$], 15, 10, 9, ${\textbf{23}}$, 22, 18, 17), ${\textrm{lk}}(23)$ = $C_{14}([{\textbf{9}}$, 8, 13, 14, 20, 21, ${\textbf{22}}$], 18, 17, 4, ${\textbf{3}}$, 2, 15, 10), ${\textrm{lk}}(4)$ = $C_{14}([{\textbf{5}}$, 6, 7, 0, 1, 2, ${\textbf{3}}$], 23, 22, 18, ${\textbf{17}}$, 16, 20, 21), ${\textrm{lk}}(17)$ = $C_{14}($[{\textbf{16}}$, 15, 10, 11, 12, 19, ${\textbf{18}}$], 22, 23, 3, ${\textbf{4}}$, 5, 21, 20$), ${\textrm{lk}}(9)= C_{14}([{\textbf{8}}$, 13, 14, 20, 21, 22, ${\textbf{23}}$], 3, 2, 15, ${\textbf{10}}$, 11, 0, 7), ${\textrm{lk}}(10)$ = $C_{14}([{\textbf{15}}$, 16, 17, 18, 19, 12, ${\textbf{11}}$], 0, 7, 8, ${\textbf{9}}$, 23, 3, 2), ${\textrm{lk}}(16)$ = $C_{14}([{\textbf{15}}$, 10, 11, 12, 19, 18, ${\textbf{17}}$], 4, 5, 21, ${\textbf{20}}$, 14, 1, 2), ${\textrm{lk}}(20)$ = $C_{14}([{\textbf{21}}$, 22, 23, 9, 8, 13, ${\textbf{14}}$], 1, 2, 15, ${\textbf{16}}$, 17, 4, 5), ${\textrm{lk}}(18)$ = $C_{14}([{\textbf{17}}$, 16, 15, 10, 11, 12, ${\textbf{19}}$], 6, 5, 21, ${\textbf{22}}$, 23, 3, 4), ${\textrm{lk}}(22)$ = $C_{14}([{\textbf{23}}$, 9, 8, 13, 14, 20, ${\textbf{21}}$], 5, 6, 19, ${\textbf{18}}$, 17, 4, 3). This map is isomorphic to $N_1(6^2, 8)$ by the map (0, 13)(1, 14)(2, 20)(3, 21)(4, 22)(5, 23)(6, 9)(7, 8)(10, 19)(11, 12)(15, 18)(16, 17). 

\smallskip

\noindent{\bf Case 3\,:} If $a = 20$ then we see that $b \in \{$13, 14, 21$\}$, $i.e.$, $(a, b) \in \{$(20, 13), (20, 14), (20, 21)$\}$. If $(a, b) = (20, 13)$ then successively we get $c = 14$, $d = 21$, $e = 22$ and $f = 23$. This implies ${\textrm{lk}}(8)$ = $C_{14}([{\textbf{20}}$, 13, 14, 21, 22, 23, ${\textbf{9}}$], 10, 11, 0, ${\textbf{7}}$, 6, h, g), where $(g, h) \in \{$(16, 15), (16, 17), (17, 16), (17, 18), (18, 17)$\}$. If $(g, h) \in \{$(17, 16), (17, 18)$\}$ then considering ${\textrm{lk}}(8)$ and ${\textrm{lk}}(6)$ successively we see that ${\textrm{lk}}(15)$ or ${\textrm{lk}}(19)$ can not be completed. For the remaining values of $(g, h)$ we have following\,:

If $(g, h) = (16, 15)$ then ${\textrm{lk}}(8)$ = $C_{14}([{\textbf{9}}$, 23, 22, 21, 14, 13, ${\textbf{20}}$], 16, 15, 6, ${\textbf{7}}$, 0, 11, 10). This implies ${\textrm{lk}}(6)$ = $C_{14}([{\textbf{5}}$, 4, 3, 2, 1, 0, ${\textbf{7}}$], 8, 20, 16, ${\textbf{15}}$, 10, 9, 23), ${\textrm{lk}}(9)$ = $C_{14}([{\textbf{23}}$, 22, 21, 14, 13, 20, ${\textbf{8}}$], 7, 0, 11, ${\textbf{10}}$, 15, 6, 5) and ${\textrm{lk}}(5)$ = $C_{14}([{\textbf{4}}$, 3, 2, 1, 0, 7, ${\textbf{6}}$], 15, 10, 9, ${\textbf{23}}$, 22, i, j) for some $i, j \in V(O_2)$. Observe that $(i, j)\in\{$(17, 18), (18, 17), (18, 19)$\}$. If $(i, j) = (17, 18)$ then considering ${\textrm{lk}}(22)$ we see 13\,21 as an edge and a non-edge both and if $(i, j) = (18, 17)$ or (18, 19) then successively considering ${\textrm{lk}}(5)$ and ${\textrm{lk}}(4)$ we see $\deg(13) > 3$ or $\deg(20) > 3$. So $(g, h) \neq (16, 15)$. If $(g, h) = (16, 17)$ then ${\textrm{lk}}(8)$ = $C_{14}([{\textbf{9}}$, 23, 22, 21, 14, 13, ${\textbf{20}}$], 16, 17, 6, ${\textbf{7}}$, 0, 11, 10) and ${\textrm{lk}}(6)$ = $C_{14}([{\textbf{5}}$, 4, 3, 2, 1, 0, ${\textbf{7}}$], 8, 20, 16, ${\textbf{17}}$, 18, $j, i)$ for some $i, j \in V(O_3)$. In this case $(j, i)\in\{$(21, 22), (22, 21), (22, 23), (23, 22)$\}$. If $(j, i) = (21, 22)$ then successively considering ${\textrm{lk}}(6)$ and ${\textrm{lk}}(5)$ we see that ${\textrm{lk}}(23)$ can not be completed. If $(j, i) = (22, 21)$ then we see that $0, 1 \in V(O_2)$. A contradiction, as $0, 1 \in V(O_1)$. If $(j, i) = (22, 23)$ then considering ${\textrm{lk}}(6)$, ${\textrm{lk}}(5)$ and ${\textrm{lk}}(4)$ successively we get $\deg(13) > 3$. If $(j, i) = (23, 22)$ then considering ${\textrm{lk}}(18)$ we see that $15\,19$ is simultaneously an edge and a non-edge of $N$. So $(g, h) \neq (16, 17)$. If $(g, h) = (18, 17)$ then ${\textrm{lk}}(8)$ = $C_{14}([{\textbf{9}}$, 23, 22, 21, 14, 13, ${\textbf{20}}$], 18, 17, 6, ${\textbf{7}}$, 0, 11, 10), ${\textrm{lk}}(6)$ = $C_{14}([{\textbf{5}}$, 4, 3, 2, 1, 0, ${\textbf{7}}$], 8, 20, 18, ${\textbf{17}}$, 16, i, j) for some $i, j \in V(O_3)$. Now proceeding as in previous case, we see that the map does not exist.

\noindent{\bf Subcase 3.1\,:} If $(a, b) = (20, 14)$ then successively we get $c = 13$, $d = 21$, $e = 22$ and $f = 23$. This implies ${\textrm{lk}}(8)$ = $C_{14}([{\textbf{20}}$, 14, 13, 21, 22, 23, ${\textbf{9}}$], 10, 11, 0, ${\textbf{7}}$, 6, $h, g)$ for some $g, h \in V(O_2)$. In this case $(g, h) \in \{$(16, 15), (16, 17), (17, 16), (17, 18), (18, 17), (18, 19)$\}$. If $(g, h) \in \{$(17, 16), (17, 18)$\}$ then considering ${\textrm{lk}}(8)$ and ${\textrm{lk}}(6)$ successively we see ${\textrm{lk}}(15)$ or ${\textrm{lk}}(19)$ can not be completed. For the remaining values of $(g, h)$ we have following\,:

\noindent{\bf Subcase 3.1.1\,:} If $(g, h) = (16, 15)$ then ${\textrm{lk}}(8)$ = $C_{14}([{\textbf{20}}$, 14, 13, 21, 22, 23, ${\textbf{9}}$], 10, 11, 0, ${\textbf{7}}$, 6, 15, 16), ${\textrm{lk}}(7)$ = $C_{14}([{\textbf{0}}$, 1, 2, 3, 4, 5, ${\textbf{6}}$], 15, 16, 20, ${\textbf{8}}$, 9, 10, 11), ${\textrm{lk}}(6)$ = $C_{14}([{\textbf{5}}$, 4, 3, 2, 1, 0, ${\textbf{7}}$], 8, 20, 16, ${\textbf{15}}$, 10, 9, 23), ${\textrm{lk}}(15)$ = $C_{14}([{\textbf{16}}$, 17, 18, 19, 12, 11, ${\textbf{10}}$], 9, 23, 5, ${\textbf{6}}$, 7, 8, 20), and ${\textrm{lk}}(5)$ = $C_{14}([{\textbf{4}}$, 3, 2, 1, 0, 7, ${\textbf{6}}$],15, 10, 9, ${\textbf{23}}$, 22, $i, j)$ for some $i, j \in V(O_2)$. Observe that $(i, j)\in\{$(18, 19), (19, 18)$\}$. If $(i, j) = (19, 18)$ then considering ${\textrm{lk}}(5)$ we see that ${\textrm{lk}}(22)$ can not be completed. On the other hand when $(i, j) = (18, 19)$ then ${\textrm{lk}}(23)$ = $C_{14}([{\textbf{9}}$, 8, 20, 14, 13, 21, ${\textbf{22}}$], 18, 19, 4, ${\textbf{5}}$, 6, 15, 10), completing successively we get ${\textrm{lk}}(1)$ = $C_{14}([{\textbf{0}}$, 7, 6, 5, 4, 3, ${\textbf{2}}$], 17, 16, 20, ${\textbf{14}}$, 13, 12, 11), ${\textrm{lk}}(14)$ = $C_{14}([{\textbf{20}}$, 8, 9, 23, 22, 21, ${\textbf{13}}$], 12, 11, 0, ${\textbf{1}}$, 2, 17, 16), ${\textrm{lk}}(2)$ = $C_{14}([{\textbf{3}}$, 4, 5, 6, 7, 0, ${\textbf{1}}$], 14, 20, 16, ${\textbf{17}}$, 18, 22, 21), ${\textrm{lk}}(17)$ = $C_{14}([{\textbf{16}}$, 15, 10, 11, 12, 19, ${\textbf{18}}$], 22, 21, 3, ${\textbf{2}}$, 1, 14, 20), ${\textrm{lk}}$ $(3)$ = $C_{14}([{\textbf{4}}$, 5, 6, 7, 0, 1, ${\textbf{2}}$], 17, 18, 22, ${\textbf{21}}$, 13, 12, 19), ${\textrm{lk}}(21)$ = $C_{14}([{\textbf{13}}$, 14, 20, 8, 9, 23, ${\textbf{22}}$], 18, 17, 2, ${\textbf{3}}$, 4, 19, 12), ${\textrm{lk}}(4)$ = $C_{14}([{\textbf{5}}$, 6, 7, 0, 1, 2, ${\textbf{3}}$], 21, 13, 12, ${\textbf{19}}$, 18, 22, 23), ${\textrm{lk}}(19)$ = $C_{14}([{\textbf{12}}$, 11, 10, 15, 16, 17, ${\textbf{18}}$], 22, 23, 5, ${\textbf{4}}$, 3, 21, 13), ${\textrm{lk}}(9)$ = $C_{14}([{\textbf{8}}$, 20, 14, 13, 21, 22, ${\textbf{23}}$], 5, 6, 15, ${\textbf{10}}$, 11, 0, 7), ${\textrm{lk}}(10)$ = $C_{14}([{\textbf{15}}$, 16, 17, 18, 19, 12, ${\textbf{11}}$], 0, 7, 8, ${\textbf{9}}$, 23, 5, 6), ${\textrm{lk}}(12)$ = $C_{14}([{\textbf{11}}$, 10, 15, 16, 17, 18, ${\textbf{19}}$], 4, 3, 21, ${\textbf{13}}$, 14, 1, 0), ${\textrm{lk}}(13)$ = $C_{14}([{\textbf{14}}$, 20, 8, 9, 23, 22, ${\textbf{21}}$], 3, 4, 19, ${\textbf{12}}$, 11, 0, 1), ${\textrm{lk}}(16)$ $= C_{14}([{\textbf{15}}$, 10, 11, 12, 19, 18, ${\textbf{17}}$], 2, 1, 14, ${\textbf{20}}$, 8, 7, 6), ${\textrm{lk}}(20)$ = $C_{14}([{\textbf{8}}$, 9, 23, 22, 21, 13, ${\textbf{14}}$], 1, 2, 17, ${\textbf{16}}$, 15, 6, 7), ${\textrm{lk}}(18)$ = $C_{14}([{\textbf{17}}$, 16, 15, 10, 11, 12, ${\textbf{19}}$], 4, 5, 23, ${\textbf{22}}$, 21, 3, 2), ${\textrm{lk}}(22)$ = $C_{14}([{\textbf{23}}$, 9, 8, 20, 14, 13, ${\textbf{21}}$], 3, 2, 17, ${\textbf{18}}$, 19, 4, 5). This is isomorphic to $N_1(6^2, 8)$ by the map (0, 14, 18, 4, 23, 10)(1, 20, 19, 3, 22, 15, 7, 13, 17, 5, 9, 11)(2, 21, 16, 6, 8, 12).

\noindent{\bf Subcase 3.1.2\,:} If $(g, h) = (16, 17)$ then ${\textrm{lk}}(8)$ = $C_{14}([{\textbf{20}}$, 14, 13, 21, 22, 23, ${\textbf{9}}$], 10, 11, 0, ${\textbf{7}}$, 6, 17, 16), ${\textrm{lk}}(7)$ = $C_{14}([{\textbf{0}}$, 1, 2, 3, 4, 5, ${\textbf{6}}$], 17, 16, 20, ${\textbf{8}}$, 9, 10, 11), ${\textrm{lk}}(6)$ = $C_{14}([{\textbf{5}}$, 4, 3, 2, 1, 0, ${\textbf{7}}$], 8, 20, 16, ${\textbf{17}}$, 18, $i, j)$ for some $i, j \in V(O_3)$. We see easily that $(i, j)\in\{$(21, 22), (22, 21), (22, 23), (23, 22)$\}$. If $(i, j) = (21, 22)$ then successively considering ${\textrm{lk}}(6)$ and ${\textrm{lk}}(5)$ we see that ${\textrm{lk}}(23)$ can not be completed. If $(i, j) = (23, 22)$ then considering ${\textrm{lk}}(6)$ we see that ${\textrm{lk}}(18)$ can not be completed. If $(i, j) = (22, 23)$ then successively considering ${\textrm{lk}}(6), ${\textrm{lk}}(5) and ${\textrm{lk}}(4)$ we get $\deg(14) > 3$, a contradiction. If $(i, j) = (22, 21)$ then ${\textrm{lk}}(17)$ = $C_{14}([{\textbf{16}}$, 15, 10, 11, 12, 19, ${\textbf{18}}$], 22, 21, 5, ${\textbf{6}}$, 7, 8, 20), completing successively we get ${\textrm{lk}}(5)$ = $C_{14}([{\textbf{4}}$, 3, 2, 1, 0, 7, ${\textbf{6}}$], 17, 18, 22, ${\textbf{21}}$, 13, 12, 19), ${\textrm{lk}}(21)$ = $C_{14}([{\textbf{13}}$, 14, 20, 8, 9, 23, ${\textbf{22}}$], 18, 17, 6, ${\textbf{5}}$, 4, 19, 12), ${\textrm{lk}}(4)$ = $C_{14}([{\textbf{5}}$, 6, 7, 0, 1, 2, ${\textbf{3}}$], 23, 22, 18, ${\textbf{19}}$, 12, 13, 21), ${\textrm{lk}}(19)$ = $C_{14}([{\textbf{12}}$, 11, 10, 15, 16, 17, ${\textbf{18}}$], 22, 23, 3, ${\textbf{4}}$, 5, 21, 13), ${\textrm{lk}}(1)$ = $C_{14}([{\textbf{0}}$, 7, 6, 5, 4, 3, ${\textbf{2}}$], 15, 16, 20, ${\textbf{14}}$, 13, 12, 11), ${\textrm{lk}}$ $(14)$ = $C_{14}([{\textbf{20}}$, 8, 9, 23, 22, 21, ${\textbf{13}}$], 12, 11, 0, ${\textbf{1}}$, 2, 15, 16), ${\textrm{lk}}(2)$ = $C_{14}([{\textbf{3}}$, 4, 5, 6, 7, 0, ${\textbf{1}}$], 14, 20, 16, ${\textbf{15}}$, 10, 9, 23), ${\textrm{lk}}(15)$ = $C_{14}([{\textbf{16}}$, 17, 18, 19, 12, 11, ${\textbf{10}}$], 9, 23, 3, ${\textbf{2}}$, 1, 14, 20), ${\textrm{lk}}(3)$ = $C_{14}([{\textbf{4}}$, 5, 6, 7, 0, 1, ${\textbf{2}}$], 15, 10, 9, ${\textbf{23}}$, 22, 18, 19), ${\textrm{lk}}(23)$ = $C_{14}([{\textbf{9}}$, 8, 20, 14, 13, 21, ${\textbf{22}}$], 18, 19, 4, ${\textbf{3}}$, 2, 15, 10), ${\textrm{lk}}(9)$ = $C_{14}([{\textbf{8}}$, 20, 14, 13, 21, 22, ${\textbf{23}}$], 3, 2, 15, ${\textbf{10}}$, 11, 0, 7), ${\textrm{lk}}(10)$ = $C_{14}([{\textbf{15}}$, 16, 17, 18, 19, 12, ${\textbf{11}}$], 0, 7, 8, ${\textbf{9}}$, 23, 3, 2), ${\textrm{lk}}(12)$ = $C_{14}([{\textbf{11}}$, 10, 15, 16, 17, 18, ${\textbf{19}}$], 4, 5, 21, ${\textbf{13}}$, 14, 1, 0), ${\textrm{lk}}(13)$ $= C_{14}([{\textbf{14}}$, 20, 8, 9, 23, 22, ${\textbf{21}}$], 5, 4, 19, ${\textbf{12}}$, 11, 0, 1), ${\textrm{lk}}(16)$ = $C_{14}([{\textbf{15}}$, 10, 11, 12, 19, 18, ${\textbf{17}}$], 6, 7, 8, ${\textbf{20}}$, 14, 1, 2), ${\textrm{lk}}(20)$ = $C_{14}([{\textbf{8}}$, 9, 23, 22, 21, 13, ${\textbf{14}}$], 1, 2, 15, ${\textbf{16}}$, 17, 6, 7), ${\textrm{lk}}(18)$ = $C_{14}([{\textbf{17}}$, 16, 15, 10, 11, 12, ${\textbf{19}}$], 4, 3, 23, ${\textbf{22}}$, 21, 5, 6), ${\textrm{lk}}(22)$ = $C_{14}([{\textbf{23}}$, 9, 8, 20, 14, 13, ${\textbf{21}}$], 5, 6, 17, ${\textbf{18}}$, 19, 4, 3). This is $N_2(6^2, 8)$ as given in Section \ref{example 4.1}.

\noindent{\bf Subcase 3.1.3\,:} If $(g, h) = (18, 17)$ then ${\textrm{lk}}(8)$ = $C_{14}([{\textbf{20}}$, 14, 13, 21, 22, 23, ${\textbf{9}}$], 10, 11, 0, ${\textbf{7}}$, 6, 17, 18), ${\textrm{lk}}(7)$ = $C_{14}([{\textbf{0}}$, 1, 2, 3, 4, 5, ${\textbf{6}}$], 17, 18, 20, ${\textbf{8}}$, 9, 10, 11). This implies ${\textrm{lk}}(6)$ = $C_{14}([{\textbf{5}}$, 4, 3, 2, 1, 0, ${\textbf{7}}$], 8, 20, 18, ${\textbf{17}}$, 16, $i, j)$ for some $i, j \in V(O_3)$, observe that $(i, j)\in\{$(21, 22), (22, 21), (22, 23)$\}$. If $(i, j) = (21, 22)$ then considering ${\textrm{lk}}(6)$ we see that ${\textrm{lk}}(5)$ and ${\textrm{lk}}(23)$ can not be completed. If $(i, j) = (22, 21)$ then successively considering ${\textrm{lk}}(6)$, ${\textrm{lk}}(5)$, ${\textrm{lk}}(21)$, ${\textrm{lk}}(13)$, ${\textrm{lk}}(12)$, ${\textrm{lk}}(19)$ and ${\textrm{lk}}(4)$ we get $\deg(14) > 3$. A contradiction. If $(i, j) = (22, 23)$ then ${\textrm{lk}}(6)$ = $C_{14}([{\textbf{7}}$, 0, 1, 2, 3, 4, ${\textbf{5}}$], 23, 22, 16, ${\textbf{17}}$, 18, 20, 8), ${\textrm{lk}}(17)$ = $C_{14}([{\textbf{16}}$, 15, 10, 11, 12, 19, ${\textbf{18}}$], 20, 8, 7, ${\textbf{6}}$, 5, 23, 22), completing successively we get ${\textrm{lk}}(1)$ = $C_{14}([{\textbf{0}}$, 7, 6, 5, 4, 3, ${\textbf{2}}$], 19, 18, 20, ${\textbf{14}}$, 13, 12, 11), ${\textrm{lk}}(14)$ = $C_{14}([{\textbf{20}}$, 8, 9, 23, 22, 21, ${\textbf{13}}$], 12, 11, 0, ${\textbf{1}}$, 2, 19, 18), ${\textrm{lk}}(2)$ = $C_{14}([{\textbf{3}}$, 4, 5, 6, 7, 0, ${\textbf{1}}$], 14, 20, 18, ${\textbf{19}}$, 12, 13, 21), ${\textrm{lk}}(19)$ = $C_{14}([{\textbf{12}}$, 11, 10, 15, 16, 17, ${\textbf{18}}$], 20, 14, 1, ${\textbf{2}}$, 3, 21, 13), ${\textrm{lk}}(3)$ = $C_{14}([{\textbf{4}}$, 5, 6, 7, 0, 1, ${\textbf{2}}$], 19, 12, 13, ${\textbf{21}}$, 22, 16, 15), ${\textrm{lk}}(21)$ = $C_{14}([{\textbf{13}}$, 14, 20, 8, 9, 23, ${\textbf{22}}$], 16, 15, 4, ${\textbf{3}}$, 2, 19, 12), ${\textrm{lk}}(4)$ = $C_{14}([{\textbf{5}}$, 6, 7, 0, 1, 2, ${\textbf{3}}$], 21, 22, 16, ${\textbf{15}}$, 10, 9, 23), ${\textrm{lk}}(15)$ = $C_{14}([{\textbf{16}}$, 17, 18, 19, 12, 11, ${\textbf{10}}$], 9, 23, 5, ${\textbf{4}}$, 3, 21, 22), ${\textrm{lk}}(5)$ = $C_{14}([{\textbf{4}}$, 3, 2, 1, 0, 7, ${\textbf{6}}$], 17, 16, 22, ${\textbf{23}}$, 9, 10, 15), ${\textrm{lk}}(23)$ = $C_{14}([{\textbf{9}}$, 8, 20, 14, 13, 21, ${\textbf{22}}$], 16, 17, 6, ${\textbf{5}}$, 4, 15, 10), ${\textrm{lk}}(9)$ = $C_{14}([{\textbf{8}}$, 20, 14, 13, 21, 22, ${\textbf{23}}$], 5, 4, 15, ${\textbf{10}}$, 11, 0, 7), ${\textrm{lk}}(10)$ = $C_{14}([{\textbf{15}}$, 16, 17, 18, 19, 12, ${\textbf{11}}$], 0, 7, 8, ${\textbf{9}}$, 23, 5, 4), ${\textrm{lk}}(12)$ = $C_{14}([{\textbf{11}}$, 10, 15, 16, 17, 18, ${\textbf{19}}$], 2, 3, 21, ${\textbf{13}}$, 14, 1, 0), ${\textrm{lk}} (13)$ = $C_{14}([{\textbf{14}}$, 20, 8, 9, 23, 22, ${\textbf{21}}$], 3, 2, 19, ${\textbf{12}}$, 11, 0, 1), ${\textrm{lk}}(16)$ = $C_{14}([{\textbf{15}}$, 10, 11, 12, 19, 18, ${\textbf{17}}$], 6, 5, 23, ${\textbf{22}}$, 21, 3, 4), ${\textrm{lk}}(18)$ = $C_{14}([{\textbf{17}}$, 16, 15, 10, 11, 12, ${\textbf{19}}$], 2, 1, 14, ${\textbf{20}}$, 8, 7, 6), ${\textrm{lk}}(20)$ = $C_{14}([{\textbf{8}}$, 9, 23, 22, 21, 13, ${\textbf{14}}$], 1, 2, 19, ${\textbf{18}}$, 17, 6, 7), ${\textrm{lk}}(22)$ $= C_{14}([{\textbf{23}}$, 9, 8, 20, 14, 13, ${\textbf{21}}$], 3, 4, 15, ${\textbf{16}}$, 17, 6, 5). This is isomorphic to $N_1(6^2, 8)$ by the map (0, 21, 10, 6, 14, 16)(1, 22, 11, 5, 13, 15, 7, 20, 17)(2, 23, 12, 4, 8, 18)(3, 9, 19).

\noindent{\bf Subcase 3.1.4\,:} If $(g, h) = (18, 19)$ then ${\textrm{lk}}(8)$ = $C_{14}([{\textbf{20}}$, 14, 13, 21, 22, 23, ${\textbf{9}}$], 10, 11, 0, ${\textbf{7}}$, 6, 19, 18), ${\textrm{lk}}(7)$ = $C_{14}([{\textbf{0}}$, 1, 2, 3, 4, 5, ${\textbf{6}}$],19, 18, 20, ${\textbf{8}}$, 9, 10, 11), ${\textrm{lk}}(6)$ = $C_{14}([{\textbf{5}}$, 4, 3, 2, 1, 0, ${\textbf{7}}$], 8, 20, 18, ${\textbf{19}}$, 12, 13, 21), ${\textrm{lk}}(19)$ = $C_{14}([{\textbf{12}}$, 11, 10, 15, 16, 17, ${\textbf{18}}$], 20, 8, 7, ${\textbf{6}}$, 5, 21, 13) and ${\textrm{lk}}(5)$ = $C_{14}([{\textbf{4}}$, 3, 2, 1, 0, 7, ${\textbf{6}}$], 19, 12, 13, ${\textbf{21}}$, 22, $i, j)$ for some $i, j \in V(O_2)$. It is easy to see that $(i, j) \in\{$(15, 16), (16, 15), (16, 17), (17, 16)$\}$. If $(i, j) = (17, 16)$ then considering ${\textrm{lk}}(17)$ we see that 14\,23 is simultaneously an edge and a non-edge of $N$. If $(i, j) = (15, 16)$ then considering ${\textrm{lk}}(5)$ we see that ${\textrm{lk}}(4)$ and ${\textrm{lk}}(17)$ can not be completed. If $(i, j) = (16, 17)$ then successively considering ${\textrm{lk}}(5)$ and ${\textrm{lk}}(4)$ we get $\deg(14) > 3$. A contradiction. If $(i, j) = (16, 15)$ then ${\textrm{lk}}(21)$ = $C_{14}([{\textbf{13}}$, 14, 20, 8, 9, 23, ${\textbf{22}}$], 16, 15, 4, ${\textbf{5}}$, 6, 19, 12), completing successively we get ${\textrm{lk}}(1)$ = $C_{14}([{\textbf{0}}$, 7, 6, 5, 4, 3, ${\textbf{2}}$], 17, 18, 20, ${\textbf{14}}$, 13, 12, 11), ${\textrm{lk}}(14)$ = $C_{14}([{\textbf{20}}$, 8, 9, 23, 22, 21, ${\textbf{13}}$], 12, 11, 0, ${\textbf{1}}$, 2, 17, 18), ${\textrm{lk}}(2)$ = $C_{14}([{\textbf{3}}$, 4, 5, 6, 7, 0, ${\textbf{1}}$], 14, 20, 18, ${\textbf{17}}$, 16, 22, 23), ${\textrm{lk}}(17)$ = $C_{14}([{\textbf{16}}$, 15, 10, 11, 12, 19, ${\textbf{18}}$], 20, 14, 1, ${\textbf{2}}$, 3, 23, 22), ${\textrm{lk}}(3)$ = $C_{14}([{\textbf{4}}$, 5, 6, 7, 0, 1, ${\textbf{2}}$], 17, 16, 22, ${\textbf{23}}$, 9, 10, 15), ${\textrm{lk}}(23)$ = $C_{14}([{\textbf{9}}$, 8, 20, 14, 13, 21, ${\textbf{22}}$], 16, 17, 2, ${\textbf{3}}$, 4, 15, 10), ${\textrm{lk}}(4)$ = $C_{14}([{\textbf{5}}$, 6, 7, 0, 1, 2, ${\textbf{3}}$], 23, 9, 10, ${\textbf{15}}$, 16, 22, 21), ${\textrm{lk}}(15)$ = $C_{14}([{\textbf{16}}$, 17, 18, 19, 12, 11, ${\textbf{10}}$],9, 23, 3, ${\textbf{4}}$, 5, 21, 22), ${\textrm{lk}}(9)$ = $C_{14}([{\textbf{8}}$, 20, 14, 13, 21, 22, ${\textbf{23}}$], 3, 4, 15, ${\textbf{10}}$, 11, 0, 7), ${\textrm{lk}}(10)$ = $C_{14}$ $([{\textbf{15}}$, 16, 17, 18, 19, 12, ${\textbf{11}}$], 0, 7, 8, ${\textbf{9}}$, 23, 3, 4), ${\textrm{lk}}(12)$ = $C_{14}([{\textbf{11}}$, 10, 15, 16, 17, 18, ${\textbf{19}}$], 6, 5, 21, ${\textbf{13}}$, 14, 1, 0), ${\textrm{lk}}(13)$ = $C_{14}([{\textbf{14}}$, 20, 8, 9, 23, 22, ${\textbf{21}}$], 5, 6, 19, ${\textbf{12}}$, 11, 0, 1), ${\textrm{lk}}(16)$ = $C_{14}([{\textbf{15}}$, 10, 11, 12, 19, 18, ${\textbf{17}}$], 2, 3, 23, ${\textbf{22}}$, 21, 5, 4), ${\textrm{lk}}(18)$ = $C_{14}([{\textbf{17}}$, 16, 15, 10, 11, 12, ${\textbf{19}}$], 6, 7, 8, ${\textbf{20}}$, 14, 1, 2), ${\textrm{lk}}(20)$ = $C_{14}([{\textbf{8}}$, 9, 23, 22, 21, 13, ${\textbf{14}}$], 1, 2, 17, ${\textbf{18}}$, 19, 6, 7), ${\textrm{lk}}(22)$ = $C_{14}([{\textbf{23}}$, 9, 8, 20, 14, 13, ${\textbf{21}}$], 5, 4, 15, ${\textbf{16}}$, 17, 2, 3). This is isomorphic to $N_2(6^2, 8)$ by the map (0, 1, 2, 3, 4, 5, 6, 7)(8, 11, 14, 15, 21, 17, 23, 19)(9, 12, 20, 10, 13, 16, 22, 18).

\noindent{\bf Subcase 3.2\,:} If $(a, b) = (20, 21)$ then we see that $c \in \{$13, 14, 22$\}$.

If $c = 14$ then completing successively we get $d = 13$, $e = 22$, $f = 23$ and $(g, h) \in \{$(16, 15), (16, 17), (17, 16), (17, 18), (18, 17), (18, 19)$\}$. If $(g, h) = (16, 15)$ then considering ${\textrm{lk}}(8)$, ${\textrm{lk}}(9)$ we see that ${\textrm{lk}}(5)$ and ${\textrm{lk}}(22)$ can not be completed. If $(g, h) = (16, 17)$ then ${\textrm{lk}}(8)$ = $C_{14}([{\textbf{9}}$, 23, 22, 13, 14, 21, ${\textbf{20}}$], 16, 17, 6, ${\textbf{7}}$, 0, 11, 10) and ${\textrm{lk}}(6)$ = $C_{14}([{\textbf{5}}$, 4, 3, 2, 1, 0, ${\textbf{7}}$], 8, 20, 16, ${\textbf{17}}$, 18, $i, j)$ for some $i, j \in V(O_3)$. Observe that $(i, j)\in\{$(22, 23), (23, 22)$\}$. If $(i, j) = (22, 23)$ then successively considering ${\textrm{lk}}(5)$, ${\textrm{lk}}(22)$, ${\textrm{lk}}(13)$, ${\textrm{lk}}(4)$, ${\textrm{lk}}(12)$, ${\textrm{lk}}(19)$, ${\textrm{lk}}(18)$ and ${\textrm{lk}}(23)$ we see that ${\textrm{lk}}(2)$ and ${\textrm{lk}}(3)$ can not be completed. If $(i, j) = (23, 22)$ then successively considering ${\textrm{lk}}(5)$, ${\textrm{lk}}(22)$, ${\textrm{lk}}(13)$, ${\textrm{lk}}(4)$, ${\textrm{lk}}(12)$ and ${\textrm{lk}}(18)$ we see that ${\textrm{lk}}(2)$ and ${\textrm{lk}}(3)$ can not be completed. If $(g, h) = (17, 16)$ then ${\textrm{lk}}(8)$ = $C_{14}([{\textbf{9}}$, 23, 22, 13, 14, 21, ${\textbf{20}}$], 17, 16, 6, ${\textbf{7}}$, 0, 11, 10), ${\textrm{lk}}(6)$ = $C_{14}([{\textbf{5}}$, 4, 3, 2, 1, 0, ${\textbf{7}}$], 8, 20, 17, ${\textbf{16}}$, 15, $j, i)$ and ${\textrm{lk}}(15)$ = $C_{14}([{\textbf{10}}$, 11, 12, 19, 18, 17, ${\textbf{16}}$], 6, 5, $i$, ${\textit {\textbf j}}$, $k$, 23, 9) for some $i, j, k \in V(O_3)$. Now, it is easy to see that $i$, $j$, $k$ have no values in $V$ so that ${\textrm{lk}}(15)$ can be completed. In case $(g, h) = (17, 18)$ then considering ${\textrm{lk}}(8)$ we see ${\textrm{lk}}(19)$ can not be completed. If $(g, h) = (18, 19)$ then considering ${\textrm{lk}}(8)$ and ${\textrm{lk}}(6)$  we see that ${\textrm{lk}}(23)$ can not be completed. If $(g, h) = (18, 17)$ then ${\textrm{lk}}(8)$ = $C_{14}([{\textbf{9}}$, 23, 22, 13, 14, 21, ${\textbf{20}}$], 18, 17, 6, ${\textbf{7}}$, 0, 11, 10), ${\textrm{lk}}(6)$ = $C_{14}([{\textbf{5}}$, 4, 3, 2, 1, 0, ${\textbf{7}}$], 8, 20, 18, ${\textbf{17}}$, 16, $i, j)$ for some $i, j \in V(O_3)$. In this case $(i, j)\in\{$(22, 23), (23, 22)$\}$. If $(i, j) = (22, 23)$ then successively considering ${\textrm{lk}}(6)$, ${\textrm{lk}}(5)$, ${\textrm{lk}}(9)$, ${\textrm{lk}}(10)$, ${\textrm{lk}}(23)$ and ${\textrm{lk}}(4)$ we get $\deg(13) > 3$. If $(i, j) = (23, 22)$ then successively considering ${\textrm{lk}}(6)$, ${\textrm{lk}}(5)$, ${\textrm{lk}}(4)$, ${\textrm{lk}}(13)$, ${\textrm{lk}}(12)$,  ${\textrm{lk}}(19)$, ${\textrm{lk}}(18)$ and ${\textrm{lk}}(20)$ we see that length of ${\textrm{lk}}(2)$ is less than 14. So $c \neq 14$.

\smallskip

\noindent{\bf Subcase 3.2.1\,:} If $c = 13$ then $d = 14$, $e = 22$, $f = 23$ and $(g, h) \in \{$(16, 15), (16, 17), (17, 16), (17, 18), (18, 17), (18, 19)$\}$. If $(g, h) \in \{$(17, 16), (17, 18)$\}$ then considering ${\textrm{lk}}(8)$ and ${\textrm{lk}}(6)$ we see ${\textrm{lk}}(15)$ or ${\textrm{lk}}(19)$ can not be completed. For the remaining values of $(g, h)$ we have following subcases.

\noindent{\bf Subcase 3.2.1.1\,:} If $(g, h) = (18, 17)$ then ${\textrm{lk}}(8)$ = $C_{14}([{\textbf{9}}$, 23, 22, 14, 13, 21, ${\textbf{20}}$], 18, 17, 6, ${\textbf{7}}$, 0, 11, 10) and ${\textrm{lk}}(6)$ = $C_{14}([{\textbf{5}}$, 4, 3, 2, 1, 0, ${\textbf{7}}$], 8, 20, 18, ${\textbf{17}}$, 16, $i, j)$ for some $i, j \in V(O_3)$. Observe that $i \in \{$22, 23$\}$. If $i = 22$ then $j = 23$, now successively considering ${\textrm{lk}}(5), ${\textrm{lk}}(10), ${\textrm{lk}}(9)$ and ${\textrm{lk}}(4)$ we see $\deg(14) > 3$. A contradiction. If $i = 23$ then $j = 22$ this implies length of ${\textrm{lk}}(23)$ is less than 14. A contradiction. So $(g, h) \neq (18, 17)$.

\noindent{\bf Subcase 3.2.1.2\,:} If $(g, h) = (16, 15)$ then ${\textrm{lk}}(8)$ = $C_{14}([{\textbf{9}}$, 23, 22, 21, 13, 14, ${\textbf{20}}$], 16, 15, 6, ${\textbf{7}}$, 0, 11, 10), ${\textrm{lk}}(7)$ = $C_{14}([{\textbf{0}}$, 1, 2, 3, 4, 5, ${\textbf{6}}$], 15, 16, 20, ${\textbf{8}}$, 9, 10, 11), ${\textrm{lk}}(6)$ = $C_{14}([{\textbf{7}}$, 0, 1, 2, 3, 4, ${\textbf{5}}$],23, 9, 10, ${\textbf{15}}$, 16, 20, 8), ${\textrm{lk}}(15)$ = $C_{14}([{\textbf{16}}$, 17, 18, 19, 12, 11, ${\textbf{10}}$], 9, 23, 5, ${\textbf{6}}$, 7, 8, 20), ${\textrm{lk}}(5)$ = $C_{14}$ $([{\textbf{4}}$, 3, 2, 1, 0, 7, ${\textbf{6}}$],15, 10, 9, ${\textbf{23}}$, 22, $i, j)$, ${\textrm{lk}}(22)$ = $C_{14}([{\textbf{23}}$, 9, 8, 20, 21, 13, ${\textbf{14}}$], 21, 20, $k$, ${\textit {\textbf i}}$, $j$, 4, 5) for some $i, j, k \in V(O_2)$. In this case $i \in \{$17, 18, 19$\}$. If $i = 17$ then $j = 18$, now considering ${\textrm{lk}}(17)$ we see that
14\,21 is both an edge and a non-edge. If $i = 19$ then $j = 18$ and $k = 12$, now considering ${\textrm{lk}}(22)$ we see that 13\,14 is both an edge and a non-edge of $N$. If $i = 18$ then $j = 17$ and $k = 19$. This implies
${\textrm{lk}}(23)$ = $C_{14}([{\textbf{9}}$, 8, 20, 21, 13, 14, ${\textbf{22}}$], 18, 17, 4, ${\textbf{5}}$, 6, 15, 10), ${\textrm{lk}}(18)$ = $C_{14}([{\textbf{17}}$, 16, 15, 10, 11, 12, ${\textbf{19}}$], 2, 1, 14, ${\textbf{22}}$, 23, 5, 4). Now completing successively we get ${\textrm{lk}}(1)$ = $C_{14}([{\textbf{0}}$, 7, 6, 5, 4, 3, ${\textbf{2}}$], 19, 18, 22, ${\textbf{14}}$, 13, 12, 11), ${\textrm{lk}}(14)$ = $C_{14}([{\textbf{13}}$, 21, 20, 8, 9, 23, ${\textbf{22}}$],18, 19, 2, ${\textbf{1}}$, 0, 11, 12), ${\textrm{lk}}(2)$ = $C_{14}([{\textbf{3}}$, 4, 5, 6, 7, 0, ${\textbf{1}}$], 14, 22, 18, ${\textbf{19}}$, 12, 13, 21), ${\textrm{lk}}(19)$ = $C_{14}$ $([{\textbf{12}}$, 11, 10, 15, 16, 17, ${\textbf{18}}$], 22, 14, 1, ${\textbf{2}}$, 3, 21, 13), ${\textrm{lk}}(3)$ = $C_{14}([{\textbf{4}}$, 5, 6, 7, 0, 1, ${\textbf{2}}$], 19, 12, 13, ${\textbf{21}}$, 20, 16, 17), ${\textrm{lk}}(21)$ = $C_{14}([{\textbf{13}}$, 14, 22, 23, 9, 8, ${\textbf{20}}$], 16, 17, 4, ${\textbf{3}}$, 2, 19, 12),
${\textrm{lk}}(4)$ = $C_{14}([{\textbf{5}}$, 6, 7, 0, 1, 2, ${\textbf{3}}$], 21, 20, 16, ${\textbf{17}}$, 18, 22, 23), ${\textrm{lk}}(17)$ = $C_{14}([{\textbf{16}}$, 15, 10, 11, 12, 19, ${\textbf{18}}$], 22, 23, 5, ${\textbf{4}}$, 3, 21, 20), ${\textrm{lk}}(9)$ = $C_{14}([{\textbf{8}}$, 20, 21, 13, 14, 22, ${\textbf{23}}$], 5, 6, 15, ${\textbf{10}}$, 11, 0, 7), ${\textrm{lk}}(10)$ = $C_{14}([{\textbf{15}}$, 16, 17, 18, 19, 12, ${\textbf{11}}$], 0, 7, 8, ${\textbf{9}}$, 23, 5, 6), ${\textrm{lk}}(12)$ = $C_{14}([{\textbf{11}}$, 10, 15, 16, 17, 18, ${\textbf{19}}$], 2, 3, 21, ${\textbf{13}}$, 14, 1, 0), ${\textrm{lk}}(13)$ = $C_{14}([{\textbf{14}}$, 22, 23, 9, 8, 20, ${\textbf{21}}$], 3, 2, 19, ${\textbf{12}}$, 11, 0, 1), ${\textrm{lk}}(16)$ = $C_{14}([{\textbf{15}}$, 10, 11, 12, 19, 18, ${\textbf{17}}$], 4, 3, 21, ${\textbf{20}}$, 8, 7, 6), ${\textrm{lk}}$ $(20)$ = $C_{14}([{\textbf{8}}$, 9, 23, 22, 14, 13, ${\textbf{21}}$], 3, 4, 17, ${\textbf{16}}$, 15, 6, 7). This is isomorphic to $N_1(6^2, 8)$ by the map (0, 10, 23, 2, 12, 8, 4, 18, 14)(1, 11, 9, 3, 19, 13, 7, 15, 22)(5, 17, 20)(6, 16, 21).

\noindent{\bf Subcase 3.2.1.3\,:} If $(g, h) = (16, 17)$ then ${\textrm{lk}}(8)$ = $C_{14}([{\textbf{9}}$, 23, 22, 14, 13, 21, ${\textbf{20}}$], 16, 17, 6, ${\textbf{7}}$, 0, 11, 10), ${\textrm{lk}}(7)$ = $C_{14}([{\textbf{0}}$, 1, 2, 3, 4, 5, ${\textbf{6}}$], 17, 16, 20, ${\textbf{8}}$, 9, 10, 11) and ${\textrm{lk}}(6)$ = $C_{14}([{\textbf{5}}$, 4, 3, 2, 1, 0, ${\textbf{7}}$], 8, 20, 16, ${\textbf{17}}$, 18, $i, j)$ for some $i, j \in V(O_3)$. Observe that $i = 22$, this implies $j = 23$. Then ${\textrm{lk}}(6)$ = $C_{14}([{\textbf{7}}$, 0, 1, 2, 3, 4, ${\textbf{5}}$], 23, 22, 18, ${\textbf{17}}$, 16, 20, 8), ${\textrm{lk}}(17)$ = $C_{14}([{\textbf{16}}$, 15, 10, 11, 12, 19, ${\textbf{18}}$], 22, 23, 5, ${\textbf{6}}$, 7, 8, 20). Now completing successively we get ${\textrm{lk}}(1)$ = $C_{14}([{\textbf{0}}$, 7, 6, 5, 4, 3, ${\textbf{2}}$], 19, 18, 22, ${\textbf{14}}$, 13, 12, 11), ${\textrm{lk}}(14)$ = $C_{14}([{\textbf{13}}$, 21, 20, 8, 9, 23, ${\textbf{22}}$], 18, 19, 2, ${\textbf{1}}$, 0, 11, 12), ${\textrm{lk}}(2)$ = $C_{14}([{\textbf{3}}$, 4, 5, 6, 7, 0, ${\textbf{1}}$], 14, 22, 18, ${\textbf{19}}$, 12, 13, 21), ${\textrm{lk}}(19)$ = $C_{14}([{\textbf{12}}$, 11, 10, 15, 16, 17, ${\textbf{18}}$], 22, 14, 1, ${\textbf{2}}$, 3, 21, 13), ${\textrm{lk}}(3)$ = $C_{14}([{\textbf{4}}$, 5, 6, 7, 0, 1, ${\textbf{2}}$], 19, 12, 13, ${\textbf{21}}$, 20, 16, 15), ${\textrm{lk}}(21)$ = $C_{14}([{\textbf{13}}$, 14, 22, 23, 9, 8, ${\textbf{20}}$], 16, 15, 4, ${\textbf{3}}$, 2, 19, 12), ${\textrm{lk}}(4)$ = $C_{14}([{\textbf{5}}$, 6, 7, 0, 1, 2, ${\textbf{3}}$], 21, 20, 16, ${\textbf{15}}$, 10, 9, 23), ${\textrm{lk}}(15)$ $= C_{14}([{\textbf{16}}$, 17, 18, 19, 12, 11, ${\textbf{10}}$], 9, 23, 5, ${\textbf{4}}$, 3, 21, 20), ${\textrm{lk}}(5)$ = $C_{14}([{\textbf{4}}$, 3, 2, 1, 0, 7, ${\textbf{6}}$], 17, 18, 22, ${\textbf{23}}$, 9, 10, 15), ${\textrm{lk}}(23)$ = $C_{14}([{\textbf{9}}$, 8, 20, 21, 13, 14, ${\textbf{22}}$], 18, 17, 6, ${\textbf{5}}$, 4, 15, 10), ${\textrm{lk}}(9)$ = $C_{14}([{\textbf{8}}$, 20, 21, 13, 14, 22, ${\textbf{23}}$], 5, 4, 15, ${\textbf{10}}$, 11, 0, 7), ${\textrm{lk}}(10)$ = $C_{14}([{\textbf{15}}$, 16, 17, 18, 19, 12, ${\textbf{11}}$], 0, 7, 8, ${\textbf{9}}$, 23, 5, 4), ${\textrm{lk}}(12)$ = $C_{14}([{\textbf{11}}$, 10, 15, 16, 17, 18, ${\textbf{19}}$], 2, 3, 21, ${\textbf{13}}$, 14, 1, 0), ${\textrm{lk}}(13)$ = $C_{14}([{\textbf{14}}$, 22, 23, 9, 8, 20, ${\textbf{21}}$], 3, 2, 19, ${\textbf{12}}$, 11, 0, 1), ${\textrm{lk}}(16)$ = $C_{14}([{\textbf{15}}$, 10, 11, 12, 19, 18, ${\textbf{17}}$], 6, 7, 8, ${\textbf{20}}$, 21, 3, 4), ${\textrm{lk}}(18)$ $= C_{14}([{\textbf{17}}$, 16, 15, 10, 11, 12, ${\textbf{19}}$], 2, 1, 14, ${\textbf{22}}$, 23, 5, 6), ${\textrm{lk}}(20)$ = $C_{14}([{\textbf{8}}$, 9, 23, 22, 14, 13, ${\textbf{21}}$], 3, 4, 15, ${\textbf{16}}$, 17, 6, 7), ${\textrm{lk}}(22)$ = $C_{14}([{\textbf{23}}$, 9, 8, 20, 21, 13, ${\textbf{14}}$],1, 2, 19, ${\textbf{18}}$, 17, 6, 5). This is isomorphic to $N_2(6^2, 8)$ by the map (0, 3, 6, 1, 4, 7, 2, 5)(8, 15)(9, 10)(11, 23)(12, 22)(13, 18)(14, 19, 21, 17)(16, 20).

\noindent{\bf Subcase 3.2.1.4\,:} If $(g, h) = (18, 19)$ then completing ${\textrm{lk}}(6)$, ${\textrm{lk}}(5)$, ${\textrm{lk}}(21)$, ${\textrm{lk}}(20)$ we get ${\textrm{lk}}(4)$ = $C_{14}([{\textbf{3}}$, 2, 1, 0, 7, 6, ${\textbf{5}}$], 21, 20, 18, ${\textbf{17}}$, 16, $i, j)$ for some $i, j \in V(O_3)$. Observe that $i = 22$, this implies $j = 23$. Then ${\textrm{lk}}(17)$ = $C_{14}([{\textbf{16}}$, 15, 10, 11, 12, 19, ${\textbf{18}}$], 20, 21, 5, ${\textbf{4}}$, 3, 23, 22). Now completing successively we get ${\textrm{lk}}(1)$ = $C_{14}([{\textbf{0}}$, 7, 6, 5, 4, 3, ${\textbf{2}}$], 15, 16, 22, ${\textbf{14}}$, 13, 12, 11), ${\textrm{lk}}(14)$ $= C_{14}([{\textbf{13}}$, 21, 20, 8, 9, 23, ${\textbf{22}}$], 16, 15, 2, ${\textbf{1}}$, 0, 11, 12), ${\textrm{lk}}(2)$ = $C_{14}([{\textbf{3}}$, 4, 5, 6, 7, 0, ${\textbf{1}}$], 14, 22, 16, ${\textbf{15}}$, 10, 9, 23), ${\textrm{lk}}(15)$ = $C_{14}([{\textbf{16}}$, 17, 18, 19, 12, 11, ${\textbf{10}}$], 9, 23, 3, ${\textbf{2}}$, 1, 14, 22), ${\textrm{lk}}(3)$ = $C_{14}([{\textbf{4}}$, 5, 6, 7, 0, 1, ${\textbf{2}}$], 15, 10, 9, ${\textbf{23}}$, 22, 16, 17), ${\textrm{lk}}(23)$ = $C_{14}([{\textbf{9}}$, 8, 20, 21, 13, 14, ${\textbf{22}}$], 16, 17, 4, ${\textbf{3}}$, 2, 15, 10), ${\textrm{lk}}(7)$ = $C_{14}([{\textbf{0}}$, 1, 2, 3, 4, 5, ${\textbf{6}}$], 19, 18, 20, ${\textbf{8}}$, 9, 10, 11), ${\textrm{lk}}(9)$ = $C_{14}([{\textbf{8}}$, 20, 21, 13, 14, 22, ${\textbf{23}}$], 3, 2, 15, ${\textbf{10}}$, 11, 0, 7), ${\textrm{lk}}(10)$ = $C_{14}([{\textbf{15}}$, 16, 17, 18, 19, 12, ${\textbf{11}}$], 0, 7, 8, ${\textbf{9}}$, 23, 3, 2), ${\textrm{lk}}(12)$ = $C_{14}([{\textbf{11}}$, 10, 15, 16, 17, 18, ${\textbf{19}}$], 6, 5, 21, ${\textbf{13}}$, 14, 1, 0), ${\textrm{lk}}(13)$ = $C_{14}([{\textbf{14}}$, 22, 23, 9, 8, 20, ${\textbf{21}}$], 5, 6, 19, ${\textbf{12}}$, 11, 0, 1), ${\textrm{lk}}(16)$ = $C_{14}([{\textbf{15}}$, 10, 11, 12, 19, 18, ${\textbf{17}}$], 4, 3, 23, ${\textbf{22}}$, 14, 1, 2), ${\textrm{lk}}(18)$ = $C_{14}([{\textbf{17}}$, 16, 15, 10, 11, 12, ${\textbf{19}}$], 6, 7, 8, ${\textbf{20}}$, 21, 5, 4), ${\textrm{lk}}(19)$ = $C_{14}([{\textbf{12}}$, 11, 10, 15, 16, 17, ${\textbf{18}}$], 20, 8, 7, ${\textbf{6}}$, 5, 21, 13), ${\textrm{lk}}$ $(22)$ = $C_{14}([{\textbf{23}}$, 9, 8, 20, 21, 13, ${\textbf{14}}$], 1, 2, 15, ${\textbf{16}}$, 17, 4, 3). This is isomorphic to $N_1(6^2, 8)$ by the map (0, 20)(1, 21, 7, 14, 5, 8)(2, 22, 4, 9)(3, 23)(6, 13)(10, 17)(11, 18)(12, 19)(15, 16).

\smallskip

\noindent{\bf Subcase 3.2.2\,:} If $c = 22$ then we have $d \in \{$13, 14, 23$\}$.

If $d = 13$ then successively we get $e = 14$, $f = 23$ and $(g, h) \in \{$(16, 15), (16, 17), (17, 16), (17, 18), (18, 17)$\}$. If $(g, h) \in \{$(17, 16), (17, 18)$\}$ then considering ${\textrm{lk}}(8)$ we see that ${\textrm{lk}}(6)$ can not be completed. If $(g, h) = (16, 15)$ then successively considering ${\textrm{lk}}(8)$, ${\textrm{lk}}(6)$ and ${\textrm{lk}}(9)$ we see that ${\textrm{lk}}(14)$ and ${\textrm{lk}}(23)$ can not be completed. If $(g, h) = (16, 17)$ then ${\textrm{lk}}(8)$ = $C_{14}([{\textbf{9}}$, 23, 14, 13, 22, 21, ${\textbf{20}}$], 16, 17, 6, ${\textbf{7}}$, 0, 11, 10), ${\textrm{lk}}(6)$ = $C_{14}([{\textbf{5}}$, 4, 3, 2, 1, 0, ${\textbf{7}}$], 8, 20, 16, ${\textbf{17}}$, 18, $i, j)$ for some $i, j \in V(O_3)$. Observe that $(i, j)\in\{$(21, 22), (22, 21)$\}$. If $(i, j) = (21, 22)$ then successively considering ${\textrm{lk}}(6)$, ${\textrm{lk}}(5)$ and ${\textrm{lk}}(4)$ we see that $\deg(20) > 3$. A contradiction. If $(i, j) = (22, 21)$ then successively we get ${\textrm{lk}}(6)$ = $C_{14}([{\textbf{5}}$, 4, 3, 2, 1, 0, ${\textbf{7}}$], 8, 20, 16, ${\textbf{17}}$, 18, 22, 21), ${\textrm{lk}}(5)$ = $C_{14}([{\textbf{4}}$, 3, 2, 1, 0, 7, ${\textbf{6}}$], 17, 18, 22, ${\textbf{21}}$, 20, 16, 15), ${\textrm{lk}}(4)$ = $C_{14}([{\textbf{5}}$, 6, 7, 0, 1, 2, ${\textbf{3}}$], 23, 9, 10, ${\textbf{15}}$, 16, 20, 21) and ${\textrm{lk}}(9)$ = $C_{14}([{\textbf{8}}$, 20, 21, 22, 13, 14, ${\textbf{23}}$], 3, 4, 15, ${\textbf{10}}$, 11, 0, 7). This implies $C_9(0, 1, 14, 23, 3, 4, 5, 6, 7) \subseteq {\textrm{lk}}(2)$, a contradiction.
If $(g, h) = (18, 17)$ then ${\textrm{lk}}(8)$ = $C_{14}([{\textbf{9}}$, 23, 14, 13, 22, 21, ${\textbf{20}}$], 18, 17, 6, ${\textbf{7}}$, 0, 11, 10), ${\textrm{lk}}(6)$ = $C_{14}([{\textbf{5}}$, 4, 3, 2, 1, 0, ${\textbf{7}}$], 8, 20, 18, ${\textbf{17}}$, 16, $i, j)$ for some $i, j \in V(O_3)$. In this case, $(i, j)\in\{$(21, 22), (22, 21)$\}$. Now proceeding further we get a contradiction for each value of $(i, j)$. So $d \neq 13$.

If $d = 14$ then $e = 13$, $f = 23$ and $(g, h) \in \{$(16, 15), (16, 17), (17, 16), (17, 18), (18, 17), (18, 19)$\}$. If $(g, h) \in \{$(17, 16), (17, 18)$\}$ then considering ${\textrm{lk}}(8)$ we see ${\textrm{lk}}(15)$ or ${\textrm{lk}}(19)$ can not be completed. For the remaining values of $(g, h)$ we have following subcases.

\noindent{\bf Subcase A:} If $(g, h) = (16, 15)$ then successively we get ${\textrm{lk}}(8)$ = $C_{14}([{\textbf{9}}$, 23, 13, 14, 22, 21, ${\textbf{20}}$], 16, 15, 6, ${\textbf{7}}$, 0, 11, 10), ${\textrm{lk}}(7)$ = $C_{14}([{\textbf{0}}$, 1, 2, 3, 4, 5, ${\textbf{6}}$], 15, 16, 20, ${\textbf{8}}$, 9, 10, 11), ${\textrm{lk}}(6)$ $= C_{14}([{\textbf{5}}$, 4, 3, 2, 1, 0, ${\textbf{7}}$], 8, 20, 16, ${\textbf{15}}$, 10, 9, 23), ${\textrm{lk}}(15)$ = $C_{14}([{\textbf{10}}$, 11, 12, 19, 18, 17, ${\textbf{16}}$], 20, 8, 7, ${\textbf{6}}$, 5, 23, 9), ${\textrm{lk}}(23)$ = $C_{14}([{\textbf{9}}$, 8, 20, 21, 22, 14, ${\textbf{13}}$], 12, 19, 4, ${\textbf{5}}$, 6, 15, 10), ${\textrm{lk}}(5)$ = $C_{14}([{\textbf{4}}$, 3, 2, 1, 0, 7, ${\textbf{6}}$], 15, 10, 9, ${\textbf{23}}$, 13, 12, 19), ${\textrm{lk}}(4)$ = $C_{14}([{\textbf{3}}$, 2, 1, 0, 7, 6, ${\textbf{5}}$], 23, 13, 12, ${\textbf{19}}$, 18, $i, j)$ for some $i, j \in V(O_3)$. Observe that $i \in \{$21, 22$\}$. If $i = 21$ then $j = 22$. This implies $C_9(0, 1, 14, 22, 3, 4, 5, 6, 7) \subseteq {\textrm{lk}}(2)$. A contradiction. So $i = 22$ then $j = 21$. This implies ${\textrm{lk}}(4)$ = $C_{14}([{\textbf{5}}$, 6, 7, 0, 1, 2, ${\textbf{3}}$], 21, 22, 18, ${\textbf{19}}$, 12, 13, 23), ${\textrm{lk}}(19)$ = $C_{14}([{\textbf{12}}$, 11, 10, 15, 16, 17, ${\textbf{18}}$], 22, 21, 3, ${\textbf{4}}$, 5, 23, 13). Now completing successively we get ${\textrm{lk}}(1)$ = $C_{14}([{\textbf{0}}$, 7, 6, 5, 4, 3, ${\textbf{2}}$], 17, 18, 22, ${\textbf{14}}$, 13, 12, 11), ${\textrm{lk}}(14)$ = $C_{14}([{\textbf{13}}$, 23, 9, 8, 20, 21, ${\textbf{22}}$], 18, 17, 2, ${\textbf{1}}$, 0, 11, 12), ${\textrm{lk}}(2)$ = $C_{14}([{\textbf{3}}$, 4, 5, 6, 7, 0, ${\textbf{1}}$], 14, 22, 18, ${\textbf{17}}$, 16, 20, 21), ${\textrm{lk}}(17)$ = $C_{14}([{\textbf{16}}$, 15, 10, 11, 12, 19, ${\textbf{18}}$], 22, 14, 1, ${\textbf{2}}$, 3, 21, 20), ${\textrm{lk}}(3)$ $= C_{14}([{\textbf{4}}$, 5, 6, 7, 0, 1, ${\textbf{2}}$], 17, 16, 20, ${\textbf{21}}$, 22, 18, 19), ${\textrm{lk}}(21)$ = $C_{14}([{\textbf{22}}$, 14, 13, 23, 9, 8, ${\textbf{20}}$], 16, 17, 2, ${\textbf{3}}$, 4, 19, 18), ${\textrm{lk}}(9)$ = $C_{14}([{\textbf{8}}$, 20, 21, 22, 14, 13, ${\textbf{23}}$], 5, 6, 15, ${\textbf{10}}$, 11, 0, 7), ${\textrm{lk}}(10)$ = $C_{14}([{\textbf{15}}$, 16, 17, 18, 19, 12, ${\textbf{11}}$], 0, 7, 8, ${\textbf{9}}$, 23, 5, 6), ${\textrm{lk}}(12)$ = $C_{14}([{\textbf{11}}$, 10, 15, 16, 17, 18, ${\textbf{19}}$], 4, 5, 23, ${\textbf{13}}$, 14, 1, 0), ${\textrm{lk}}(13)$ = $C_{14}([{\textbf{14}}$, 22, 21, 20, 8, 9, ${\textbf{23}}$], 5, 4, 19, ${\textbf{12}}$, 11, 0, 1), ${\textrm{lk}}(16)$ = $C_{14}([{\textbf{15}}$, 10, 11, 12, 19, 18, ${\textbf{17}}$], 2, 3, 21, ${\textbf{20}}$, 8, 7, 6), ${\textrm{lk}}(18)$ = $C_{14}([{\textbf{17}}$, 16, 15, 10, 11, 12, ${\textbf{19}}$], 4, 3, 21, ${\textbf{22}}$, 14, 1, 2), ${\textrm{lk}}(20)$ $= C_{14}([{\textbf{8}}$, 9, 23, 13, 14, 22, ${\textbf{21}}$], 3, 2, 17, ${\textbf{16}}$, 15, 6, 7), ${\textrm{lk}}(22)$ = $C_{14}([{\textbf{14}}$, 13, 23, 9, 8, 20, ${\textbf{21}}$], 3, 4, 19, ${\textbf{18}}$, 17, 2, 1). This is isomorphic to $N_1(6^2, 8)$ by the map (0, 1, 2, 3, 4, 5, 6, 7)(8, 11, 14, 17, 23, 19, 21, 15)(9, 12, 20, 10, 13, 18, 22, 16).

\noindent{\bf Subcase B\,:} If $(g, h) = (16, 17)$ then ${\textrm{lk}}(8)$ = $C_{14}([{\textbf{9}}$, 23, 13, 14, 22, 21, ${\textbf{20}}$], 16, 17, 6, ${\textbf{7}}$, 0, 11, 10), ${\textrm{lk}}(7)$ = $C_{14}([{\textbf{0}}$, 1, 2, 3, 4, 5, ${\textbf{6}}$], 17, 16, 20, ${\textbf{8}}$, 9, 10, 11), ${\textrm{lk}}(6)$ = $C_{14}([{\textbf{5}}$, 4, 3, 2, 1, 0, ${\textbf{7}}$], 8, 20, 16, ${\textbf{17}}$, 18, $i, j)$ for some $i, j \in V(O_3)$. In this case $i \in \{$21, 22$\}$. If $i = 21$, $j = 22$. Now considering ${\textrm{lk}}(21)$ we see 15\,19 as an edge and a non-edge both. If $i = 22$ then $j = 21$. This implies ${\textrm{lk}}(6)$ = $C_{14}$ $([{\textbf{7}}$, 0, 1, 2, 3, 4, ${\textbf{5}}$], 21, 22, 18, ${\textbf{17}}$, 16, 20, 8), ${\textrm{lk}}(17)$ = $C_{14}([{\textbf{16}}$, 15, 10, 11, 12, 19, ${\textbf{18}}$], 22, 21, 5, ${\textbf{6}}$, 7, 8, 20). Now completing successively we get ${\textrm{lk}}(1)$ = $C_{14}([{\textbf{0}}$, 7, 6, 5, 4, 3, ${\textbf{2}}$], 19, 18, 22, ${\textbf{14}}$, 13, 12, 11), ${\textrm{lk}}(14)$ = $C_{14}([{\textbf{13}}$, 23, 9, 8, 20, 21, ${\textbf{22}}$], 18, 19, 2, ${\textbf{1}}$, 0, 11, 12), ${\textrm{lk}}(2)$ = $C_{14}([{\textbf{3}}$, 4, 5, 6, 7, 0, ${\textbf{1}}$], 14, 22, 18, ${\textbf{19}}$, 12, 13, 23), ${\textrm{lk}}(19)$ = $C_{14}([{\textbf{12}}$, 11, 10, 15, 16, 17, ${\textbf{18}}$], 22, 14, 1, ${\textbf{2}}$, 3, 23, 13), ${\textrm{lk}}(3)$ = $C_{14}([{\textbf{4}}$, 5, 6, 7, 0, 1, ${\textbf{2}}$], 19, 12, 13, ${\textbf{23}}$, 9, 10, 15), ${\textrm{lk}}$ $(23)$ = $C_{14}([{\textbf{9}}$, 8, 20, 21, 22, 14, ${\textbf{13}}$], 12, 19, 2, ${\textbf{3}}$, 4, 15, 10), ${\textrm{lk}}(4)$ = $C_{14}([{\textbf{5}}$, 6, 7, 0, 1, 2, ${\textbf{3}}$], 23, 9, 10, ${\textbf{15}}$, 16, 20, 21), ${\textrm{lk}}(15)$ = $C_{14}([{\textbf{16}}$, 17, 18, 19, 12, 11, ${\textbf{10}}$], 9, 23, 3, ${\textbf{4}}$, 5, 21, 20), ${\textrm{lk}}(21)$ = $C_{14}$ $([{\textbf{22}}$, 14, 13, 23, 9, 8, ${\textbf{20}}$], 16, 15, 4, ${\textbf{5}}$, 6, 17, 18), ${\textrm{lk}}(5)$ = $C_{14}([{\textbf{4}}$, 3, 2, 1, 0, 7, ${\textbf{6}}$], 17, 18, 22, ${\textbf{21}}$, 20, 16, 15), ${\textrm{lk}}(9)$ = $C_{14}([{\textbf{8}}$, 20, 21, 22, 14, 13, ${\textbf{23}}$], 3, 4, 15, ${\textbf{10}}$, 11, 0, 7), ${\textrm{lk}}(10)$ = $C_{14}([{\textbf{15}}$, 16, 17, 18, 19, 12, ${\textbf{11}}$], 0, 7, 8, ${\textbf{9}}$, 23, 3, 4), ${\textrm{lk}}(12)$ = $C_{14}([{\textbf{11}}$, 10, 15, 16, 17, 18, ${\textbf{19}}$], 2, 3, 23, ${\textbf{13}}$, 14, 1, 0), ${\textrm{lk}}$ $(13)= C_{14}([{\textbf{14}}$, 22, 21, 20, 8, 9, ${\textbf{23}}$], 3, 2, 19, ${\textbf{12}}$, 11, 0, 1), ${\textrm{lk}}(16)$ = $C_{14}([{\textbf{15}}$, 10, 11, 12, 19, 18, ${\textbf{17}}$], 6, 7, 8, ${\textbf{20}}$, 21, 5, 4), ${\textrm{lk}}(18)$ = $C_{14}([{\textbf{17}}$, 16, 15, 10, 11, 12, ${\textbf{19}}$], 2, 1, 14, ${\textbf{22}}$, 21, 5, 6), ${\textrm{lk}}(20)$ = $C_{14}([{\textbf{8}}$, 9, 23, 13, 14, 22, ${\textbf{21}}$], 5, 4, 15, ${\textbf{16}}$, 17, 6, 7), ${\textrm{lk}}(22)$ = $C_{14}([{\textbf{14}}$, 13, 23, 9, 8, 20, ${\textbf{21}}$], 5, 6, 17, ${\textbf{18}}$, 19, 2, 1). This is isomorphic to $N_1(6^2, 8)$ by the map (0, 2, 4, 6)(1, 3, 5, 7)(8, 14, 23, 21)(9, 20, 13, 22)(10, 18)(11, 17)(12, 16)(15, 19).

\noindent{\bf Subcase C\,:} If $(g, h) = (18, 17)$ then successively we get ${\textrm{lk}}(6)$ = $C_{14}([{\textbf{5}}$, 4, 3, 2, 1, 0, ${\textbf{7}}$], 8, 20, 18, ${\textbf{17}}$, 16, $i, j)$ for some $i, j \in V(O_3)$. In this case $(i, j)\in\{$(21, 22), (22, 21)$\}$. If $(i, j) = (21, 22)$ then considering ${\textrm{lk}}(21)$ we see that 15\,19 is both an edge and a non-edge. So $(i, j) = (22, 21)$ then ${\textrm{lk}}(6)$ = $C_{14}([{\textbf{7}}$, 0, 1, 2, 3, 4, ${\textbf{5}}$], 21, 22, 16, ${\textbf{17}}$, 18, 20, 8), ${\textrm{lk}}(17)$ = $C_{14}([{\textbf{16}}$, 15, 10, 11, 12, 19, ${\textbf{18}}$], 20, 8, 7, ${\textbf{6}}$, 5, 21, 22). Now completing successively we get ${\textrm{lk}}(1)$ = $C_{14}([{\textbf{0}}$, 7, 6, 5, 4, 3, ${\textbf{2}}$], 15, 16, 22, ${\textbf{14}}$, 13, 12, 11), ${\textrm{lk}}(14)$ = $C_{14}([{\textbf{13}}$, 23, 9, 8, 20, 21, ${\textbf{22}}$], 16, 15, 2, ${\textbf{1}}$, 0, 11, 12), ${\textrm{lk}}(2)$ = $C_{14}([{\textbf{3}}$, 4, 5, 6, 7, 0, ${\textbf{1}}$], 14, 22, 16, ${\textbf{15}}$, 10, 9, 23), ${\textrm{lk}}(15)$ = $C_{14}([{\textbf{16}}$, 17, 18, 19, 12, 11, ${\textbf{10}}$], 9, 23, 3, ${\textbf{2}}$, 1, 14, 22), ${\textrm{lk}}(23)$ = $C_{14}([{\textbf{9}}$, 8, 20, 21, 22, 14, ${\textbf{13}}$], 12, 19, 4, ${\textbf{3}}$, 2, 15, 10), ${\textrm{lk}}(3)$ = $C_{14}([{\textbf{4}}$, 5, 6, 7, 0, 1, ${\textbf{2}}$], 15, 10, 9, ${\textbf{23}}$, 13, 12, 19), ${\textrm{lk}}(19)$ = $C_{14}([{\textbf{12}}$, 11, 10, 15, 16, 17, ${\textbf{18}}$],20, 21, 5, ${\textbf{4}}$, 3, 23, 13), ${\textrm{lk}}(4)$ = $C_{14}([{\textbf{5}}$, 6, 7, 0, 1, 2, ${\textbf{3}}$], 23, 13, 12, ${\textbf{19}}$, 18, 20, 21), ${\textrm{lk}}(5)$ = $C_{14}([{\textbf{4}}$, 3, 2, 1, 0, 7, ${\textbf{6}}$], 17, 16, 22, ${\textbf{21}}$, 20, 18, 19), ${\textrm{lk}}(21)$ = $C_{14}([{\textbf{22}}$, 14, 13, 23, 9, 8, ${\textbf{20}}$], 18, 19, 4, ${\textbf{5}}$, 6, 17, 16), ${\textrm{lk}}(7)$ = $C_{14}([{\textbf{0}}$, 1, 2, 3, 4, 5, ${\textbf{6}}$], 17, 18, 20, ${\textbf{8}}$, 9, 10, 11), ${\textrm{lk}}(9)$ = $C_{14}([{\textbf{8}}$, 20, 21, 22, 14, 13, ${\textbf{23}}$],3, 2, 15, ${\textbf{10}}$, 11, 0, 7), ${\textrm{lk}}(10)$ = $C_{14}([{\textbf{15}}$, 16, 17, 18, 19, 12, ${\textbf{11}}$], 0, 7, 8, ${\textbf{9}}$, 23, 3, 2), ${\textrm{lk}}(12)$ = $C_{14}$ $([{\textbf{11}}$, 10, 15, 16, 17, 18, ${\textbf{19}}$], 4, 3, 23, ${\textbf{13}}$, 14, 1, 0), ${\textrm{lk}}(13)$ = $C_{14}([{\textbf{14}}$, 22, 21, 20, 8, 9, ${\textbf{23}}$], 3, 4, 19, ${\textbf{12}}$, 11, 0, 1), ${\textrm{lk}}(16)$ = $C_{14}([{\textbf{15}}$, 10, 11, 12, 19, 18, ${\textbf{17}}$], 6, 5, 21, ${\textbf{22}}$, 14, 1, 2), ${\textrm{lk}}(18)$ = $C_{14}([{\textbf{17}}$, 16, 15, 10, 11, 12, ${\textbf{19}}$], 4, 5, 21, ${\textbf{20}}$, 8, 7, 6), ${\textrm{lk}}(20)$ = $C_{14}([{\textbf{8}}$, 9, 23, 13, 14, 22, ${\textbf{21}}$], 5, 4, 19, ${\textbf{18}}$, 17, 6, 7), ${\textrm{lk}}(22)$ = $C_{14}([{\textbf{14}}$, 13, 23, 9, 8, 20, ${\textbf{21}}$], 5, 6, 17, ${\textbf{16}}$, 15, 2, 1). This is isomorphic to $N_2(6^2, 8)$ by the map (0, 7, 6, 5, 4, 3, 2, 1)(8, 17, 21, 19, 23, 15, 14, 11)(9, 16, 13, 10, 20, 18, 22, 12).

\noindent{\bf Subcase D\,:} If $(g, h) = (18, 19)$ then successively we get ${\textrm{lk}}(8)$ = $C_{14}([{\textbf{9}}$, 23, 13, 14, 22, 21, ${\textbf{20}}$], 18, 19, 6, ${\textbf{7}}$, 0, 11, 10), ${\textrm{lk}}(7)$ = $C_{14}([{\textbf{0}}$, 1, 2, 3, 4, 5, ${\textbf{6}}$], 19, 18, 20, ${\textbf{8}}$, 9, 10, 11), ${\textrm{lk}}(19)$ = $C_{14}([{\textbf{12}}$, 11, 10, 15, 16, 17, ${\textbf{18}}$], 20, 8, 7, ${\textbf{6}}$, 5, 23, 13), ${\textrm{lk}}(6)$ = $C_{14}([{\textbf{5}}$, 4, 3, 2, 1, 0, ${\textbf{7}}$], 8, 20, 18, ${\textbf{19}}$, 12, 13, 23), ${\textrm{lk}}(5)$ = $C_{14}([{\textbf{4}}$, 3, 2, 1, 0, 7, ${\textbf{6}}$], 19, 12, 13, ${\textbf{23}}$, 9, 10, 15), ${\textrm{lk}}(23)$ = $C_{14}([{\textbf{9}}$, 8, 20, 21, 22, 14, ${\textbf{13}}$], 12, 19, 6, ${\textbf{5}}$, 4, 15, 10), ${\textrm{lk}}(4)$ = $C_{14}([{\textbf{3}}$, 2, 1, 0, 7, 6, ${\textbf{5}}$], 23, 9, 10, ${\textbf{15}}$, 16, $i, j)$ for some $i, j \in V(O_3)$. In this case we see, $(i, j) = (22, 21)$. Then ${\textrm{lk}}(4)$ = $C_{14}([{\textbf{5}}$, 6, 7, 0, 1, 2, ${\textbf{3}}$], 21, 22, 16, ${\textbf{15}}$, 10, 9, 23), ${\textrm{lk}}(15)$ = $C_{14}([{\textbf{16}}$, 17, 18, 19, 12, 11, ${\textbf{10}}$], 9, 23, 5, ${\textbf{4}}$, 3, 21, 22), completing successively we get ${\textrm{lk}}(1)$ = $C_{14}([{\textbf{0}}$, 7, 6, 5, 4, 3, ${\textbf{2}}$], 17, 16, 22, ${\textbf{14}}$, 13, 12, 11), ${\textrm{lk}}(14)$ = $C_{14}([{\textbf{13}}$, 23, 9, 8, 20, 21, ${\textbf{22}}$], 16, 17, 2, ${\textbf{1}}$, 0, 11, 12), ${\textrm{lk}}(17)$ = $C_{14}([{\textbf{16}}$, 15, 10, 11, 12, 19, ${\textbf{18}}$], 20, 21, 3, ${\textbf{2}}$, 1, 14, 22), ${\textrm{lk}}(2)$ = $C_{14}([{\textbf{3}}$, 4, 5, 6, 7, 0, ${\textbf{1}}$], 14, 22, 16, ${\textbf{17}}$, 18, 20, 21), ${\textrm{lk}}(21)$ = $C_{14}([{\textbf{22}}$, 14, 13, 23, 9, 8, ${\textbf{20}}$], 18, 17, 2, ${\textbf{3}}$, 4, 15, 16), ${\textrm{lk}}(3)$ = $C_{14}([{\textbf{4}}$, 5, 6, 7, 0, 1, ${\textbf{2}}$], 17, 18, 20, ${\textbf{21}}$, 22, 16, 15), ${\textrm{lk}}(9)$ = $C_{14}$ $([{\textbf{8}}$, 20, 21, 22, 14, 13, ${\textbf{23}}$], 5, 4, 15, ${\textbf{10}}$, 11, 0, 7), ${\textrm{lk}}(10)$ = $C_{14}([{\textbf{15}}$, 16, 17, 18, 19, 12, ${\textbf{11}}$], 0, 7, 8, ${\textbf{9}}$, 23, 5, 4), ${\textrm{lk}}(12)$ = $C_{14}([{\textbf{11}}$, 10, 15, 16, 17, 18, ${\textbf{19}}$], 6, 5, 23, ${\textbf{13}}$, 14, 1, 0), ${\textrm{lk}}(13)$ = $C_{14}([{\textbf{14}}$, 22, 21, 20, 8, 9, ${\textbf{23}}$], 5, 6, 19, ${\textbf{12}}$, 11, 0, 1), ${\textrm{lk}}(16)$ = $C_{14}([{\textbf{15}}$, 10, 11, 12, 19, 18, ${\textbf{17}}$], 2, 1, 14, ${\textbf{22}}$, 21, 3, 4), ${\textrm{lk}}(18)$ = $C_{14}([{\textbf{17}}$, 16, 15, 10, 11, 12, ${\textbf{19}}$], 6, 7, 8, ${\textbf{20}}$, 21, 3, 2), ${\textrm{lk}}(20)$ = $C_{14}([{\textbf{8}}$, 9, 23, 13, 14, 22, ${\textbf{21}}$], 3, 2, 17, ${\textbf{18}}$, 19, 6, 7), ${\textrm{lk}}(22)$ = $C_{14}([{\textbf{14}}$, 13, 23, 9, 8, 20, ${\textbf{21}}$], 3, 4, 15, ${\textbf{16}}$, 17, 2, 1). This is isomorphic to $N_2(6^2, 8)$ by the map (0, 2, 4, 6)(1, 3, 5, 7)(8, 14, 23)(9, 20, 13)(10, 16, 18, 12)(11, 15, 17, 19).

\noindent{\bf Subcase 3.2.2.2\,:} If $d = 23$ then $(e, f) = (13, 14)$. This implies ${\textrm{lk}}(14)$ = $C_{14}([{\textbf{9}}$, 8, 20, 21, 22, 23, ${\textbf{13}}$], 12, 11, 0, ${\textbf 1}$, 2, 15, 10), ${\textrm{lk}}(1)$ = $C_{14}([{\textbf{0}}$, 7, 6, 5, 4, 3, ${\textbf{2}}$], 15, 10, 9, ${\textbf{14}}$, 13, 12, 11), ${\textrm{lk}}(9)$ $= C_{14}([{\textbf{8}}$, 20, 21, 22, 23, 13, ${\textbf{14}}$], 1, 2, 15, ${\textbf{10}}$, 11, 0, 7), ${\textrm{lk}}(10)$ = $C_{14}([{\textbf{11}}$, 12, 19, 18, 17, 16, ${\textbf{15}}$], 2, 1, 14, ${\textbf{9}}$, 8, 7, 0) and ${\textrm{lk}}(8)$ = $C_{14}([{\textbf{20}}$, 21, 22, 23, 13, 14, ${\textbf{9}}$], 10, 11, 0, ${\textbf{7}}$, 6, $h, g)$, where $(g, h) \in \{$(16, 17), (17, 16), (17, 18), (18, 17), (18, 19)$\}$. If $(g, h) = (17, 16)$ or $(17, 18)$ then considering ${\textrm{lk}}(8)$ we see that ${\textrm{lk}}(15)$ or ${\textrm{lk}}(19)$ can not be completed. Also, $(16, 17)\cong(18, 19)$ by the map (0, 9)(1, 14)(2, 13)(3, 23)(4, 22)(5, 21)(6, 20)(7, 8)(10, 11)(12, 15)(16, 19)(17, 18). So we search the map for $(g, h) \in \{$(16, 17), (18, 17)$\}$.

\noindent{\bf Subcase A\,:} If $(g, h) = (16, 17)$ then ${\textrm{lk}}(8)$ = $C_{14}([{\textbf{9}}$, 14, 13, 23, 22, 21, ${\textbf{20}}$], 16, 17, 6, ${\textbf{7}}$, 0, 11, 10), ${\textrm{lk}}(7)$ = $C_{14}$ $([{\textbf{0}}$, 1, 2, 3, 4, 5, ${\textbf{6}}$], 17, 16, 20, ${\textbf{8}}$, 9, 10, 11) and ${\textrm{lk}}(6)$ = $C_{14}([{\textbf{5}}$, 4, 3, 2, 1, 0, ${\textbf{7}}$], 8, 20, 16, ${\textbf{17}}$, 18, $i, j)$ for some $i, j \in V(O_3)$. It is easy to see that $i = 22$ or 23. If $i = 23$ then $j = 22$. Now considering ${\textrm{lk}}(5)$ we see 34 as an edge and a non-edge both. So $i = 22$ then $j = 23$. This implies ${\textrm{lk}}(6)$ = $C_{14}([{\textbf{7}}$, 0, 1, 2, 3, 4, ${\textbf{5}}$], 23, 22, 18, ${\textbf{17}}$, 16, 20, 8), ${\textrm{lk}}(17)$ = $C_{14}([{\textbf{16}}$, 15, 10, 11, 12, 19, ${\textbf{18}}$], 22, 23, 5, ${\textbf{6}}$, 7, 8, 20), completing successively we get ${\textrm{lk}}(2)$ = $C_{14}([{\textbf{3}}$, 4, 5, 6, 7, 0, ${\textbf{1}}$], 14, 9, 10, ${\textbf{15}}$, 16, 20, 21), ${\textrm{lk}}(15)$ = $C_{14}([{\textbf{16}}$, 17, 18, 19, 12, 11, ${\textbf{10}}$], 9, 14, 1, ${\textbf{2}}$, 3, 21, 20), ${\textrm{lk}}(21)$ = $C_{14}([{\textbf{22}}$, 23, 13, 14, 9, 8, ${\textbf{20}}$], 16, 15, 2, ${\textbf{3}}$, 4, 19, 18), ${\textrm{lk}}(3)$ = $C_{14}([{\textbf{4}}$, 5, 6, 7, 0, 1, ${\textbf{2}}$], 15, 16, 20, ${\textbf{21}}$, 22, 18, 19), ${\textrm{lk}}(4)$ = $C_{14}([{\textbf{5}}$, 6, 7, 0, 1, 2, ${\textbf{3}}$], 21, 22, 18, ${\textbf{19}}$, 12, 13, 23), ${\textrm{lk}}(23)$ $= C_{14}([{\textbf{22}}$, 21, 20, 8, 9, 14, ${\textbf{13}}$], 12, 19, 4, ${\textbf{5}}$, 6, 17, 18), ${\textrm{lk}}(5)$ = $C_{14}([{\textbf{4}}$, 3, 2, 1, 0, 7, ${\textbf{6}}$], 17, 18, 22, ${\textbf{23}}$, 13, 12, 19), ${\textrm{lk}}(12)$ = $C_{14}([{\textbf{11}}$, 10, 15, 16, 17, 18, ${\textbf{19}}$], 4, 5, 23, ${\textbf{13}}$, 14, 1, 0), ${\textrm{lk}}(13)$ = $C_{14}([{\textbf{14}}$, 9, 8, 20, 21, 22, ${\textbf{23}}$], 5, 4, 19, ${\textbf{12}}$, 11, 0, 1), ${\textrm{lk}}(20)$ = $C_{14}([{\textbf{8}}$, 9, 14, 13, 23, 22, ${\textbf{21}}$], 3, 2, 15, ${\textbf{16}}$, 17, 6, 7), ${\textrm{lk}}(16)$ = $C_{14}([{\textbf{15}}$, 10, 11, 12, 19, 18, ${\textbf{17}}$], 6, 7, 8, ${\textbf{20}}$, 21, 3, 2), ${\textrm{lk}}(22)$ = $C_{14}([{\textbf{23}}$, 13, 14, 9, 8, 20, ${\textbf{21}}$], 3, 4, 19, ${\textbf{18}}$, 17, 6, 5), ${\textrm{lk}}(18)$ = $C_{14}([{\textbf{17}}$, 16, 15, 10, 11, 12, ${\textbf{19}}$], 4, 3, 21, ${\textbf{22}}$, 23, 5, 6), ${\textrm{lk}}(19)$ = $C_{14}([{\textbf{12}}$, 11, 10, 15, 16, 17, ${\textbf{18}}$], 22, 21, 3, ${\textbf{4}}$, 5, 23, 13). This is isomorphic to $N_1(6^2, 8)$ by the map (0, 4)(1, 3)(5, 7)(8, 21, 14, 23)(9, 22, 13)(10, 16, 18, 12)(15, 17, 19, 11).

\noindent {\bf Subcase B\,:} If $(g, h) = (18, 17)$ then ${\textrm{lk}}(8)$ = $C_{14}([{\textbf{9}}$, 14, 13, 23, 22, 21, ${\textbf{20}}$], 17, 18, 6, ${\textbf{7}}$, 0, 11, 10), ${\textrm{lk}}(7)$ = $C_{14}([{\textbf{0}}$, 1, 2, 3, 4, 5, ${\textbf{6}}$], 17, 18, 20, ${\textbf{8}}$, 9, 10, 11). This implies ${\textrm{lk}}(6)$ = $C_{14}([{\textbf{5}}$, 4, 3, 2, 1, 0, ${\textbf{7}}$], 8, 20, 18, ${\textbf{17}}$, 16, j, i), ${\textrm{lk}}(17)$ = $C_{14}([{\textbf{16}}$, 15, 10, 11, 12, 19, ${\textbf{18}}$], 20, 8, 7, ${\textbf{6}}$, 5, 21, 22) and ${\textrm{lk}}(16)$ = $C_{14}([{\textbf{15}}$, 10, 11, 12, 19, 18, ${\textbf{17}}$], 6, 5, $i$, ${\textit {\textbf j}}$, $k$, 3, 2) for some $i, j, k \in V(O_3)$. Then we see, $j \in\{$21, 22, 23$\}$. If $j = 21$ then either $i = 20$ or $k = 20$. But in both cases, $\deg(20) > 3$. A contradiction. If $j = 23$ then $i = 13$ or $k = 13$, again in both cases, we see $\deg(13) > 3$. If $j = 22$ then $i \in \{$21, 23$\}$.

If $i = 21$ then $k = 23$. This implies ${\textrm{lk}}(16)$ = $C_{14}([{\textbf{15}}$, 10, 11, 12, 19, 18, ${\textbf{17}}$], 6, 5, 21, ${\textbf{22}}$, 23, 3, 2), ${\textrm{lk}}(22)$ = $C_{14}([{\textbf{23}}$, 13, 14, 9, 8, 20, ${\textbf{21}}$], 5, 6, 17, ${\textbf{16}}$, 15, 2, 3), completing successively we get ${\textrm{lk}}(2)$ = $C_{14}([{\textbf{3}}$, 4, 5, 6, 7, 0, ${\textbf{1}}$], 14, 9, 10, ${\textbf{15}}$, 16, 22, 23), ${\textrm{lk}}(15)$ = $C_{14}([{\textbf{16}}$, 17, 18, 19, 12, 11, ${\textbf{10}}$], 9, 14, 1, ${\textbf{2}}$, 3, 23, 22),  ${\textrm{lk}}(23)$ = $C_{14}([{\textbf{22}}$, 21, 20, 8, 9, 14, ${\textbf{13}}$], 12, 19, 4, ${\textbf{3}}$, 2, 15, 16), ${\textrm{lk}}(3)$ = $C_{14}([{\textbf{4}}$, 5, 6, 7, 0, 1, ${\textbf{2}}$], 15, 16, 22, ${\textbf{23}}$, 13, 12, 19), ${\textrm{lk}}(19)$ = $C_{14}([{\textbf{12}}$, 11, 10, 15, 16, 17, ${\textbf{18}}$], 20, 21, 5, ${\textbf{4}}$, 3, 23, 13), ${\textrm{lk}}(4)$ = $C_{14}([{\textbf{5}}$, 6, 7, 0, 1, 2, ${\textbf{3}}$], 23, 13, 12, ${\textbf{19}}$, 18, 20, 21), ${\textrm{lk}}(21)$ = $C_{14}([{\textbf{22}}$, 23, 13, 14, 9, 8, ${\textbf{20}}$], 18, 19, 4, ${\textbf{5}}$, 6, 17, 16), ${\textrm{lk}}(5)$ = $C_{14}([{\textbf{4}}$, 3, 2, 1, 0, 7, ${\textbf{6}}$], 17, 16, 22, ${\textbf{21}}$, 20, 18, 19), ${\textrm{lk}}$ $(12)$ = $C_{14}([{\textbf{11}}$, 10, 15, 16, 17, 18, ${\textbf{19}}$], 4, 3, 23, ${\textbf{13}}$, 14, 1, 0), ${\textrm{lk}}(13)$ = $C_{14}([{\textbf{14}}$, 9, 8, 20, 21, 22, ${\textbf{23}}$], 3, 4, 19, ${\textbf{12}}$, 11, 0, 1), ${\textrm{lk}}(18)$ = $C_{14}([{\textbf{17}}$, 16, 15, 10, 11, 12, ${\textbf{19}}$], 4, 5, 21, ${\textbf{20}}$, 8, 7, 6), ${\textrm{lk}}(20)$ = $C_{14}$ $([{\textbf{8}}$, 9, 14, 13, 23, 22, ${\textbf{21}}$], 5, 4, 19, ${\textbf{18}}$, 17, 6, 7). This is isomorphic to $N_1(6^2, 8)$ by the map (0, 13, 6, 20, 2, 9)(1, 8)(3, 23, 5, 21)(4, 22)(7, 14)(10, 11, 12, 19, 18, 17, 16, 15).

If $i = 23$ then $k = 21$. This implies ${\textrm{lk}}(16)$ = $C_{14}([{\textbf{15}}$, 10, 11, 12, 19, 18, ${\textbf{17}}$], 6, 5, 23, ${\textbf{22}}$, 21, 3, 2), ${\textrm{lk}}(22)$ = $C_{14}([{\textbf{23}}$, 13, 14, 9, 8, 20, ${\textbf{21}}$], 3, 2, 15, ${\textbf{16}}$, 17, 6, 5), completing successively we get ${\textrm{lk}}(2)$ = $C_{14}([{\textbf{3}}$, 4, 5, 6, 7, 0, ${\textbf{1}}$], 14, 9, 10, ${\textbf{15}}$, 16, 22, 21), ${\textrm{lk}}(15)$ = $C_{14}([{\textbf{16}}$, 17, 18, 19, 12, 11, ${\textbf{10}}$], 9, 14, 1, ${\textbf{2}}$, 3, 21, 22),  ${\textrm{lk}}(21)$ = $C_{14}([{\textbf{22}}$, 23, 13, 14, 9, 8, ${\textbf{20}}$], 18, 19, 4, ${\textbf{3}}$, 2, 15, 16), ${\textrm{lk}}(3)$ = $C_{14}([{\textbf{4}}$, 5, 6, 7, 0, 1, ${\textbf{2}}$], 15, 16, 22, ${\textbf{21}}$, 20, 18, 19), ${\textrm{lk}}(19)$ = $C_{14}([{\textbf{12}}$, 11, 10, 15, 16, 17, ${\textbf{18}}$], 20, 21, 3, ${\textbf{4}}$, 5, 23, 13$ $), ${\textrm{lk}}(4)$ = $C_{14}([{\textbf{5}}$, 6, 7, 0, 1, 2, ${\textbf{3}}$], 21, 20, 18, ${\textbf{19}}$, 12, 13, 23), ${\textrm{lk}}(23)$ = $C_{14}([{\textbf{22}}$, 21, 20, 8, 9, 14, ${\textbf{13}}$], 12, 19, 4, ${\textbf{5}}$, 6, 17, 16), ${\textrm{lk}}(5)$ = $C_{14}([{\textbf{4}}$, 3, 2, 1, 0, 7, ${\textbf{6}}$], 17, 16, 22, ${\textbf{23}}$, 13, 12, 19), ${\textrm{lk}}$ $(12)$ = $C_{14}([{\textbf{11}}$, 10, 15, 16, 17, 18, ${\textbf{19}}$], 4, 5, 23, ${\textbf{13}}$, 14, 1, 0), ${\textrm{lk}}(13)$ = $C_{14}([{\textbf{14}}$, 9, 8, 20, 21, 22, ${\textbf{23}}$], 5, 4, 19, ${\textbf{12}}$, 11, 0, 1), ${\textrm{lk}}(18)$ = $C_{14}([{\textbf{17}}$, 16, 15, 10, 11, 12, ${\textbf{19}}$], 4, 3, 21, ${\textbf{20}}$, 8, 7, 6), ${\textrm{lk}}(20)$ = $C_{14}$ $([{\textbf{8}}$, 9, 14, 13, 23, 22, ${\textbf{21}}$], 3, 4, 19, ${\textbf{18}}$, 17, 6, 7). This is isomorphic to $N_2(6^2, 8)$ by the map (0, 5)(1, 4)(2, 3)(6, 7)(8, 17)(9, 18, 20, 16)(10, 22)(11, 21, 15, 23)(12, 13)(14, 19). Thus the Lemma \ref{lem3} is proved. $\hfill\Box$

\smallskip

Table below shows a list of semi-equivelar maps on the surface of Euler characteristic $-1$ obtained in this article and in \cite{upadhyay1}.


\begin{center}
{\textbf{Table : Semi-equivelar maps on the surface of Euler characteristic $-1$}}
 \end{center}

\begin{center}
\renewcommand{\arraystretch}{1}
\begin{tabular}{|p{.8cm}|p{2cm}|p{1.7cm}|p{2cm}|p{3cm}|}
\hline
S.No. & SEM-Type & Exist (Yes/No) & Transitive or Not & Number of SEMs

\\
\hline
1&$(3^5, 4)$& YES & No & 3 $(K_1, K_2, K_3)$
\\
\hline
2&$(3^4, 4^2)$& No &$-$&$-$
\\
\hline
3&$(3^4, 8)$& No &$-$&$-$
\\
\hline
4&$(3^2, 4, 3, 6)$& No &$-$&$-$
\\
\hline
5&$(3, 4^4)$& No &$-$&$-$
\\
\hline
6&$(3, 4, 8, 4)$& Yes & No & 2 ($K_1(3, 4, 8, 4)$, $K_2(3, 4, 8, 4)$)
\\
\hline
7&$(3, 6, 4, 6)$& No &$-$&$-$
\\
\hline
8&$(4^3, 6)$& No &$-$&$-$
\\
\hline
9&$(4, 6, 16)$& Yes & No & 2 ($M_1(4, 6, 16)$, $M_2(4, 6, 16)$)
\\
\hline
10&$(4, 8, 12)$& No &$-$&$-$
\\
\hline
11&$(6^2, 8)$& Yes & No & 2 ($N_1(6^2, 8)$, $N_2(6^2, 8)$)
\\
\hline

%
%
%
%
%
%
%
%
%
%
%
%
%
%
%
%
%

\end{tabular}
\end{center}
\section{Acknowledgement}

Work of first author is partially supported by SERB, DST grant No. SR/S4/MS:717/10 and second author is a CSIR SRF (Award No. 09/1023(0003)/2010-EMR-I). The work reported here forms part of the PhD thesis of the second author.


{\small

\begin{thebibliography}{99}

\bibitem{altbrehm(1986)} Altshuler, A.; Brehm, U.; {The weakly neighborly polyhedral maps on the 2-manifolds with chi -1} \emph{Discrete Comput. Geom.} {\bf 1986}, 1,  355-369.


\bibitem{brehm(2008)} Brehm, U.; K$\ddot{\mbox u}$hnel, W. {Equivelar maps on torus}\emph{Europ. J. Combi.} {\bf 2008}, 29,  1843-1861.
\bibitem{brehm(1997)} Brehm, U.; Schulte, E. \emph{Polyhedral Maps, Handbook of Discrete and Computational Geometry}; CRC Press: Boca Raton, 1997.
\bibitem{babai(1991)} Babai, L. {Vertex transitive graphs and vertex transitive maps}, \emph{J. Graph Th.} {\bf 1991}, 15, 587-627.

\bibitem{dattau(2005)} Datta, B.; Upadhyay, A. K., Degree regular triangulation of torus and Klein bottle \emph{Proc. Indian Acad. Sci.} {\bf 2005}, 115, 279-307.

\bibitem{dattau(2006)} Datta, B.; Upadhyay, A. K., Degree regular triangulation of double torus \emph{Forum Math.} {\bf 2006}, 18, 1011-1025.
\bibitem{coxeter(1980)} Coxeter, H. S. M.; Moser, W. O. J.; \emph{Generators and Relations for Discrete Groups}; Springer-Verlag: Berlin, 1980.

\bibitem{conder(1995)} Conder, M. D. E.; Everitt B.; {Regular maps on non-orientable surfaces} \emph{Geom. Dedicata} {\bf 1995}, 56,  209-219.










\bibitem{conder(2009)} Conder, M. D. E.; {Regular maps and hyper maps of chi -1 to -200}, \emph{J. Comb. Theory} {\bf 2009}, 99,  455-459.


\bibitem{datta(2005)} Datta, B.; {A note on existence of {k, k}-equivelar polyhedral maps} \emph{Beitr. Algebra Geom.} {\bf 2005}, 46,  537-544.

\bibitem{grunbaum(1987)} Gr$\ddot{\mbox u}$nbaum, B.; Shephard G. C. \emph{Tilings and Patterns}; W. H. Freeman and com.: New York, 1987.
\bibitem{jones(1978)} Jones, G. A.; Singerman, D.; {Theory of maps on orientable surfaces}\emph{Proc. London Math. Soc.} {\bf 1978}, 37, 273-307.


\bibitem{karabas(2007)} Karabas, J.; Nedela, R.; {Archimedean solids of genus 2} \emph{Electronic Notes in Discrete Math.} {\bf 2007}, 28,  331-339.


\bibitem{karabas(2012)} Karabas, J.; Nedela R.; {Archimedean maps of higher genera} \emph{Mathematics of Comput.} {\bf 2012}, 81, 569-583.




\bibitem{karabas(web)} Karabas, \emph{Archimedean solids}, 2011, http://www.savbb.sk/$^{\thicksim}$ karabas/science.html.





\bibitem{datta_nandini} Datta, B.; Nilakantan, N., Equivelar polyhedral with few vertices \emph{Discrete Comput. Geom.} {\bf 2001}, 26,  429-461.




\bibitem{negami(2001)} Negami, S.; Nakamoto, A.; {Triangulations on closed surfaces covered by vertices of given degree} \emph{Graphs comb.} {\bf 2001}, 17,  529-537.


\bibitem{lutz(thesis)} Lutz, F. H. \emph{Triangulated Manifolds with Few Vertices and Vertex-Transitive Group Actions}; Shaker Verlag: Aachen, 1999.










\bibitem{frank(2010)} Lutz, F. H.; Sulanke, T.; Tiwari, A. K.; Upadhyay, A. K. \emph{arXiv:1001.2777[math.CO]} {\bf 2010}.











\bibitem{spanier(1981)} Spanier, E. H. \emph{Algebraic Topology}; Springer-Verlag: New York, 1981.



















%
%






%
%





%
%
%
%



\bibitem{bondy_murthy}
Bondy J. A. and Murthy U. S. R.,: {\em Graph Theory}, GTM 244, Springer, 2008.
\bibitem{upadhyay0}
Upadhyay A. K. and Tiwari A. K., An Enumeration of Semi-Equivelar Maps on Torus and Klein bottle(preprint)
\bibitem{upadhyay1}
Upadhyay A. K., Tiwari A. K. and Maity D.,: Semi-equivelar maps, {\em Beitr Algebra Geom}, DOI 10.1007/s13366-012-0130-6
\end{thebibliography}
}
\end{document}